\documentclass[a4paper,10pt]{article}
\usepackage[scaled]{helvet}              
\usepackage[T1]{fontenc}                 
\usepackage{geometry}  
\usepackage{setspace}                 
\usepackage{exscale}                  
\usepackage{a4wide}                   
\usepackage{framed}                   
\usepackage{subscript}                
\usepackage{authblk}                  
\usepackage{appendix}                 
\usepackage{colordvi}                 
\usepackage{enumitem}                 
\usepackage[american]{babel}
\usepackage{fancyhdr}                   
\usepackage[colorlinks=true,linkcolor=black,citecolor=black,filecolor=black,urlcolor=black,hypertexnames=false]{hyperref} 
\usepackage{aliascnt}      
\usepackage{cleveref}      
\usepackage[babel,german=guillemets]{csquotes}
\usepackage[backend=bibtex,style=numeric-comp,useprefix=true,hyperref=true,firstinits=true]{biblatex}   
\bibliography{../../../Import/References/references.bib}
\usepackage[fleqn]{amsmath} 
\usepackage{mathtools}      
\usepackage{amssymb}        
\usepackage{latexsym}       
\usepackage{dsfont}         
\numberwithin{equation}{section}      
\input{../../../Import/Environments/thm-environments.tex}
\input{../../../Import/Macros/mathmacros-global.tex}
\newcommand{\concNl}{number concentration}

\newcommand{\concRol}{mass concentration}

\newcommand{\massfluxRol}{mass flux}
\newcommand*{\tx}{(t,x)}
\newcommand*{\timeEnd}{T_0}
\newcommand*{\constBoltz}{\ensuremath{ k_b }}           
\newcommand*{\constCharge}{\ensuremath{ e_0 }}          
\newcommand*{\constEL}{\ensuremath{ \epsilon_0 }}       
\newcommand*{\constMG}{\ensuremath{ \mu_0 }}            
\newcommand{\entropy}{ s }                                                  
\newcommand{\entropyPure}{ s^{0} }                                          
\newcommand{\entropyMixl}[1][l]{ s^{mix}_{#1} }                             
\newcommand{\entropyFlux}{\ensuremath{ \vecj_{s} }}                         
\newcommand{\entropyPureFlux}{\ensuremath{ \vecj^{0}_{s} }}                 
\newcommand{\entropyMixFlux}{\ensuremath{ \vecj^{mix}_{s} }}                
\newcommand*{\diss}{\ensuremath{ \sigma }}                                  
\newcommand*{\dissPure}{\ensuremath{ \sigma^{0} }}                          
\newcommand*{\dissMix}{\ensuremath{ \sigma^{mix} }}                         
\newcommand*{\temp}{\ensuremath{ T }}                                       
\newcommand*{\volume}{\ensuremath{ V }}                                     
\newcommand*{\volumeSpecific}{\ensuremath{ v }}                             
\newcommand{\fieldF}{\ensuremath{ \vecu }}                                  
\newcommand*{\permeabH}{\ensuremath{ K }}                                   
\newcommand*{\chargeEL}[1][]{\ensuremath{ \rho_{f{#1}} }}                   
\newcommand*{\chargeELspecific}[1][]{\ensuremath{ \rho^{spec}_{f{#1}} }}    
\newcommand*{\currentEL}[1][]{ \ensuremath{ \veci_{f{#1}} } }               
\newcommand{\potEL}{\ensuremath{ \Phi }}                                    
\newcommand{\fieldEL}{\ensuremath{ \vecE }}                                 
\newcommand{\fieldMG}{\ensuremath{ \vecB }}                                 
\newcommand*{\ml}[1][l]{\ensuremath{ m_{{#1}} }}                                                    
\newcommand*{\Dl}[1][l]{\ensuremath{ d_{{#1}} }}                                                    
\newcommand*{\zl}[1][l]{\ensuremath{ z_{{#1}} }}                                                    
\newcommand*{\mobl}[1][l]{\ensuremath{ \omega_{{#1}} }}                                             
\newcommand{\energyTotl}[1][l]{ e^{tot}_{{#1}} }                                         
\newcommand{\energyTotlFlux}[1][]{\ensuremath{ \vecj_{e{#1}} }}                          
\newcommand{\energyIntl}[1][l]{ e^{int}_{{#1}} }                                         
\newcommand{\energyIntPurel}[1][l]{ e^{0}_{{#1}} }                                       
\newcommand{\energyIntPureAnsatz}{ \hat{e} }                                             
\newcommand{\energyIntMixl}[1][l]{ e^{mix}_{{#1}} }                                      
\newcommand{\stressTotl}[1][l]{ \vecT_{{#1}} }                                           
\newcommand{\stressViscl}[1][l]{ \vectau_{{#1}} }                                        
\newcommand{\pressHydrl}[1][l]{ p_{{#1}} }                                               
\newcommand{\pressTotl}[1][l]{ P_{{#1}} }                                                
\newcommand{\chempotl}[1][l]{ \mu_{{#1}} }                                               
\newcommand{\chempotPurel}[1][l]{ \mu^0_{{#1}} }                                         
\newcommand{\chempotMixl}[1][l]{ \mu^{mix}_{{#1}} }                                      
\newcommand{\elchempotl}[1][l]{ \mu^{el}_{{#1}} }                                        
\newcommand{\elchempotMixl}[1][l]{ \mu^{mix,el}_{{#1}} }                                 
\newcommand{\rl}[1][l]{\ensuremath{ r_{{#1}} }}                                          
\newcommand{\Rl}[1][l]{\ensuremath{ R_{{#1}} }}                                          
\newcommand{\nl}[1][l]{\ensuremath{ n_{{#1}} }}                             
\newcommand{\unitnl}{ \sqbrac{m^{-3}} }                                     
\newcommand{\yl}[1][l]{\ensuremath{ y_{{#1}} }}                              
\newcommand{\rol}[1][l]{\ensuremath{ \rho_{{#1}} }}                          
\newcommand{\rolfluxrel}[1][l]{\ensuremath{ \vecj_{{#1}} }}                  
\newcommand{\unitrol}{ \sqbrac{kg\;m^{-3}} }                                 
\newcommand{\unitrolFlux}{ \sqbrac{kg\;m^{-2}s^{-1}} }                       
\newcommand{\unitrolPDE}{ \sqbrac{kg\;m^{-3}s^{-1}} }                        
\newcommand{\slj}[1][l]{\ensuremath{ s_{{#1}j} }}                                        
\newcommand{\kfj}[1][j]{\ensuremath{ k^f_{{#1}} }}                                       
\newcommand{\kbj}[1][j]{\ensuremath{ k^b_{{#1}} }}                                       
\newcommand{\Kj}[1][j]{\ensuremath{ K^{{#1}} }}                                          
\newcommand{\Rj}[1][j]{\ensuremath{ R_{{#1}} }}                                          
\newcommand{\Rfjmal}[1][i]{\ensuremath{ \kfj \prod_{\slj[{#1}]<0} \yl[i]^{-\slj[{#1}]} }}   
\newcommand{\Rbjmal}[1][i]{\ensuremath{ \kbj \prod_{\slj[{#1}]>0} \yl[i]^{ \slj[{#1}]} }}   
\newcommand{\Rjmal}[1][i]{ \ensuremath{ \Rfjmal[i] - \Rbjmal[i] }}                       

\input{../../../Import/Environments/environments.tex}
\newcounter{countModAssump}  
\newcommand{\labelA}{A}

\let\sim~ 
\begin{document}
%
%
\title{A thermodynamically consistent model for multicomponent electrolyte solutions}
\author[1]{Matthias Herz} 
\author[1]{Peter Knabner}
\renewcommand\Affilfont{\itshape\small}
\affil[1]{Department of Mathematics, University of Erlangen-N\"urnberg, Cauerstr. 11, D-91058 Erlangen, Germany}
\date{\today}
\maketitle 
%
%
\begin{abstract}
This paper presents a thermodynamically consistent model for multicomponent electrolyte solutions. The first part of this paper derives the general governing equations for nonequilibrium systems 
within the theory of nonequilibrium thermodynamics. Here, we consider electrolyte solutions as general mixtures of charged constituents. Furthermore, in this part of the paper we combine the general 
theory of nonequilibrium thermodynamics with the well-known splittings of the entropy and the energy into a pure substance part and a part due to mixing. Thereby, we successfully establish evolution 
equations for both parts. Furthermore, we derive for both parts explicit expressions of the respective entropy production rates. Hence, we provide an approach that allows to study the entropy 
of mixing independently of the pure substance entropy and vice versa. This is of great value, in particular for a better understanding of the complex phenomena due to mixing in multicomponent systems. 
\par
In the second part of this paper, we close the system of general balance equations by applying constitutive laws. This is the crucial step in the modeling procedure. For this reason, we thermodynamically 
validate every involved constitutive law, i.e., we show that every constitutive law is in accordance with the second law of thermodynamics. Thus, the contribution of \cref{chapter:modelES} is to present 
thermodynamically consistent mathematical models for electrolyte solutions. Most importantly, the choices of the constitutive laws are motivated by the goal to obtain a model that contains the \pnp\ with 
convection, which is the classical and widely used mathematical model for electrolyte solutions. Hence, in \cref{chapter:modelES} we firstly provide for this classical model a thermodynamical verification, 
and secondly we clearly reveal the limitations of this model. Finally, by means of the general model for electrolyte solutions, we present a thermodynamically consistent extension of the \pnp.
\par
This paper is a revised and updated version of the previous preprint no.~375 with the same title.
\\[2.0mm]
\textbf{Keywords:} Nonequilibrium thermodynamics, mixture theory, thermodynamically consistent model, electrolyte solutions, electrohydrodynamics, \pnp. 
\end{abstract}
\fancyhead[L]{Introduction}
\section{Introduction}\label{sec:NonEqThermo-introduction}
This paper splits into two parts. In \cref{chapter:NonEqThermo} we consider general nonequilibrium systems, which are ubiquitous in biology, engineering, and hydrodynamics. 
Here, \enquote{nonequilibrium} means that not all macroscopic state variables are given constants in space and time. In fact, the characteristic feature of nonequilibrium systems are ongoing spatio-temporal dynamics, 
which lead to various kind of fluxes and production rates. Macroscopically, these spatio-temporal dynamics can be described by continuum mechanical densities, which vary in space and time.
However, nonequilibrium systems are not necessarily captured by continuum densities and commonly, the attribute that nonequilibrium systems possess locally well-defined densities is known 
as \textit{local thermodynamical equilibrium assumption} (LTE). Hence, the LTE assumption characterizes exactly those nonequilibrium systems, for which the macroscopic continuum mechanical description 
based on densities applies. In this paper, we confine ourselves to this class of nonequilibrium systems.
\par
Furthermore, to capture the ongoing spatio-temporal dynamics, the continuum mechanical description is based on the following general balance laws in differential form
\begin{align*}
 \dert a + \grad\cdot\vecj_a = r_a~,
\end{align*}
which naturally involve the flux~$\vecj_a$ and production rate~$r_a$ of the considered quantity~$a$. However, these fluxes and production rates are purely abstract balance quantities and, 
in particular, they contain no information about how the underlying processes, which drive the flux~$\vecj_a$ and the production rate~$r_a$, in reality take place. This information is provided 
in a second step by applying constitutive laws, i.e., by assuming that the fluxes and the production rates are given by a specific functional expression. 
Thus, the crucial modeling step is the choice of the constitutive laws, as this step transforms the abstract balance laws into 
real physical equations. Consequently, it is essential to ensure that the involved constitutive laws are in accordance with general physical principles. 
\medskip
\par
The main task of \textit{nonequilibrium thermodynamics}%
\footnote{Strictly speaking, there are different schools of nonequilibrium thermodynamics, e.g., see \cite{Lavenda-book} for a more detailed overview. 
          However, since all these schools are concerned with nonequilibrium processes, we henceforth subsume them under the common name \enquote{nonequilibrium thermodynamics}.}
is to provide criteria for validating constitutive laws. More precisely, nonequilibrium thermodynamics is based on the fact 
that energy is conserved (first law of thermodynamics), and that entropy never decreases (second law of thermodynamics). Usually these fundamental principles are formulated in terms of a 
vanishing energy production rate~$r_e=0$ and a nonnegative entropy production rate~$\diss\geq0$. This means that in nonequilibrium thermodynamics the first law and the second law are stated 
in differential from as balance equations for the total energy density~$\rol[]\energyTotl[]$ and the entropy density~$\rol[]\entropy$ 
\begin{align*}
 &\dert\brac{\rol[]\energyTotl[]} +\grad\cdot\energyTotlFlux = r_e \quad\text{ with } \quad r_e=0,        &~&\text{[first law of thermodynamics]}\\
 &\dert\brac{\rol[]\entropy} +\grad\cdot\entropyFlux = \diss       \quad\text{ with } \quad \diss\geq 0.  &~&\text{[second law of thermodynamics]}
\end{align*}
In a subsequent step, an ansatz for the total energy density~$\rol[]\energyTotl[]$ is chosen, and beginning from this ansatz, 
a specific formula for the entropy production rate~$\diss$ in terms of the remaining fluxes and production rates is derived in a long procedure. 
Thus, the goal of nonequilibrium thermodynamics is to establish a functional dependency in the form 
\begin{align*}
  \diss=f(\text{\enquote{fluxes}, \enquote{production rates}}),
\end{align*}
which yields the constraint $f(\text{\enquote{fluxes}, \enquote{production rates}})\geq0$ for the fluxes and production rates. This is exactly the minimal criterion, which must continue to hold true, 
when a constitutive law is applied to a flux or a production rate. Hence, nonequilibrium thermodynamics allows to validate constitutive laws in the sense that it can be shown, whether a constitutive law 
respects the second law of thermodynamics. Furthermore, the formula $\diss=f(\text{\enquote{fluxes}, \enquote{production rates}})$ shows which fluxes and production rates lead in which 
situations to a contribution such that we have $\diss>0$. Since $\diss>0$ characterizes irreversible processes, nonequilibrium thermodynamics allows to identify the irreversible subprocesses. 
An other important contribution of nonequilibrium thermodynamics is, that the construction of the general governing equations for a nonequilibrium system within the framework of nonequilibrium thermodynamics 
clearly reveals the involved (unrealistic) assumptions. This serves as crucial starting point for improving existing models.
\par
To summarize with the words of I. Prigogine in \cite[p. 336]{PrigogineKondepudi-book}: 
\enquote{\textit{nonequilibrium thermodynamics is founded on the explicit expression for $\diss$ in terms of the irreversible processes that we can identify and study experimentally.}}             
\medskip
\par
Historically, nonequilibrium thermodynamics was developed amongst many  others in 
\cite{AtkinCraine, Bowen1967, CompteJou, DeGrootMazur-book, Dreyer1987, Lavenda-book, Liu-book, Massoudi2007, Muller1968, MullerTomasso-book, PrigogineNicolis, PrigogineWiame, PrigogineKondepudi-book, RajagopalEtAll1986, Samohyl, TolmanFine, TruesdellToupin-book}. 
In particular the treatment of porous media was carried out amongst others in \cite{Bowen1980, Wilmanski-book, Bennethum1996-1,Bennethum1996-2} 
and \cite{Castellanos-book, DeGrootMazur-book, Masliyah-book, Samohyl, Roubicek2005-1, DreyerGuhlkeMuller} considered mixtures of charged constituents. As nonequilibrium thermodynamics typically deals with
multicomponent systems, i.e. with mixtures, nonequilibrium thermodynamics is in particular in the work of \cite{Samohyl, TruesdellToupin-book} referred to as mixture theory. 
\medskip
\par
This paper presents a derivation of the governing equations for nonequilibrium systems of \textit{charged constituents}, which are subject to LTE. The presentation is mainly taken 
from \cite{DeGrootMazur-book, PrigogineKondepudi-book, TruesdellToupin-book} and follows the fundamental principles~\ref{Modeling:step1}--\ref{Modeling:step3} below.
\par
The contribution of \cref{chapter:NonEqThermo} of this paper is to combine the well-known splitting of the internal energy in a pure substance part and a part due to mixing, cf. \cite{EckGarckeKnabner-book}, with the general procedure 
from \cite{DeGrootMazur-book}. Thereby, we derive evolution equations for the entropy of mixing and the pure substance entropy. The key point is that these equations lead to specific formulas for 
the respective entropy production rates, which allow to study each of this parts separately from each other. Furthermore, we obtain generalized versions of Dalton's law and Raoult's law. 
\medskip
\par
The rest of \cref{chapter:NonEqThermo} of this paper is organized as follows: Firstly, we list the assumptions and the fundamental modeling principles in \cref{sec:NonEqThermo-assumption}. 
Then, we proceed with the mass conservation equations in \cref{sec:NonEqThermo-mass}, the charge conservation equations in \cref{sec:NonEqThermo-charge}, the momentum conservation equations 
in \cref{sec:NonEqThermo-mom}, and the energy conservation equations in \cref{sec:NonEqThermo-energy}. Finally, we derive the evolution equation for the entropy density in \cref{sec:NonEqThermo-entropy}.
\medskip
\par
In \cref{chapter:modelES} of this paper, we apply the results of \cref{chapter:NonEqThermo} to electrolyte solutions. This is an important step, as the study of a particular electrolyte solution is at 
the heart of many applications in biology, engineering, and hydrodynamics. This task can be very challenging, as the characteristic feature of electrolyte solutions are various, simultaneously occurring 
physical phenomena. Moreover, these simultaneously occurring physical phenomena are usually mutually dependent. Consequently, in electrolyte solutions 
several coupling-induced nonlinearities arise, which still lead to new questions that are subject of present-day research. The reasons therefore are manifold:
Firstly, the coupled subprocesses lead to one overall observable resp.~measurable output. To experimentally detect from this overall output the informations 
about the respective subprocesses and their interplay is a challenging task. 
Secondly, it is even harder to control the interplay of coupled subprocesses, such that a desired output can be reliably produced. Nevertheless, this is essential for 
realizing technical applications.
Thirdly, when numerically simulating processes in electrolyte solutions, the algorithms have to cope with the coupling-induced nonlinearities. 
\medskip
\par
The basis for the before mentioned steps is a sound theoretical model, that adequately captures the characteristic features of the considered electrolyte solution. 
In \cref{chapter:modelES} of this paper, we show how to derive the governing equations for electrolyte solutions in the general framework of \textit{nonequilibrium thermodynamics}, 
which we presented in the preceding \cref{chapter:NonEqThermo}. By implementing this approach, we obtain a thermodynamical consistent model for electrolyte solutions. 
\par
More precisely, we adopt the general evolution equations from \cref{chapter:NonEqThermo}, and subsequently, we transform these abstract balance equations 
into specific physical equations by closing the resulting system of equations with the aid of constitutive laws. This is the crucial step in the modeling procedure 
and we account for this, as we thermodynamically validate every involved constitutive law. Here, thermodynamical validation means, we subsequently prove that every constitutive law 
is in accordance with the second law of thermodynamics. Altogether, the presented approach clearly reveals the construction of the resulting model.
\medskip
\par
For a detailed overview of nonequilibrium thermodynamics, we refer to the introduction of \cref{chapter:NonEqThermo}. In completion of this overview, we add that 
amongst others \cite{AtkinCraine2, BotheDreyer, Castellanos-book, CleggEtAl, DeGrootMazur-book, EringenIngram, HornJackson, Masliyah-book, Muller1968, OnsagerFuoss, Pitzer1, PitzerKim, Samohyl, Roubicek2005-1, DreyerGuhlkeMuller} 
contributed to mixtures of charged constituents resp. mixtures of chemically reacting constituents. 
\medskip 	
\par
The contribution of \cref{chapter:modelES} of this paper is to present a thermodynamically consistent model for multicomponent electrolyte solutions. In this sense, this paper 
continues the work of \cite{Castellanos-book, Masliyah-book, Samohyl, Roubicek2005-1, DreyerGuhlkeMuller}, where similar models for electrolyte solutions 
have been established before. However, in this paper we account for the fact, that the electric phenomena are governed by relativistic Maxwell's equations, 
while the remaining conservation laws are nonrelativistic equations, cf. \cite{ArnoldMechanics-book, LifshitzLandau-book1, LifshitzLandau-book2, LopezDavalos-book, MarsdenRatiu-book, TruesdellToupin-book}. 
For this reason, we derive in this paper the nonrelativistic limit of Maxwell's equations similar to \cite{Castellanos-book}. 
Moreover, we consider reactive electrolyte solutions and we model the reaction kinetics by means of the fundamental mass action law kinetics. Thus, 
the subsequent presentation includes a thermodynamical justification of the mass action law kinetics. Finally, the main contribution of \cref{chapter:modelES} of this paper is to show that the presented model can be 
reduced to the classical Poisson-Nernst--Planck model with convection. More precisely, we choose the constitutive laws exactly such that the resulting model contains the Poisson-Nernst--Planck 
model with convection. Hence,  we provide a thermodynamical verification of this classical model. Furthermore, we clearly reveal the involved assumptions of these classical models and by means 
of the general model from \cref{sec:modelES-model} we present a possible extension to more general situations.  
\medskip
\par  
The rest of \cref{chapter:modelES} of this paper is organized as follows: Firstly, in \cref{sec:modelES-maxwell}, we derive the electrostatic limit of Maxwell's equations. Then, in \cref{sec:modelES-pdes}, 
we summarize the general governing equations and in \cref{sec:modelES-energy}, we introduce the ansatzes for the internal energies. We close these system of equations in 
\cref{sec:modelES-reactions}--\cref{sec:modelES-heat}, by introducing and thermodynamically validating several constitutive laws. In \cref{sec:modelES-model}, 
we summarize the resulting mathematical model for electrolyte solutions and in \cref{sec:modelES-pnp} we show, that this model contains the famous and widely used Poisson-Nernst--Planck model.
%
%
%
%
\part{Nonequilibrium Thermodynamics}\label{chapter:NonEqThermo}
\fancyhead[L]{Assumptions}
\section{Assumptions and Fundamental Principles}\label{sec:NonEqThermo-assumption}
We now introduce the following assumptions, which we henceforth suppose to hold true  
\begin{enumerate}[label=({A}\arabic*), ref=({A}\arabic*), itemsep=0.0mm,  start=\value{countModAssump}+1]
 \item \textbf{Domain: } For $n\in\setN$, we henceforth consider a bounded domain~$\Omega\subset\setR^n$ with boundary~$\Gamma:=\partial\Omega$. 
       Furthermore, we suppose that this domain is a pure fluid domain, which is fully saturated with the considered mixture.%
       \label{Assump:Domain}
 \item \textbf{Time: } We assume that we observe the mixture over a certain time interval~$[0,\timeEnd]$.%
       \label{Assump:Time}
 \item \textbf{Temperature: } We assume thermal equilibrium inside the mixture. Hence, we have a unique temperature~$\temp$ for all constituents of the mixture.%
       \label{Assump:Temp}
 \item \textbf{Constituents: } For $L\in\setN$, we assume that the mixture consists of $L$~different constituents, which represent $L$ different chemical species. 
        For the chemical species, we use the index $l\in\{1,\ldots,L\}$. We suppose that we have one solvent, which is indexed such that this solvent is the $L$th chemical species.%
       \label{Assump:Constituents}
 \item \textbf{Charged constituents: } We allow for charged chemical species carrying the charges~$\constCharge\zl\sim[C]$. 
        Here, $\constCharge$ is the elementary charge and $\zl$ is the valency. Thus electrically neutral chemical species are included via $\zl=0$.%
       \label{Assump:ChargedConstituents}
 \item \textbf{Mass conservation in chemical reactions: } We assume that the sum of all mass production rates~$\rl$ vanishes, i.e., $\sum_l \rl = 0$.%
        \label{Assump:reactions}   
 \item  \textbf{Charge conservation in chemical reactions: } We suppose that mass production rates~$\rl$ are subject to $\sum_l \frac{\constCharge\zl}{\ml}\rl = 0.$%
        \label{Assump:reactionsAndCharges}
 \item  \textbf{Conservation of momentum: } We assume $\sum_l \sqbrac{\vecF^{int}_l +\rl\fieldF_l} =0$ for the internal interaction forces $\vecF^{int}_l$ 
        and the momentum transfer due to chemical reactions~$\rl\fieldF_l$.%
        \label{Assump:momentum} 
  \item  \textbf{Definition of the total pressures: } The partial total pressures~$\pressTotl$ are defined with partial stress tensors~$\stressTotl$ 
         by $\pressTotl:=-\frac{1}{n}\trace{\stressTotl}$. Analogous, we define the total pressure~$\pressTotl[]$ with the mixture stress tensor~$\stressTotl[]$
         by $\pressTotl[]:=-\frac{1}{n}\trace{\stressTotl[]}$.%
        \label{Assump:thermoPressure}     
 \item  \textbf{Decomposition of the stress tensors: } We assume that the stress tensors~$\stressTotl$ are given by%
        \footnote{We denote the unit matrix by $\mathds{1}\in\setR^{n\times n}$ and the trace of a matrix $A \in\setR^{n\times n}$ by $\trace{A}=\sum_i A_{ii}$. }
        $\stressTotl=-\pressHydrl\mathds{1}+\stressViscl$. For the viscous stress tensors~$\stressViscl$, we suppose symmetry, i.e., $\stressViscl=\stressViscl^\top$.%
        \label{Assump:stresses}   
 \item  \textbf{Electrostatics: } We assume that the electric phenomena are captured by an electric field~$\fieldEL$ and an electrostatic potential~$\potEL$, which are conneted via $\fieldEL=-\grad\potEL$.%
        \label{Assump:electrostatics}  
\setcounter{countModAssump}{\value{enumi}}
\end{enumerate}
\begin{remark}
  Assumption~\ref{Assump:Domain} means, that in the context of porous media, we are on the pore scale, looking inside a single pore. 
  See \cite{bear-book, Hornung-book, Muhammad-book} for further details and an introduction to the modeling of porous media.
\hfill$\square$
\end{remark}
\begin{remark}\label{remark:NonEqThermo-distinguishPressures}
  In this paper, we distinguish between the partial hydrostatic pressures~$\pressHydrl$ and the partial total pressures~$\pressTotl$. This is similar to, e.g., \cite{BotheDreyer}. 
  Note that assumption~\ref{Assump:stresses} and \ref{Assump:thermoPressure} reveal the connection 
  \begin{align*}
   \pressTotl = -\frac{1}{n}\trace{\stressTotl} = -\frac{1}{n}\trace{-\pressHydrl\mathds{1}+\stressViscl} = \pressHydrl - \frac{1}{n}\trace{\stressViscl}~. 
  \end{align*}
  Hence, the partial total pressures~$\pressTotl$ and the partial hydrostatic pressures~$\pressHydrl$ coincide, provided we have traceless partial viscous stress tensors~$\stressViscl$.
\hfill$\square$
\end{remark}
Henceforth, we derive the governing equations for mixtures of charged constituents based on the following fundamental principles, cf. \cite{TruesdellToupin-book}
\begin{enumerate}[label=(M\arabic*), ref=(M\arabic*), itemsep=0.0mm]
 \item Firstly, we postulate abstract conservation laws for each constituent of the mixture.%
       \label{Modeling:step1}
 \item Secondly, the conservation laws of the mixture as a whole are derived by summing over the corresponding conservation laws of the constituents. 
       This procedure reveals how the behavior of the mixture depends on the behavior of the constituents.%
       \label{Modeling:step2}
 \item Thirdly, the conservation equations for the mixture should have the same form as the corresponding conservation law of a single medium. This is ensured by defining the physical quantities of the mixture 
       such that, in the end, the conservation laws of the mixture look like single medium equations. Compared to the corresponding physical quantities of the constituents, this step leads to a generalized notion
       for some physical quantities of the mixture .%
       \label{Modeling:step3}
\end{enumerate} 
%
%
%
\fancyhead[L]{Mass conservation}
\section{Mass Conservation Equations}\label{sec:NonEqThermo-mass}
In this section, we briefly derive abstract mass conservation equations, which govern the kinetics on continuum scales. Here, we characterize continuum scales by simultaneously considering 
a large number of particles of a given chemical species. This approach leads to averaged kinetics, which are formulated in terms of the following quantities: 
\begin{enumerate}[align=left, label=(\roman*), leftmargin=*, topsep=2.0mm, itemsep=-1.2mm]
  \item In a representative elementary volume (REV)~$\volume$~$[m^{3}]$, we assume that~$N_l$ particles of the $l$th chemical species are present.
	To simultaneously track these particles, we define the \concNl~$\nl$ and the \concRol~$\rol$, cf.~\cite[Chaper~6]{Masliyah-book}, by
        \begin{align}\label{eq:NonEqThermo-mass-defConc}
           \nl:= N_l\volume^{-1} \quad\sim\unitnl \qquad\text{and}\qquad  \rol:= \ml\nl \quad\sim\unitrol.
        \end{align}
        Here, $\ml~[kg]$ are the molecular masses. Moreover, we henceforth identify the given chemical species with their concentrations~$\rol$ resp. $\nl$. 
        Next, we note that summing over all chemical species, defines the total \concRol
        \begin{align}\label{eq:NonEqThermo-mass-defToTConc}
         \rol[]:=\sum_l \rol \qquad\sim\unitrol~.
        \end{align} 
        Furthermore, we introduce the dimensionless mass fractions
        \begin{align}\label{eq:NonEqThermo-mass-defMassFrac}
           \yl:= \frac{\rol}{\rol[]}, \qquad\text{for which } \eqref{eq:NonEqThermo-mass-defToTConc} \text{ reads as } \sum_l\yl =1.
        \end{align} 
  \item To simultaneously describe the movement of a large number of molecules of the $l$th chemical species, we suppose that the \concRol~$\rol$ moves with the averaged velocity field~$\fieldF_l$. 
	Furthermore, in each REV the center of total mass moves with the so-called barycentric velocity field~$\fieldF$. This is the velocity, that is visible to an external observer at rest and it is defined by
        \begin{align*}
         \fieldF := \sum_l \yl \fieldF_l \qquad\sim\sqbrac{ms^{-1}}.
        \end{align*}
  \item The relative movement of the particles of the $l$th chemical species with respect to the barycentric velocity field~$\fieldF$ is described by the so-called drift velocity~$(\fieldF_l-\fieldF)$.
        Furthermore, with this drift velocity field we define the so-called drift \massfluxRol~$\rolfluxrel$ by
        \begin{align}\label{eq:NonEqThermo-mass-defDriftFlux}
          \rolfluxrel:= \rol(\fieldF_l-\fieldF) \qquad\sim\unitrolFlux.
        \end{align}       
        Thus, $\rolfluxrel$ describes the relative movement of the $l$th chemical species~$\rol$ with respect to the barycentric velocity field~$\fieldF$. From the definition of $\rolfluxrel$ follows that
        \begin{align}\label{eq:NonEqThermo-mass-sumCondition}
	  \sum_l \rolfluxrel = \sum_l \rol (\fieldF_l -\fieldF) = \rol[]\sum_l \yl\fieldF_l - \sum_l\rol\fieldF = \rol[]\fieldF - \rol[]\fieldF = 0.
	\end{align}
\end{enumerate}
Following \cite{EckGarckeKnabner-book, Gurtin-book, oden-book}, we formulate with the just defined quantities the mass continuity equations for 
the constituents of the mixture. More precisely, we consider a moving REV~$\volume(t)$ and for this REV we claim the general balance statement: 
\textit{The temporal change of total mass in $\volume(t)$ arises due to mass production.} Mathematically, this balance statement is for each constituent~$\rol$, $l\in\cbrac{1,\ldots,L}$, 
with the mass production rates~$\rl\sim[kg/(m^3s)]$ given by
\begin{align*}
 \derr \Intdx[\volume(t)]{ \rol ~} = \Intdx[\volume(t)]{ \rl ~}~.
\end{align*}
Application of Reynold's transport theorem, cf. \cite{EckGarckeKnabner-book}, yields the equivalent equation
\begin{align*}
 \Intdx[\volume(t)]{ \dert\rol + \grad\cdot\brac{\rol\fieldF_l} ~} = \Intdx[\volume(t)]{ \rl ~}~.
\end{align*}
We assume that for $l\in\cbrac{1,\ldots,L}$,  the involved quantities are sufficiently smooth, such that we can \enquote{drop} the integrals 
and formulate the preceding integral equations as pointwise differential equations, which are exactly the mass conservation equations for the constituents. 
These equations are given with \eqref{eq:NonEqThermo-mass-defDriftFlux} by
\begin{alignat}{2}\label{eq:NonEqThermo-mass-massBalance}
\dert\rol + \grad\cdot\brac{ \rol\fieldF + \rolfluxrel} ~=~\rl       && \qquad\sim\unitrolPDE ~ .
\end{alignat}
To establish the mass conservation equation for the total \concRol~$\rol[]$, we sum equations~\eqref{eq:NonEqThermo-mass-massBalance} over the index~$l$. Thereby, we obtain with \ref{Assump:reactions}, 
\eqref{eq:NonEqThermo-mass-defToTConc}, and \eqref{eq:NonEqThermo-mass-sumCondition}
\begin{align}\label{eq:NonEqThermo-mass-massBalanceTot}
	    & \dert \sum_l \rol + \grad\cdot\brac{ \sum_l \rol\fieldF + \sum_l\rolfluxrel } = \sum_l \rl  &&\sim\unitrolPDE, \nonumber\\[2.0mm]
\Equivalent~& \dert\rol[] + \grad\cdot\brac{ \rol[] \fieldF } = 0                                         &&\sim\unitrolPDE~.
\end{align}
Next, we derive the equations for the mass fractions~$\yl$ by modifying the left-hand side of equations~\eqref{eq:NonEqThermo-mass-massBalance} 
by using \eqref{eq:NonEqThermo-mass-massBalanceTot} and the material derivative~$\derm$ from~\eqref{eq:NonEqThermo-mom-matDer}.This yields
\begin{align*}
    \dert\rol + \grad\cdot\sqbrac{ \rol\fieldF + \rolfluxrel} 
  = \rol[]\brac{\dert\yl + \fieldF\cdot\grad\yl} +  \grad\cdot\rolfluxrel + \yl\sqbrac{\dert\rol[] + \grad\cdot\brac{ \rol[]\fieldF}} 
  = \rol[]\derm\yl +  \grad\cdot\rolfluxrel~. 
\end{align*}
Hence, for $l\in\cbrac{1,\ldots,L}$, the mass fractions~$\yl$ solve the following equations, which we equivalently write in conservative form and nonconservative form
\begin{subequations}
\begin{align}
 & \dert(\yl\rol[]) + \grad\cdot\brac{\yl\rol[]\fieldF + \rolfluxrel}  =  \rl  \qquad     &\text{in  }&\Omega \qquad\sim\unitrolPDE ~,   \label{eq:NonEqThermo-mass-massFracBalance-Cons}    \\[2.0mm]
 & \rol[]\derm\yl = - \grad\cdot\rolfluxrel + \rl  \qquad                                 &\text{in  }&\Omega \qquad\sim\unitrolPDE ~.   \label{eq:NonEqThermo-mass-massFracBalance-NonCons}
\end{align}
\end{subequations} 
Finally, we return to \eqref{eq:NonEqThermo-mass-massBalanceTot}, which we multiply by $-\rol[]^{-1}$. Thereby, we arrive at
\begin{align*}
  0  = -\rol[]^{-1}\dert\rol[] - \rol[]^{-1}\grad\cdot\brac{\rol[]\fieldF} 
     = \rol[]\dert\!\brac{\rol[]^{-1}} + \fieldF\cdot\grad\brac{\rol[]^{-1}}  - \grad\cdot\fieldF.
\end{align*}
Thus, with the material derivate~$\derm$ from \eqref{eq:NonEqThermo-mom-matDer}, the specific volume~$\volumeSpecific:=\rol[]^{-1}$ solves the equation
\begin{align}\label{eq:NonEqThermo-mass-specificVolPDE}
  & \rol[]\derm\volumeSpecific = \grad\cdot\fieldF &\text{in  }&\Omega \qquad\sim\sqbrac{s^{-1}}.
\end{align}
\begin{remark}[Independent drift \massfluxRol{es}]\label{remark:NonEqThermo-independentDriftFluxes}
  By rearranging \eqref{eq:NonEqThermo-mass-sumCondition}, we obtain 
  \begin{align}\label{eq:NonEqThermo-mass-sumCondition2}
    \rolfluxrel[L] = - \sum_{l=1}^{L-1} \rolfluxrel.
  \end{align}
  Hence, we have only $L-1$ independent drift \massfluxRol{es}~$\rolfluxrel$. This dependency reflects the physical fact, that moving solute molecules always collide with the solvent molecules. 
  Thus, among other reasons, the drift \massfluxRol~$\rolfluxrel[L]$ of the solvent is always caused by the drift \massfluxRol{es}~$\rolfluxrel[l]$, $l\in\cbrac{1,\ldots,L-1}$, of the solutes. 
  However, the preceding equation shows, that at the same time, this is the only reason. This means that the drift motion of the solutes completely determine the drift motion of the solvent.  
\hfill$\square$
\end{remark}
\begin{remark}[Independent variables]\label{remark:NonEqThermo-independentVariables}
  We derived equation~\eqref{eq:NonEqThermo-mass-massBalanceTot} by summing over equations~\eqref{eq:NonEqThermo-mass-massBalance}. Hence, equation~\eqref{eq:NonEqThermo-mass-massBalanceTot} 
  is a linear combination and does not provide new information. This shows, that we have the $L$~primal unknowns $\cbrac{\rol[1],\ldots,\rol[L]}$, which respectively solve the 
  $L$~equations~\eqref{eq:NonEqThermo-mass-massBalance}. Furthermore, these primal unknowns determine $\rol[]$ via \eqref{eq:NonEqThermo-mass-defToTConc}.
  \par
  On the other hand, inserting \eqref{eq:NonEqThermo-mass-sumCondition2} into the equation for the solvent~$\rol[L]$ may lead to a complicated equation. 
  Thus, it is often more convenient to drop the solvent~$\rol[L]$ and to consider the $L$~primal unknowns $\cbrac{\rol[1],\ldots,\rol[L-1],\rol[]}$, 
  which respectively solve equations~\eqref{eq:NonEqThermo-mass-massBalance} for $l\in\cbrac{1,\ldots,L-1}$ and equation~\eqref{eq:NonEqThermo-mass-massBalanceTot}. 
  The \concRol~$\rol[L]$ of the solvent is obtained from these choice of primal unknowns via  
  \begin{align*}
  \rol[L] = \rol[] -\sum_{l=1}^{L-1} \rol~. \\[-14.0mm]
  \end{align*} 
\hfill$\square$
\end{remark}
\fancyhead[L]{Charge conservation}
\section{Charge Conservation Equations}\label{sec:NonEqThermo-charge}
We already mentioned in the introduction, that we consider mixtures that consist of charged constituents. According to \ref{Assump:ChargedConstituents}, the respective charges of the constituents 
are given by $\constCharge\zl\sim[C]$. Thus, electrically neutral constituents are included via $\zl=0$.
\par
For any kind of charged or neutral chemical species, the evolution of the \concRol{s}~$\rol$ are governed by equations~\eqref{eq:NonEqThermo-mass-massBalance}. 
However, for charged chemical species transport of mass is equivalent to the transport of charges, i.e., with electric currents. This means, we have to care about both: mass transport 
and charge transport (or equivalently electric currents). To account for this, we multiply the mass transport equations~\eqref{eq:NonEqThermo-mass-massBalance} by the constants~$\frac{\constCharge\zl}{\ml}$. 
Thereby, we obtain for $l\in\cbrac{1,\ldots,L}$ with \eqref{eq:NonEqThermo-mass-defConc} the corresponding charge transport equations
\begin{align}\label{eq:NonEqThermo-charge-chargeBalance} 
 & \dert(\constCharge\zl\nl) + \grad\cdot\brac{\constCharge\zl\nl\fieldF + \frac{\constCharge\zl}{\ml}\rolfluxrel}  =  \frac{\constCharge\zl}{\ml}\rl  \qquad    &\text{in  }&\Omega \qquad\sim\sqbrac{Cm^{-3}s^{-1}}.   
\end{align} 
The physical interpretation of these equations is 
\begin{align*} 
  \underbrace{ \dert(\constCharge\zl\nl) }_{ \begin{subarray}{l} \text{temporal change of } \\ \text{the $l$th charge density} \end{subarray} } 
+ \underbrace{ \grad\cdot\brac{\constCharge\zl\nl\fieldF} }_{ \begin{subarray}{l} \text{electric current} \\ \text{due to barycentric} \\ \text{convection} \end{subarray} } 
+ \underbrace{ \grad\cdot\brac{\frac{\constCharge\zl}{\ml}\rolfluxrel} }_{ \begin{subarray}{l} \text{electric current} \\ \text{due to mixing} \\ \text{(=drift motion)} \end{subarray} }  
= \underbrace{ \frac{\constCharge\zl}{\ml}\rl~. }_{ \begin{subarray}{l} \text{charge production} \end{subarray} }   
\end{align*} 
Next, we sum over equations~\eqref{eq:NonEqThermo-charge-chargeBalance} and we define the free charge density~$\chargeEL\sim[Cm^{-3}]$ by 
\begin{align}\label{eq:NonEqThermo-charge-DefFreeCharge}
 \chargeEL:=\sum_l\frac{\constCharge\zl}{\ml}\rol=\sum_l\constCharge\zl\nl ~. 
\end{align}
Furthermore, we introduce the free current density~$\currentEL\sim[A/m^2]$ by
\begin{align}\label{eq:NonEqThermo-charge-DefFreeCurrent}
 \currentEL:=\sum_l \frac{\constCharge\zl}{\ml}\rolfluxrel~. 
\end{align}
This leads us with \ref{Assump:reactionsAndCharges} to the charge conservation equation
\begin{align}\label{eq:NonEqThermo-charge-chargeBalanceTot} 
             & \dert\sqbrac{\sum_l\constCharge\zl\nl} + \grad\cdot\sqbrac{\sum_l\constCharge\zl\nl~\fieldF + \sum_l\frac{\constCharge\zl}{\ml}\rolfluxrel}  =  \sum_l\frac{\constCharge\zl}{\ml}\rl   \qquad\text{in  }\Omega \nonumber\\[4.0mm]
 \Equivalent & \dert\chargeEL + \grad\cdot\brac{\chargeEL\fieldF + \currentEL}  =  0    \qquad\text{in  }\Omega \qquad\sim\sqbrac{Cm^{-3}s^{-1}}.   
\end{align} 
%
\begin{remark}[Ohmic currents]\label{remark:NonEqThermo-ohmicCurrents}
 We now assume an electric field~$\fieldEL$ generates the free current density~$\currentEL$ such that the free current density is proportional to this generating electric field. In this case the free current density
 is given according to Ohm's law, cf.~\cite{Masliyah-book, LifshitzLandau-book2, PrigogineKondepudi-book}, by
 \begin{align*}
  \currentEL=\constEL\eps_r\fieldEL,
 \end{align*}
 where the constant of proportionality is the so-called electric permittivity~$\constEL\eps_r$.  Furthermore, in stationary situations ($\dert\chargeEL=0$) with vanishing barycentric flow ($\fieldF\equiv0$),
 the charge conservation equation\eqref{eq:NonEqThermo-charge-chargeBalanceTot} reduces to Gauss's law, cf.~\cite{Masliyah-book, LifshitzLandau-book2, LopezDavalos-book},  
 \begin{align*}
  &\grad\cdot\brac{\constEL\eps_r\fieldEL}=0 \qquad  &\text{in  }&\Omega \qquad\sim\sqbrac{Cm^{-3}s^{-1}}.  \\[-14.0mm]
 \end{align*} 
 \hfill$\square$
\end{remark} 
\fancyhead[L]{Momentum conservation}
\section{Momentum Conservation Equations}\label{sec:NonEqThermo-mom}
%
Firstly, we recall from \cite{EckGarckeKnabner-book, oden-book}, that for vector fields~$\vecv$ the material derivatives $\derm[t_l]$, $\derm$ with respect to 
the transporting velocity fields~$\fieldF_l$, $\fieldF$ are defined by
\begin{subequations}
\begin{align}
 \derm[t_l]\vecv &:= \dert\vecv + \fieldF_l\cdot\grad\vecv, \label{eq:NonEqThermo-mom-matDerl} \\[2.0mm]
 \derm\vecv      &:= \dert\vecv + \fieldF\cdot\grad\vecv.   \label{eq:NonEqThermo-mom-matDer}
\end{align}
\end{subequations}
Moreover, we obtain%
\footnote{Here, we denote the dyadic product of two vectors $\vecv,\vecw\in\setR^n$ by $\vecv\otimes\vecw\in\setR^{n\times n}$. 
          Furthermore, we use the product rule $\grad\cdot\brac{\vecv\otimes\vecw}=\vecv\cdot\grad\vecw + \vecw \grad\cdot\vecv$.  } 
with \eqref{eq:NonEqThermo-mass-massBalance}
\begin{align*}
 \rol\derm[t_l]\vecv  = \dert(\rol\vecv) + \grad\cdot(\rol\fieldF_l\otimes\vecv) - \vecv\sqbrac{ \dert\rol + \grad\cdot\brac{\rol\fieldF_l} } 
                      = \dert(\rol\vecv) + \grad\cdot(\rol\fieldF_l\otimes\vecv) - \rl\vecv~,
\end{align*}
and we come with \eqref{eq:NonEqThermo-mass-massBalanceTot} to
\begin{align*}
   \rol[]\derm\vecv = \dert(\rol[]\vecv) + \grad\cdot\brac{\rol[]\fieldF\otimes\vecv} - \vecv\sqbrac{\dert\rol[] + \grad\cdot\brac{\rol[]\fieldF}} 
                    = \dert(\rol[]\vecv) + \grad\cdot\brac{\rol[]\fieldF\otimes\vecv}.                  
\end{align*}
Hence, we have the identities
\begin{subequations}
\begin{align}
 \rol\derm[t_l]\vecv &= \dert(\rol\vecv) + \grad\cdot(\rol\fieldF_l\otimes\vecv) - \rl\vecv,  \label{eq:NonEqThermo-mom-matDerl-Identity} \\[2.0mm]
 \rol[]\derm\vecv    &= \dert(\rol[]\vecv) + \grad\cdot\brac{\rol[]\fieldF\otimes\vecv} .     \label{eq:NonEqThermo-mom-matDer-Identiy}
\end{align}
\end{subequations}
We now turn to the basic equations for momentum conservation of the $l$th constituent. These equations are given by Newton's second law in Eulerian coordinates%
\footnote{We rewrite Newton's second law $\ml[]\vecx^{\prime\prime}=\vecF$ with the velocity~$\vecv$ as $\ml[]\vecv^\prime=\vecF$. Then, we switch from Lagrange coordinates to Eulerian coordinates. 
          Thereby ${}^\prime$ transforms to $\derm$. Finally, by using the mass density~$\rol[]$ instead of $\ml[]$, we obtain $\rol[]\derm\vecv=\vecF$.},  
\begin{subequations}
\begin{align}\label{eq:NonEqThermo-mom-newton1}
    \overbrace{\rol\derm[t_l]\fieldF_l}^{\begin{subarray}{l}\text{total change in} \\ \text{momentum} \end{subarray} } 
  = \overbrace{\vecF_l}^{\begin{subarray}{l} \text{total force} \\ \text{acting on the} \\ \text{$l$th constituent} \end{subarray} }  
    \qquad \sim\sqbrac{Nm^{-3}}. 
\end{align}
Newton's second law can be equivalently written with~\eqref{eq:NonEqThermo-mom-matDerl-Identity} in conservative from
\begin{align}\label{eq:NonEqThermo-mom-newton2}
   \overbrace{\dert(\rol\fieldF_l)}^{\begin{subarray}{l}\text{temporal change in} \\ \text{momentum} \end{subarray}}  
 + \grad\cdot\overbrace{\sqbrac{\rol\fieldF_l\otimes\fieldF_l} }^{\begin{subarray}{l}\text{momentum flux} \\ \text{tensor due to} \\ \text{convection} \end{subarray} } 
 = \overbrace{\vecF_l}^{\begin{subarray}{l} \text{total force} \\ \text{acting on the} \\ \text{$l$th constituent} \end{subarray} }   
 + \overbrace{ \rl\fieldF_l }^{\begin{subarray}{l}\text{Momentum} \\ \text{transfer due to} \\ \text{mass production} \end{subarray}}
   \qquad \sim\sqbrac{Nm^{-3}}. 
\end{align}
\end{subequations}
These equations are the general momentum conservation equations for the momentum densities~$\rol\fieldF_l$. We now extend the list of assumptions to introduce the ansatzes for the forces densities.  
\begin{enumerate}[align=left, leftmargin=*, topsep=2.0mm, label={(\labelA\arabic*)}, start=\value{countModAssump}+1]
  \item \label{Assump:Forces}%
        \textbf{Force density: } the $l$th total force density additively consist of the following contributions 
        \begin{align*}
          \vecF_l =   \overbrace{ \vecF^{\mathrm{stress}}_l }^{\begin{subarray}{l} \text{Stresses resp.} \\ \text{pressure forces} \end{subarray}} 
                    + \overbrace{ \vecF^{\mathrm{el}}_l }^{\begin{subarray}{l}\text{Electric} \\ \text{forces} \end{subarray}} 
                    + \overbrace{ \vecF^{\mathrm{int}}_l }^{\begin{subarray}{l}\text{Molecular} \\ \text{interaction} \\ \text{forces} \end{subarray}} 
                      \qquad\sim\sqbrac{Nm^{-3}}.
        \end{align*}
        Here, we suppose that the contributions due to stresses can be modeled with the stress tensors~$\stressTotl$, i.e., $\vecF^{\mathrm{stress}}_l:= \grad\cdot\stressTotl$.
        Furthermore, we assume that the electric force contributions arise due an electric field~$\fieldEL$, which is present inside the mixture. Thus, 
        we suppose $\vecF^{\mathrm{el}}_l = \constCharge\zl\nl\fieldEL$ with the valency~$\zl$ and the elementary charge~$\constCharge$. Altogether, we obtain the ansatzes
        \begin{align*}
          \vecF_l + \fieldF_l\rl = \grad\cdot\stressTotl + \constCharge\zl\nl\fieldEL + \vecF^{\mathrm{int}}_l + \rl\fieldF_l \qquad\sim\sqbrac{Nm^{-3}}.
        \end{align*}     
\setcounter{countModAssump}{\value{enumi}}%
\end{enumerate}
\begin{remark}[Internal interaction forces]\label{remark:NonEqThermo-internalForces}
  We note that the internal interaction forces between the constituents consist of two parts: The first contribution is due intermolecular interaction forces~$\vecF^{\mathrm{int}}_l$ 
  and the second contribution arises due to momentum transfer $\rl\fieldF_l$ between the constituents during mass production. According to \ref{Assump:momentum}, the total contribution 
  $\sum_l\vecF^{\mathrm{int}}_l+\rl\fieldF_l$ of internal interaction forces vanishes. Otherwise the mixture would have the unphysical ability to intrinsically produce or reduce its own momentum. 
  Besides this restriction, we do not involve further assumptions about the intermolecular interaction forces~$\vecF^{\mathrm{int}}_l$.%
\hfill$\square$
\end{remark} 
Substituting assumption~\ref{Assump:Forces} into Newton's second law~\eqref{eq:NonEqThermo-mom-newton1} and~\eqref{eq:NonEqThermo-mom-newton1}, leads to the momentum conservation equations for the constituents 
in nonconservative form 
\begin{align*}
  \rol\derm[t_l]\fieldF_l= \grad\cdot\stressTotl + \constCharge\zl\nl\fieldEL + \vecF^{int}_l \quad\sim\sqbrac{Nm^{-3}}             
\end{align*}
and in conservative form
\begin{align}\label{eq:NonEqThermo-mom-momBalance}
  \dert(\rol\fieldF_l) + \grad\cdot\brac{\rol\fieldF_l\otimes\fieldF_l} 
              = \grad\cdot\stressTotl + \constCharge\zl\nl\fieldEL + \vecF^{\mathrm{int}}_l + \fieldF_l\rl \quad\sim\sqbrac{Nm^{-3}}.                
\end{align}
We henceforth refer to these equations as the \textit{momentum conservation equations} for the constituents of the mixture.
\medskip
\par
Next, we derive the momentum conservation equations for the barycentric momentum density~$\rol[]\fieldF$ of the mixture. For that purpose, we sum over equations~\eqref{eq:NonEqThermo-mom-momBalance}. 
Together with \ref{Assump:momentum} and the free charge density~$\chargeEL:=\sum_l\constCharge\zl\nl$, we thereby arrive at  
\begin{align*}
	      &  \sum_l\sqbrac{\dert(\rol\fieldF_l) + \grad\cdot(\rol\fieldF_l\otimes\fieldF_l) }  
	       = \sum_l\sqbrac{\grad\cdot\stressTotl + \constCharge\zl\nl\fieldEL + \vecF^{\mathrm{int}}_l + \fieldF_l\rl }    \\[0.2cm]
\Equivalent~ &   \dert\brac{\rol[]\fieldF} + \grad\cdot\sqbrac{ \sum_l \rol\fieldF_l\otimes\fieldF_l } 
	      =  \chargeEL\fieldEL + \grad\cdot\sum_l\stressTotl        . 
\end{align*} 
In particular, for the sum of the momentum flux density tensors, we obtain with \eqref{eq:NonEqThermo-mass-sumCondition}
\begin{align*}
   \underbrace{ \sum_l \rol\fieldF_l\otimes\fieldF_l }_{ \begin{subarray}{l}\text{total momentum flux}  \\ \text{due to convection} \end{subarray} }
 = \underbrace{ \rol[]\fieldF\otimes\fieldF }_{ \begin{subarray}{l}\text{momentum flux} \\ \text{due to} \\ \text{barycentric flow} \end{subarray} } 
 + \underbrace{ \sum_l \rol(\fieldF_l-\fieldF)\otimes(\fieldF_l-\fieldF) }_{ \begin{subarray}{l}\text{momentum flux due to} \\ \text{mixing (=due to drift velocities)} \end{subarray} }.
\end{align*}
Substituting this identity in the above equations for $\rol[]\fieldF$, leads to
\begin{align*}
     \dert(\rol[]\fieldF) + \grad\cdot\brac{\rol[]\fieldF\otimes\fieldF} 
  =  \chargeEL\fieldEL + \grad\cdot\sum_l\brac{\stressTotl -\rol(\fieldF_l-\fieldF)\otimes(\fieldF_l-\fieldF)}    ~~\sim\sqbrac{Nm^{-3}}. 
\end{align*}
We now introduce the stress tensor~$\stressTotl[]$ of the mixture together with an additional assumption about the structure of~$\stressTotl[]$. 
\begin{enumerate}[align=left, leftmargin=*, topsep=2.0mm, label={(\labelA\arabic*)}, start=\value{countModAssump}+1]
  \item \label{Assump:mixtureStress1}%
        \textbf{Mixture stress tensor 1: } Following \cite{TruesdellToupin-book}, we define the mixture stress tensor by  
        \begin{align*}
          \stressTotl[] := \sum_l\brac{\stressTotl - \rol(\fieldF_l-\fieldF)\otimes(\fieldF_l-\fieldF)}.
        \end{align*}
 \item \label{Assump:mixtureStress2}
        \textbf{Mixture stress tensor 2: } Analogous to \ref{Assump:stresses}, we assume that a symmetric viscous mixture stress tensor~$\stressViscl[]$ and a hydrostatic mixture pressure~$\pressHydrl[]$
        exists such that the mixture stress tensor~$\stressTotl[]$ is given by $\stressTotl[]:=-\pressHydrl[]\mathds{1}+\stressViscl[]$.
\setcounter{countModAssump}{\value{enumi}}%
\end{enumerate}
Finally, substituting \ref{Assump:mixtureStress1} into the preceding equation, leads to
\begin{align}\label{eq:NonEqThermo-mom-momBalanceTot}
    \dert(\rol[]\fieldF) + \grad\cdot\brac{\rol[]\fieldF\otimes\fieldF} = \grad\cdot\stressTotl[] + \chargeEL\fieldEL \qquad\sim\sqbrac{Nm^{-3}}. 
\end{align}
These equations are the momentum conservation equations for the barycentric momentum density~$\rol[]\fieldF$ of the mixture.
\begin{remark}[Cauchy's second law of motion]
  We note that due to  \ref{Assump:stresses} and \ref{Assump:mixtureStress1} the partial stress tensors~$\stressTotl$ and the mixture stress tensor~$\stressTotl[]$ are symmetric. 
  This is exactly Cauchy's second law of motion, cf. \cite{TruesdellToupin-book}. However, according to \cite{TruesdellToupin-book} it would be sufficient to have a symmetric mixture 
  stress tensor~$\stressTotl[]$. Together with \ref{Assump:mixtureStress1}, this would allow for nonsymmetric partial stress tensors~$\stressTotl$ as long as their sum remains symmetric, i.e., 
  \begin{align*}
     \sum_l\stressTotl = \brac{\sum_l\stressTotl}^\top.   \\[-12.0mm]
  \end{align*}

\hfill$\square$
\end{remark} 
Finally, we investigate some consequences of assumptions~\ref{Assump:mixtureStress1} and \ref{Assump:mixtureStress2}.
\begin{remark}[Extended Dalton's law and extended Raoult's law]\label{remark:NonEqThermo-dalton}
  We confine ourselves to $\trace{\stressViscl[]}=0$ and $\trace{\stressViscl}=0$. Here, the total pressure~$\pressTotl[]$ coincides with the hydrostatic pressure~$\pressHydrl[]$. 
  Combining \ref{Assump:mixtureStress2} and \ref{Assump:mixtureStress1} with the definition of the partial total pressures~$\pressTotl$ and the total mixture pressure~$\pressTotl[]$ 
  in \ref{Assump:thermoPressure} shows with $\pressTotl=\pressHydrl$ and $\pressTotl[]=\pressHydrl[]$
  \begin{align}\label{eq:NonEqThermo-mom-extendedDaltonsLaw}
  \pressHydrl[]  &= -\frac{1}{n}\trace{ \stressTotl[] } 
	          = -\frac{1}{n}\sum_l\trace{\stressTotl - \rol(\fieldF_l-\fieldF)\otimes(\fieldF_l-\fieldF) }  \nonumber\\
	         &=  \sum_l \pressHydrl + \frac{1}{n}\sum_l\trace{ \rol(\fieldF_l-\fieldF)\otimes(\fieldF_l-\fieldF) } 
	          =  \sum_l \pressHydrl + \frac{1}{n}\sum_l\rol\abs{\fieldF_l-\fieldF}^2 .
  \end{align}
  This identity is a generalized version of Dalton's law for the pressure~$\pressHydrl[]$ of mixtures, cf. \cite{Atkins-book, PrigogineKondepudi-book}. 
  Furthermore, assuming that the partial pressure~$\pressHydrl$ of a constituent in the mixture is given by~$\pressHydrl =\pressHydrl^\ast \yl$,
  where $\pressHydrl^\ast$ is the partial pure substance total pressure, the above equation leads to the following extended version of Raoult's law, cf. \cite{Atkins-book, PrigogineKondepudi-book}
  \begin{align*}
  \pressHydrl[]  = \sum_l \pressHydrl^\ast\yl + \frac{1}{n}\sum_l\rol\abs{\fieldF_l-\fieldF}^2 .\\[-12.0mm]
  \end{align*}
\hfill$\square$
\end{remark}
In \ref{Assump:mixtureStress1}, we repeated from \cite{TruesdellToupin-book} the fundamental definition of the total mixture stress tensor~$\stressTotl[]$ in terms of previously defined quantities of the constituents. 
Analogously, we can adopt extended Dalton's law~\eqref{eq:NonEqThermo-mom-extendedDaltonsLaw} as reasonable definition for the hydrostatic mixture pressure~$\pressHydrl[]$. Inserting this formula 
into \ref{Assump:mixtureStress2} leads with \ref{Assump:stresses} and \ref{Assump:mixtureStress1} for the viscous mixture tensor to 
\begin{align*}
             &  -\sum_l\pressHydrl\mathds{1} + \sum_l\brac{\stressViscl - \rol(\fieldF_l-\fieldF)\otimes(\fieldF_l-\fieldF)} \\
             &= \stressTotl[] = -\pressHydrl[]\mathds{1} +\stressViscl[] = -\sum_l\pressHydrl\mathds{1} - \frac{1}{n}\sum_l\rol\abs{\fieldF_l-\fieldF}^2\mathds{1} +\stressViscl[] ~.
\end{align*}
Hence, corresponding to extended Dalton's law~\eqref{eq:NonEqThermo-mom-extendedDaltonsLaw} for hydrostatic mixture pressure~$\pressHydrl[]$, 
we obtain for the viscous mixture stress tensor~$\stressViscl[]$ the definition  
\begin{align}
  \stressViscl[] &:= \sum_l \brac{\stressViscl - \rol(\fieldF_l-\fieldF)\otimes(\fieldF_l-\fieldF)} + \frac{1}{n}\sum_l\rol\abs{\fieldF_l-\fieldF}^2\mathds{1} \label{eq:NonEqThermo-mom-viscStressMix}~.
\end{align}
From this equation follows, that the trace of $\stressViscl[]$ solely depends on the sum of the traces $\trace{\stressViscl}$, since we have
\begin{align*}
  \trace{ \stressViscl[] } &= \sum_l\trace{ \stressViscl - \rol(\fieldF_l-\fieldF)\otimes(\fieldF_l-\fieldF)} + \frac{1}{n}\sum_l\rol\abs{\fieldF_l-\fieldF}^2\trace{\mathds{1}}  \\
	                   &= \frac{1}{n}\trace{\sum_l \stressViscl} - \sum_l\rol\abs{\fieldF_l-\fieldF}^2 + \sum_l\rol\abs{\fieldF_l-\fieldF}^2 
                            = \frac{1}{n}\trace{\sum_l \stressViscl}~.
\end{align*}
We now combine this identity with \ref{Assump:thermoPressure}, \ref{Assump:mixtureStress1}, \ref{Assump:mixtureStress2}, \eqref{eq:NonEqThermo-mom-extendedDaltonsLaw}, and \eqref{eq:NonEqThermo-mom-viscStressMix}.
Thereby, we obtain between the total mixture pressure~$\pressTotl[]$ of the mixture and hydrostatic mixture pressure~$\pressHydrl[]$ the connection  
\begin{align}\label{eq:NonEqThermo-mom-pressTotpressHydr}
\pressTotl[]  &= -\frac{1}{n}\trace{ \stressTotl[] }
               = -\frac{1}{n}\sum_l\trace{ -\pressHydrl\mathds{1} +\stressViscl -\rol(\fieldF_l-\fieldF)\otimes(\fieldF_l-\fieldF) }  \nonumber\\
              &= \sum_l\pressHydrl + \frac{1}{n}\sum_l\rol\abs{\fieldF_l-\fieldF}^2 -\frac{1}{n}\sum_l\trace{ \stressViscl } 
               = \pressHydrl[] -\frac{1}{n}\trace{ \stressViscl[] } .
\end{align}
Thus, analogously to \cref{remark:NonEqThermo-distinguishPressures} the hydrostatic mixture pressure~$\pressHydrl[]$ and the total mixture pressure~$\pressTotl[]$ coincide, 
provided we have a traceless viscous mixture stress tensor~$\stressViscl[]$. 
\fancyhead[L]{Energy conservation}
\section{Energy Conservation Equations}\label{sec:NonEqThermo-energy}
First of all, we introduce a fundamental assumption about the energy densities of each constituent, cf. \cite{DeGrootMazur-book}.
\begin{enumerate}[align=left, leftmargin=*, topsep=2.0mm, label={(\labelA\arabic*)}, start=\value{countModAssump}+1]
  \item \label{Assump:energy1}%
        \textbf{Energy contributions of the constituents: } We suppose that for every constituent the energy contribution~$\rol\energyTotl$ additively decomposes into three parts: 
         firstly the kinetic energy~$\frac{1}{2}\rol\abs{\fieldF_l}^2$, 
         secondly the electric potential energy~$\constCharge\zl\nl\potEL$, 
         and thirdly the internal energy~$\rol\energyIntl$. This means, we assume the fundamental ansatzes
	 \begin{align}\label{eq:NonEqThermo-energy-ansatz}
	  \rol\energyTotl = \frac{1}{2}\rol\abs{\fieldF_l}^2 + \constCharge\zl\nl\potEL + \rol\energyIntl \qquad\sim\sqbrac{Jm^{-3}}.
	 \end{align}
 \item \label{Assump:energy1a}%
        \textbf{Structure of the specific internal energies: } We suppose that the specific internal energies~$\energyIntl$ of the constituents additively decompose 
        into a pure substance part and a part due to mixing, i.e., we assume 
	 \begin{align}\label{eq:NonEqThermo-energy-ansatzPureMixing}
	     \rol\energyIntl = \rol\energyIntPurel + \rol\energyIntMixl \qquad\sim\sqbrac{J m^{-3}}. 
	 \end{align}
\setcounter{countModAssump}{\value{enumi}}%
\end{enumerate}
To obtain the respective energy densities of the mixture, we sum over $l\in\cbrac{1,\ldots,L}$. For the electric potential energy densities this leads with \eqref{eq:NonEqThermo-charge-DefFreeCharge} to 
\begin{align*}
  \sum_l \constCharge\zl\nl\potEL = \chargeEL\potEL~.
\end{align*}
For the kinetic energy densities, we obtain with $\abs{\vecv+\vecw}^2= \abs{\vecv}^2 + \abs{\vecw}^2 + 2\vecv\cdot\vecw$ and \eqref{eq:NonEqThermo-mass-sumCondition}
\begin{align*}
  \frac{1}{2}\sum_l\rol\abs{\fieldF_l}^2 = \frac{1}{2}\rol[]\abs{\fieldF}^2 + \frac{1}{2}\sum_l\rol\abs{\fieldF_l-\fieldF}^2~.
\end{align*} 
This reveals, that the kinetic energy density of the mixture decomposes into a barycentric part and a part due to mixing. Henceforth, we consider the kinetic energy part 
due to mixing as an internal contribution and thus, we add this part to the internal energies. This means, we define the mixture total energy energy density~$\rol[]\energyTotl[]\sim[J/m^3]$, 
the mixture internal energy density~$\rol[]\energyIntl[]\sim[J/m^3]$, the mixture pure substance internal energy~$\rol[]\energyIntPurel[]\sim[J/m^3]$, 
and the mixture internal energy of mixing~$\rol[]\energyIntMixl[]\sim[J/m^3]$ by  
\begin{subequations}
\begin{align}
      \rol[]\energyTotl[] &:= \rol[]\sum_l \yl \energyTotl,                                                                                                   \\
      \rol[]\energyIntl[] &:= \rol[]\sum_l \yl \energyIntMixl + \rol[]\sum_l \yl \energyIntPurel +\frac{\rol[]}{2}\sum_l\yl\abs{\fieldF_l-\fieldF}^2~,        \\ 
  \rol[]\energyIntMixl[]  &:= \rol[]\sum_l \yl \energyIntMixl,                                                    \label{eq:NonEqThermo-energy-energyIntMix}  \\
  \rol[]\energyIntPurel[] &:= \rol[]\sum_l \yl \energyIntPurel+\frac{\rol[]}{2}\sum_l\yl\abs{\fieldF_l-\fieldF}^2 \label{eq:NonEqThermo-energy-energyIntPure}.  
\end{align}
\end{subequations}
In summary, the ansatz~\ref{Assump:energy1} for the constituent energy densities~$\rol\energyTotl$ leads for the mixture density~$\rol[]\energyTotl[]$ to
\begin{subequations}
\begin{align}
             & \overbrace{\rol[]\energyTotl[]}^{\begin{subarray}{l} \text{total} \\ \text{energy} \end{subarray}} 
              = \overbrace{\rol[]\energyIntl[]}^{\begin{subarray}{l}  \text{internal} \\ \text{energy} \end{subarray}} 
              + \overbrace{\chargeEL\potEL}^{\begin{subarray}{l}   \text{electric} \\ \text{potential} \\ \text{energy} \end{subarray}} 
              + \overbrace{\frac{1}{2}\rol[]\abs{\fieldF}^2}^{\begin{subarray}{l}  \text{barycentric} \\ \text{kinetic energy} \end{subarray}}  \label{eq:NonEqThermo-energy-ansatzTot}\\
 \Equivalent & \rol[]\energyTotl[] = \rol[]\energyIntPurel[] + \rol[]\energyIntMixl[] + \chargeEL\potEL + \frac{1}{2}\rol[]\abs{\fieldF}^2~. \label{eq:NonEqThermo-energy-ansatzTot-a}
\end{align}
\end{subequations}
Following \cite{DeGrootMazur-book, PrigogineKondepudi-book}, we now formulate the \textbf{\textit{first law of thermodynamics}}, which states that the total energy of a closed system is conserved. 
In differential form the general statement of the first law of thermodynamics is given by the balance equation
\begin{align}\label{eq:NonEqThermo-energy-totEnergyConservation}
 & \dert(\rol[]\energyTotl[]) +\grad\cdot\brac{\rol[]\energyTotl[]\fieldF +\energyTotlFlux} = 0 &\text{in  }\Omega \qquad\sim\sqbrac{J m^{-3}s^{-1}}.
\end{align}
Here, $\energyTotlFlux\sim[J/(m^2s)]$ is the energy flux. To obtain a more specific version of the first law of thermodynamics, it now remains to derive an explicit expression for the energy flux~$\energyTotlFlux$ 
in terms of the internal energy flux, the electric potential energy flux, and the kinetic energy flux. For that purpose, we subsequently derive evolution equations 
for each of the energy densities from \eqref{eq:NonEqThermo-energy-ansatzTot}. 
\begin{remark}[Energy conservation of closed systems]
  We note, that we can supplement equation \eqref{eq:NonEqThermo-energy-totEnergyConservation} with a no-flux boundary condition~$\vecJ_e\cdot\vecnu=0$ for $\vecJ_e:=\rol[]\energyTotl[]\fieldF+\energyTotlFlux$ 
  on the boundary~$\partial\Omega$. This boundary condition models a closed system, as we have no flow of energy across the boundary~$\partial\Omega$. Thus the system inside $\Omega$ is energetically separated from its 
  exterior $\setR^n\backslash\Omega$. Next, we integrate over $\Omega$ and we apply Gauss's divergence theorem, cf. \cite{EvansGariepy-book}. Thereby, we come with the above no-flux boundary condition to
  \begin{align*}
  0 &= \Intdx{\dert(\rol[]\energyTotl[])+ \grad\cdot\vecJ_e~} = \derr\Intdx{\rol[]\energyTotl[]~} + \IntdS[\partial\Omega]{\vecJ_e\cdot\vecnu~} = \derr\Intdx{\rol[]\energyTotl[]~}~.
  \end{align*}
  This is exactly the global statement of the first law of thermodynamics, i.e., of energy conservation for closed systems. Hence, the differential version~\eqref{eq:NonEqThermo-energy-totEnergyConservation}
  of the first law of thermodynamics contains the preceding global version for closed systems.
\hfill$\square$
\end{remark}
\par
In this passage, we derive the evolution equation for the electric potential energy density~$\chargeEL\potEL$ of the mixture by recalling the 
charge conservation equation~\eqref{eq:NonEqThermo-charge-chargeBalanceTot}, which we multiplying by~$\potEL$. This results in 
\begin{align*}
 \potEL\dert\chargeEL + \potEL\grad\cdot\brac{\chargeEL\fieldF + \currentEL}  =  0  \qquad\sim\sqbrac{Jm^{-3}s^{-1}}.   
\end{align*} 
Applying the product rule on the left-hand side with \ref{Assump:electrostatics} and the additional assumption of a stationary electric potential, i.e. $\dert\potEL=0$, 
shows that the evolution equation for the electric potential energy density is given by
\begin{align}\label{eq:NonEqThermo-energy-electricEnergyPDE}
 \dert\brac{\chargeEL\potEL} + \grad\cdot\brac{\chargeEL\potEL\fieldF + \potEL\currentEL}  =  - \chargeEL\fieldEL\cdot\fieldF - \currentEL\cdot\fieldEL \qquad\sim\sqbrac{Jm^{-3}s^{-1}}.   
\end{align}
\par
Next, we derive an evolution equation for the barycentric kinetic energy~$\frac{1}{2}\rol[]\abs{\fieldF}^2$ of the mixture. For that purpose, we multiply the momentum 
conservation equations~\eqref{eq:NonEqThermo-mom-momBalanceTot} by $\fieldF$. Thereby, we obtain 
\begin{align*}
  \dert(\rol[]\fieldF)\cdot\fieldF + \grad\cdot\brac{\rol[]\fieldF\otimes\fieldF}\cdot\fieldF = \brac{\grad\cdot\stressTotl[]}\cdot\fieldF + \chargeEL\fieldEL\cdot\fieldF \qquad\sim\sqbrac{Jm^{-3}s^{-1}}. 
\end{align*}
For the first term on the right-hand side, we get with the product rule 
\begin{align*}
 \brac{\grad\cdot\stressTotl[]}\cdot\fieldF=\grad\cdot\brac{\stressTotl[]\fieldF} -\stressTotl[]:\grad\fieldF, 
\end{align*}
where $:$ denotes the scalar product $A:B=\trace{A^\top B}$ of two matrices $A,B\in\setR^{n\times n}$. On the left-hand side, we receive for the first term  
\begin{align*}
    \dert(\rol[]\fieldF)\cdot\fieldF
  = \abs{\fieldF}^2\dert\rol[] + \rol[]\dert\brac{\frac{1}{2}\abs{\fieldF}^2} 
  = \dert\brac{\frac{1}{2}\rol[]\abs{\fieldF}^2} + \frac{1}{2}\abs{\fieldF}^2\dert\rol[].
\end{align*}
Furthermore, with the calculus identities 
\begin{align*}
 \grad\cdot(\vecv\otimes\vecw)=\vecv\cdot\grad\vecw+\vecw\grad\cdot\vecv \quad\text{ and  }\quad \vecv\cdot\grad\vecv=\frac{1}{2}\grad\abs{\vecv}^2-\vecv\times(\grad\times\vecv),
\end{align*}
we obtain for the second term on the left-hand side with $\vecv\bot(\vecv\times\vecw)$
\begin{align*}
      \grad\cdot\brac{\rol[]\fieldF\otimes\fieldF}\cdot\fieldF 
    = \sqbrac{ \rol[]\fieldF\cdot\grad\fieldF + \fieldF\grad\cdot\brac{\rol[]\fieldF} }\cdot\fieldF 
    = \grad\cdot\brac{\frac{1}{2}\rol[]\abs{\fieldF}^2\fieldF} + \frac{1}{2}\abs{\fieldF}^2\grad\cdot\brac{\rol[]\fieldF}~.
\end{align*}
Hence, we thereby arrive with the mass conservation equation~\eqref{eq:NonEqThermo-mass-massBalanceTot} at the evolution equation for the barycentric kinetic energy density
\begin{align}\label{eq:NonEqThermo-energy-kineticEnergyPDE}
 \dert\!\brac{\frac{1}{2}\rol[]\abs{\fieldF}^2} + \grad\cdot\!\brac{\frac{1}{2}\rol[]\abs{\fieldF}^2\fieldF - \stressTotl[]\fieldF} = - \stressTotl[]:\grad\fieldF + \chargeEL\fieldEL\cdot\fieldF
 ~~\sim\sqbrac{Jm^{-3}s^{-1}}.
\end{align}
\par
Finally, we establish the evolution equation for the internal energy density~$\rol[]\energyIntl[]$ of the mixture. The general statement of this evolution equation 
is given with the heat flux~$\vecq\sim[J/(m^2s)]$ and the internal energy production rate~$h\sim[J/(m^3s)]$ by
\begin{align}\label{eq:NonEqThermo-energy-internalEnergyPDE-prelim}
 & \dert\brac{\rol[]\energyIntl[]} +\grad\cdot\brac{\rol[]\energyIntl[]\fieldF +\vecq} = h \qquad\sim\sqbrac{J m^{-3}s^{-1}}.
\end{align}
However, this equation is an abstract balance statement and we have to set up more specific expression for the heat flux~$\vecq$ and the internal energy production rate~$h$. 
This is carried out in the next step.
\medskip
\par
We now sum up the evolution equations~\eqref{eq:NonEqThermo-energy-electricEnergyPDE}, \eqref{eq:NonEqThermo-energy-kineticEnergyPDE}, and \eqref{eq:NonEqThermo-energy-internalEnergyPDE-prelim}. Thereby, we obtain
\begin{align*}
  &  \dert\!\brac{\rol[]\energyIntl[]+\chargeEL\potEL+\rol[]\abs{\fieldF}^2} 
   + \grad\cdot\!\brac{ \sqbrac{\rol[]\energyIntl[]+\chargeEL\potEL+\frac{1}{2}\rol[]\abs{\fieldF}^2}\fieldF + \vecq + \potEL\currentEL - \stressTotl[]\fieldF } \\
  &= h - \currentEL\cdot\fieldEL - \stressTotl[]:\grad\fieldF ~.
\end{align*}
Substituting the ansatz~\eqref{eq:NonEqThermo-energy-ansatzTot} into the latter equation yields
\begin{align*}
    \dert\brac{\rol[]\energyTotl[]} 
   + \grad\cdot\!\brac{ \rol[]\energyTotl[]\fieldF + \vecq + \potEL\currentEL - \stressTotl[]\fieldF } 
  = h - \currentEL\cdot\fieldEL - \stressTotl[]:\grad\fieldF ~.
\end{align*}
Comparing this equation with the first law of thermodynamics~\eqref{eq:NonEqThermo-energy-totEnergyConservation} reveals
\begin{align*}
 \energyTotlFlux := \vecq + \potEL\currentEL - \stressTotl[]\fieldF \qquad\text{and}\qquad  h := \currentEL\cdot\fieldEL + \stressTotl[]:\grad\fieldF~.             
\end{align*}
In summary we have shown that the \textbf{\textit{first law of thermodynamics}} is given by
\begin{align}\label{eq:NonEqThermo-energy-firstLaw}
  \dert\brac{\rol[]\energyTotl} + \grad\cdot\!\brac{ \rol[]\energyTotl[]\fieldF + \vecq + \potEL\currentEL - \stressTotl[]\fieldF } = 0 \qquad\sim\sqbrac{J m^{-3}s^{-1}}.
\end{align}
\begin{remark}[Equation for the internal energy]
 The preceding definition of internal energy production rate~$h$ shows, that the evolution equation~\eqref{eq:NonEqThermo-energy-internalEnergyPDE-prelim} for the internal energy density now reads as
 \begin{align}\label{eq:NonEqThermo-energy-internalEnergyPDE}
  & \dert\brac{\rol[]\energyIntl[]} +\grad\cdot\brac{ \rol[]\energyIntl[]\fieldF + \vecq } = \currentEL\cdot\fieldEL + \stressTotl[]:\grad\fieldF ~~\sim\sqbrac{J m^{-3}s^{-1}}.
 \end{align}
  Moreover, with \eqref{eq:NonEqThermo-mass-massBalanceTot} and $\derm$ from \eqref{eq:NonEqThermo-mom-matDer}, we have the 
  identity~$\dert(\rol[]\energyIntl[]) +\grad\cdot(\rol[]\energyIntl[]\fieldF)=\rol[]\derm\energyIntl[]$.  
  Thus, we can equivalently write the evolution equation for the internal energy density in nonconservative from as
  \begin{align}\label{eq:NonEqThermo-energy-internalEnergyPDE2}
   & \rol[]\derm\energyIntl[] = - \grad\cdot\vecq + \currentEL\cdot\fieldEL + \stressTotl[]:\grad\fieldF \qquad\sim\sqbrac{J m^{-3}s^{-1}}. \\[-14.0mm] \nonumber 
  \end{align}
\hfill$\square$
\end{remark}
\fancyhead[L]{Entropy evolution}
\section{Entropy Evolution Equation}\label{sec:NonEqThermo-entropy}
So far, we established in \cref{sec:NonEqThermo-mass} -- \cref{sec:NonEqThermo-energy} several conservation equations, to which we refer as the governing equations of the considered mixtures. 
Although these equations describe the evolution of the conserved quantities, we can not deduce from these equations any restriction about admissible direction of the underlying physical processes,
cf.~\cite{FeinbergLavine, LiebYngvason1999, LifshitzLandau-book5, oden-book, Tadmor-book, TruesdellToupin-book, Wilmanski-book}. More precisely, in order to come to reasonable statements about admissible 
directions of physical processes, we have to introduce an other quantity: The specific entropy~$\entropy\sim[J/(K\,kg)]$. With this quantity we now state the following assumptions. 
\begin{enumerate}[align=left, leftmargin=*, topsep=1.0mm, itemsep=-2.0mm, label={(\labelA\arabic*)}, start=\value{countModAssump}+1]
  \item \label{Assump:energy2}%
        \textbf{Functional dependency of the specific internal energies 1: } We suppose that the specific internal energies~$\energyIntl\sim[J/kg]$ are functions of the specific entropy~$\entropy\sim[J/(K\,kg)]$, 
        the specific volume~$\volumeSpecific\sim[m^3/kg]$ and the mass fractions~$\yl$, i.e., 
	 \begin{align}\label{eq:NonEqThermo-entropy-ansatz}
	   \energyIntl :\left\{\begin{array}{c l l} \setR\times\setR\times\setR^L                                        &\rightarrow & \setR~, \\
                                                    \brac{\entropy, \volumeSpecific, \yl[1],\ldots, \yl[L]} & \mapsto    & \energyIntl\brac{\entropy, \volumeSpecific, \yl[1],\ldots, \yl[L]}~.
	                  \end{array}\right.
	 \end{align}
  \item \label{Assump:energy2a}%
        \textbf{Functional dependency of the specific internal energy 2: } Moreover, combining \ref{Assump:energy1a} and \ref{Assump:energy2}, we assume that 
        we have  
	 \begin{align}\label{eq:NonEqThermo-entropy-ansatzEnergyMixPure}
	   \energyIntl\brac{\entropy, \volumeSpecific, \yl[1],\ldots, \yl[L]} 
         = \energyIntl\brac{\entropy, \volumeSpecific, \yl}
         = \energyIntPurel\brac{\entropy, \volumeSpecific} + \energyIntMixl\brac{\yl}~.
	 \end{align} 
\setcounter{countModAssump}{\value{enumi}}%
\end{enumerate}
These assumptions have far reaching consequences and need further clarifications: 
\begin{enumerate}[align=left, label=(\roman*), leftmargin=*, topsep=2.0mm, itemsep=-1.2mm]
  \item Firstly, \ref{Assump:energy2} and \ref{Assump:energy2a} lead with \eqref{eq:NonEqThermo-energy-energyIntMix} and \eqref{eq:NonEqThermo-energy-energyIntPure} 
       for the mixture variables to
       \begin{subequations}
        \begin{align}
          \energyIntPurel[] &= \energyIntPurel[]\brac{\entropy, \volumeSpecific, \yl[1],\ldots, \yl[L]} =  \sum_l\yl\energyIntPurel\brac{\entropy, \volumeSpecific} + \frac{1}{2}\sum_l\yl\abs{\fieldF_l-\fieldF}^2~, \label{eq:NonEqThermo-entropy-energyIntPure} \\
          \energyIntMixl[]  &= \energyIntMixl[]\brac{\yl[1],\ldots, \yl[L]}  = \sum_l\yl\energyIntMixl\brac{\yl}~, \label{eq:NonEqThermo-entropy-energyIntMix}\\
          \energyIntl[]     &= \energyIntl[]\brac{\entropy, \volumeSpecific, \yl[1],\ldots, \yl[L]} = \energyIntPurel[]\brac{\entropy, \volumeSpecific, \yl[1],\ldots, \yl[L]} + \energyIntMixl[]\brac{\entropy, \volumeSpecific, \yl[1],\ldots, \yl[L]} \label{eq:NonEqThermo-entropy-energyInt} 
        \end{align}
        \end{subequations}
        Note that in equation~\eqref{eq:NonEqThermo-entropy-energyIntPure}, we treated the velocities $\fieldF_1,\ldots,\fieldF_L,\fieldF$ as parameters and for ease of readability, we omitted this parameter dependency in the notation. 
  \item Secondly, according to \eqref{eq:NonEqThermo-entropy-ansatz} the specific internal energy~$\energyIntl[]$ is a function, which is defined on the phase space 
	\begin{align*}
	\setR^{L+2}_+ =       \underbrace{\setR_+}_{ \begin{subarray}{l} \text{entropy} \\ \text{space} \end{subarray} } 
			\times\underbrace{\setR_+}_{ \begin{subarray}{l} \text{volume}  \\ \text{space} \end{subarray} }
			\times\underbrace{\setR^L_+~.}_{ \begin{subarray}{l} \text{mass fraction} \\ \text{space} \end{subarray} }
	\end{align*}
	Thus, the coordinates which determine the respective values of the specific internal energy are the specific entropy~$\entropy$, 
        the specific volume~$\volumeSpecific$, and the mass fractions~$\yl$. 
  \item Thirdly, with the usual differentiation rules for functions of several variables, the differential~$\dtot\energyIntl[]$, cf. \cite{Cartan-book, EckGarckeKnabner-book, Hirsch-book}, is given by
	\begin{align*}
	\dtot\energyIntl[] = \derd[\entropy]\energyIntl[]~\dtot\entropy + \derd[\volumeSpecific]\energyIntl[]~\dtot\volumeSpecific + \sum_l \derd[\yl]\energyIntl[]~\dtot\yl~.
	\end{align*}
	Defining the temperature~$\temp$, the pressure~$\pressHydrl[]$, and the chemical potentials~$\chempotl$ according to classical thermodynamics, 
	cf. \cite{EckGarckeKnabner-book, Lavenda-book, PrigogineKondepudi-book, Wilmanski-book}, by
	\begin{align*}
	            \temp  &:= \derd[\entropy]\energyIntl[]~ \qquad\sim\sqbrac{K},              \\
	   -\pressHydrl[]  &:= \derd[\volumeSpecific]\energyIntl[]~\qquad\sim\sqbrac{J m^{-3}},  \\
	        \chempotl  &:= \derd[\yl]\energyIntl[]~\qquad\sim\sqbrac{J \, kg^{-1}},
	\end{align*}
	and substituting these definitions into the preceding identity, we obtain the fundamental \textbf{\textit{Gibbs relation}}, cf. \cite{EckGarckeKnabner-book, Lavenda-book, PrigogineKondepudi-book, Wilmanski-book},
	\begin{align}\label{eq:NonEqThermo-entropy-gibbs}
	\dtot\energyIntl[] = \temp~\dtot\entropy - \pressHydrl[]~\dtot\volumeSpecific + \sum_l \chempotl~\dtot\yl ~.
	\end{align}
  \item Fourthly, in \cref{sec:NonEqThermo-energy} we treated the specific internal energy~$\energyIntl[]$, e.g., in \eqref{eq:NonEqThermo-energy-internalEnergyPDE} 
	as function of space and time, whereas according to \ref{Assump:energy2} the specific internal energy~$\energyIntl[]$ is a function of $(\entropy,\volumeSpecific,\yl[1],\ldots,\yl[L])$.  
	This apparent contradiction can be resolved with the so-called Nemytskii operator~$\mathcal{N}$, cf. \cite{Roubicek-book}. More precisely, in nonequilibrium systems the specific entropy~$\entropy$, 
	the specific volume~$\volumeSpecific$, and the mass fractions~$\yl$ are functions of space and time. This means, that the coordinates of the phase space, 
	on which the specific internal energy is defined, are variable coordinates in space and time. Thus, rigorously, we considered in \cref{sec:NonEqThermo-energy} the Nemytskii mapping~$\mathcal{N}\sqbrac{\energyIntl[]}\tx$
	\begin{align}\label{eq:NonEqThermo-entropy-nemytskii}
	\mathcal{N}\sqbrac{\energyIntl[]}\tx:=\energyIntl[]\brac{\entropy\tx, \volumeSpecific\tx, \yl[1]\tx,\ldots, \yl[L]\tx}~, 
	\end{align}
	and we denoted the Nemytskii mapping~$\mathcal{N}\sqbrac{\energyIntl[]}\tx$ by $\energyIntl[]\tx$. For ease of readability we henceforth denote the 
        Nemytskii mapping~$\mathcal{N}\sqbrac{\energyIntl[]}\tx$ again by $\energyIntl[]\tx$.
  \item Fifthly, with the material derivatives~$\derm$ from \eqref{eq:NonEqThermo-mom-matDer} and the usual differentiation rules for functions of several variables, we obtain from \eqref{eq:NonEqThermo-entropy-nemytskii}
	\begin{align}\label{eq:NonEqThermo-entropy-derivativeNemitskii}
	\derm\energyIntl[] = \derd[\entropy]\energyIntl[]~\derm\entropy + \derd[\volumeSpecific]\energyIntl[]~\derm\volumeSpecific + \sum_l \derd[\yl]\energyIntl[]~\derm\yl  \qquad\sim\sqbrac{J\,kg^{-1}s^{-1}}.
	\end{align}
	Assuming that the involved functions are sufficiently smooth allows to repeat the definitions 
       \begin{subequations}
	\begin{align}
	            \temp\tx  &:= \derd[\entropy]\energyIntl[]\tx~        \qquad\sim\sqbrac{K},             \label{eq:NonEqThermo-entropy-temp}   \\
	   -\pressHydrl[]\tx  &:= \derd[\volumeSpecific]\energyIntl[]\tx~ \qquad\sim\sqbrac{J m^{-3}},      \label{eq:NonEqThermo-entropy-press}  \\
	        \chempotl\tx  &:= \derd[\yl]\energyIntl[]\tx              \qquad\sim\sqbrac{J\, kg^{-1}},   \label{eq:NonEqThermo-entropy-chemPot}
	\end{align}
        \end{subequations}
	which naturally introduce the temperature~$\temp$, the pressure~$\pressHydrl[]$, and the chemical potentials~$\chempotl$ as functions of space and time. Note that we identify the pressure field~$\pressHydrl[]$ 
        with the hydrostatic pressure~$\pressHydrl[]$ from \cref{sec:NonEqThermo-mom}. Inserting these definitions in \eqref{eq:NonEqThermo-entropy-derivativeNemitskii} and multiplying by~$\rol[]$, leads to the \textbf{\textit{Gibbs relation for nonequilibrium systems}}
	\begin{align}\label{eq:NonEqThermo-entropy-gibbsNonEq}
	\rol[]\derm\energyIntl[] = \temp~\rol[]\derm\entropy - \pressHydrl[]~\rol[]\derm\volumeSpecific + \sum_l \chempotl~\rol[]\derm\yl \qquad\sim\sqbrac{Jm^{-3}s^{-1}}.
	\end{align}
\end{enumerate}
\begin{remark}[Nonautonomous ansatz]
 The autonomous ansatz~\eqref{eq:NonEqThermo-entropy-ansatz} in assumption~\ref{Assump:energy2} leads to the Nemytskii mapping~\eqref{eq:NonEqThermo-entropy-nemytskii}, 
 which was the starting point for the nonequilibrium Gibbs relation~\eqref{eq:NonEqThermo-entropy-gibbsNonEq}. However, we could choose the nonautonomous ansatz
 \begin{align} 
    \brac{t,x,\entropy, \volumeSpecific, \yl[1],\ldots, \yl[L]} ~\mapsto~ \energyIntl[]\brac{t,x,\entropy, \volumeSpecific, \yl[1],\ldots, \yl[L]}
    \tag{${\text{\ref{eq:NonEqThermo-entropy-ansatz}}}^\prime$}
 \end{align}
 instead. This leads us to the nemytskii mapping 
 \begin{align} 
   \mathcal{N}\sqbrac{\energyIntl[]}\tx:=\energyIntl[]\brac{t,x,\entropy\tx, \volumeSpecific\tx, \yl[1]\tx,\ldots, \yl[L]\tx}~,  
   \tag{${\text{\ref{eq:NonEqThermo-entropy-nemytskii}}}^\prime$}
 \end{align}
 and thus to the nonequilibrium Gibbs relation
 \begin{align}
     \rol[]\derm\mathcal{N}\sqbrac{\energyIntl[]} 
   = \rol[]\derm\energyIntl[] + \temp~\rol[]\derm\entropy - \pressHydrl[]~\rol[]\derm\volumeSpecific + \sum_l \chempotl~\rol[]\derm\yl ~.  
  \tag{${\text{\ref{eq:NonEqThermo-entropy-gibbsNonEq}}}^\prime$}
 \end{align}
 This shows, that various nonequilibrium Gibbs relations can be derived by the procedure \enquote{reasonable functional ansatz for~$\energyIntl[]$~}~$\rightarrow$~Nemytskii mapping~$\rightarrow$~Gibbs relation.
\hfill$\square$ 
\end{remark}
\begin{remark}[Chemical potential]\label{remark:NonEqThermo-chemicalPotential}
 From \eqref{eq:NonEqThermo-entropy-energyIntPure}--\eqref{eq:NonEqThermo-entropy-energyInt} and \eqref{eq:NonEqThermo-entropy-chemPot} we deduce for the chemical potentials~$\chempotl$ 
 \begin{align*}
   \chempotl := \derd[\yl]\energyIntl[] 
              = \derd[\yl]\energyIntPurel[] + \derd[\yl]\energyIntMixl[]  
              = \energyIntPurel\brac{\entropy, \volumeSpecific} + \frac{1}{2}\abs{\fieldF_l-\fieldF}^2+ \derd[\yl]\sqbrac{\yl\energyIntMixl}~.
 \end{align*}
 Thus, defining the so-called pure substance chemical potentials~$\chempotPurel$ and the chemical potentials of mixing~$\chempotMixl$ by
 \begin{align*}
   \chempotPurel :=\energyIntPurel\brac{\entropy, \volumeSpecific} + \frac{1}{2}\abs{\fieldF_l-\fieldF}^2 
   \qquad\text{and}\qquad 
   \chempotMixl := \derd[\yl]\sqbrac{\yl\energyIntMixl} \stackrel{\eqref{eq:NonEqThermo-energy-energyIntMix}}{=} \derd[\yl]\energyIntMixl[]~,
 \end{align*}
 we finally obtain for the chemical potentials~$\chempotl$ the decomposition 
 \begin{align}\label{eq:NonEqThermo-entropy-chemPotAnsatz}
   \chempotl = \energyIntPurel\brac{\entropy, \volumeSpecific} + \frac{1}{2}\abs{\fieldF_l-\fieldF}^2+ \derd[\yl]\sqbrac{\yl\energyIntMixl} = \chempotPurel + \chempotMixl~. \\[-14.0mm] \nonumber
 \end{align}
\hfill$\square$ 
\end{remark}
\begin{remark}[Electrochemical potential]\label{remark:NonEqThermo-electrochemicalPotential}
 The Gibbs relation~\eqref{eq:NonEqThermo-entropy-gibbsNonEq} does not contain the electric energy. However, as we consider mixtures of charged constituents, 
 it would be natural to involve the electric energy in the Gibbs relation. For that purpose, we introduce the so-called electrochemical potentials, cf. \cite{DeGrootMazur-book, Hunter-book, Israelachvili-book}, by 
 \begin{align}\label{eq:NonEqThermo-entropy-eletrochemPot}
   \elchempotl:=\chempotl + \frac{\constCharge\zl}{\ml}\potEL\qquad\sim\sqbrac{J \, kg^{-1}}.
 \end{align}
 With this definition, we obtain with the specific free charge~$\chargeELspecific:=\chargeEL/\rol[]\sim[C/kg]$ and \eqref{eq:NonEqThermo-mass-massBalanceTot}
 \begin{align*}
      \sum_l \chempotl~\rol[]\derm\yl 
  &= \sum_l \elchempotl~\rol[]\derm\yl - \sum_l \frac{\constCharge\zl}{\ml} \potEL ~\rol[]\derm\yl -\sum_l \frac{\constCharge\zl}{\ml} \potEL~ \overbrace{\sqbrac{\dert\rol[]+\grad\cdot(\rol[]\fieldF)}}^{=0}\yl\\
  &= \sum_l \elchempotl~\rol[]\derm\yl - \potEL \sum_l \frac{\constCharge\zl}{\ml} \sqbrac{\dert\rol + \grad\cdot(\rol\fieldF)}~.
 \end{align*}
 Furthermore, using again \eqref{eq:NonEqThermo-mass-massBalanceTot}, we have the identity $\dert\rol + \grad\cdot(\rol\fieldF) = \rol[]\derm\chargeELspecific$. This yields 
 \begin{align*}
  \sum_l \chempotl~\rol[]\derm\yl = \sum_l \elchempotl~\rol[]\derm\yl - \potEL~\rol[]\derm\chargeELspecific~.
 \end{align*}
 Hence, we equivalently rewrite the Gibbs relation for nonequilibrium systems~\eqref{eq:NonEqThermo-entropy-gibbsNonEq} with the electrochemical potentials~$\elchempotl$ as
 \begin{align}\label{eq:NonEqThermo-entropy-gibbsNonEq2}
    \rol[]\derm\energyIntl[] = \temp~\rol[]\derm\entropy - \pressHydrl[]~\rol[]\derm\volumeSpecific + \sum_l \elchempotl~\rol[]\derm\yl - \potEL~\rol[]\derm\chargeELspecific  ~\sim\sqbrac{Jm^{-3}s^{-1}}.
    \tag{\ref{eq:NonEqThermo-entropy-gibbsNonEq}a}
 \end{align}
 This version of Gibbs relation reveals the contributions of the electric energy. 
 Furthermore, inserting~\eqref{eq:NonEqThermo-entropy-chemPotAnsatz} in \eqref{eq:NonEqThermo-entropy-eletrochemPot}, 
 leads to~$\elchempotl=\chempotPurel + \chempotMixl + \frac{\constCharge\zl}{\ml}\potEL$. Thus, we can introduce the so-called electrochemical potentials of mixing~$\elchempotMixl$ by
 \begin{align}\label{eq:NonEqThermo-entropy-eletrochemPotMix}
   \elchempotMixl:=\chempotMixl + \frac{\constCharge\zl}{\ml}\potEL\qquad\sim\sqbrac{J (kg)^{-1}}. \\[-14.0mm] \nonumber
 \end{align}
\hfill$\square$ 
\end{remark}
\medskip
\par
We proceed with the introduction of the balance equation for the entropy density~$\rol[]\entropy\sim[J/(Km^3)]$. This equation is given with the entropy flux~$\entropyFlux\sim[J/(Km^2s)]$ 
and the entropy production rate~$\diss\sim[J/(Km^3s)]$ by  
\begin{align}\label{eq:NonEqThermo-entropy-entropyPDE}
 \dert\brac{\rol[]\entropy} +\grad\cdot\brac{\rol[]\entropy\fieldF+\entropyFlux} =\diss &\quad\text{in } \Omega \qquad\sqbrac{J K^{-1} m^{-3} s^{-1} }.
\end{align}
Based on this general balance equation, we formulate the \textbf{\textit{second law of thermodynamics}}, which states that \enquote{\textit{entropy can not be destroyed}},
cf. \cite{ColemanNoll, DeGrootMazur-book, Dreyer1987, Liu-book, Muller1968, MullerTomasso-book, oden-book, PrigogineKondepudi-book, Tadmor-book, Wilmanski-book}. 
Thus, mathematically the second law of thermodynamics shortly reads with the entropy production rate~$\diss$ as 
\begin{align}\label{eq:NonEqThermo-entropy-secondLaw}
 \diss\geq 0.
\end{align}
The second law of thermodynamics is precisely the missing tool, which contains the information about admissible directions of physical processes. More precisely, the evolution of every 
thermodynamic process must respect to the constraint~$\diss\geq0$. Hence, provided a process leads to entropy production, i.e., $\diss>0$, the entropy irreversibly increases,  
as entropy can not be destroyed. Consequently, this process never returns to its initial state. Such processes are called \textit{irreversible}, cf. \cite{EckGarckeKnabner-book, DreyerMullerWeiss, LiebYngvason1999, oden-book}.  
Note that \eqref{eq:NonEqThermo-entropy-entropyPDE} and \eqref{eq:NonEqThermo-entropy-secondLaw} contain the classical Clausius inequality, cf. \cite[p. 25]{DeGrootMazur-book}. 
However, regarding a detailed presentation of the classical results and the history of thermodynamics, we refer, e.g.,  to \cite{EckGarckeKnabner-book, Liu-book, Muller1968, oden-book, PrigogineKondepudi-book, Tadmor-book, TruesdellToupin-book, Wilmanski-book}. 
\par
We proceed by deriving an explicit expression for the entropy production rate~$\diss$. For that purpose, we rearrange Gibbs relation~\eqref{eq:NonEqThermo-entropy-gibbsNonEq}. Thereby, we come to
\begin{align*}
  \rol[]\derm\entropy = \frac{1}{\temp}~\rol[]\derm\energyIntl[] + \frac{\pressHydrl[]}{\temp}~\rol[]\derm\volumeSpecific - \sum_l \frac{\chempotl}{\temp}~\rol[]\derm\yl~.
\end{align*}
Inserting equations~\eqref{eq:NonEqThermo-mass-massFracBalance-NonCons}, \eqref{eq:NonEqThermo-mass-specificVolPDE}, \eqref{eq:NonEqThermo-energy-internalEnergyPDE2}, 
and substituting the mixture stress tensor~$\stressTotl[]$ by means of \ref{Assump:mixtureStress2}, yields
\begin{align*}
  \rol[]\derm\entropy &=  \frac{1}{\temp}~\sqbrac{-\grad\cdot\vecq + \currentEL\cdot\fieldEL + \stressViscl[]:\grad\fieldF}  
                          - \sum_l \frac{\chempotl}{\temp}~\sqbrac{- \grad\cdot\rolfluxrel + \rl}~.
\end{align*}
Furthermore, we transform parts of the right-hand side with $\grad\,\temp^{-1}=-\temp^{-2}\grad\,\temp$ and \ref{Assump:electrostatics}, \eqref{eq:NonEqThermo-charge-DefFreeCurrent}, 
\eqref{eq:NonEqThermo-entropy-eletrochemPot} to
\begin{align*}
   \sum_l \frac{\chempotl}{\temp}\grad\cdot\rolfluxrel + \frac{1}{\temp}\currentEL\cdot\fieldEL
 = \grad\cdot\brac{\sum_l \frac{\elchempotl}{\temp}\rolfluxrel - \frac{\potEL}{\temp}\currentEL }  - \sum_l \grad\brac{\frac{\elchempotl}{\temp}}\cdot\rolfluxrel - \frac{\potEL}{\temp^2}~\currentEL\cdot\grad\,\temp~,
\end{align*}
Analogously, we treat the term $\temp^{-1}~\grad\cdot\vecq$ on the right-hand side. Thereby, we arrive at
\begin{align*}
  \rol[]\derm\entropy &=     \grad\cdot\brac{-\frac{1}{\temp}\vecq + \sum_l \frac{\elchempotl}{\temp}\rolfluxrel - \frac{\potEL}{\temp}\currentEL } +  \frac{1}{\temp}\stressViscl[]:\grad\fieldF \\  
                      &~~~~  -\frac{1}{\temp^2}\grad\,\temp\cdot\brac{\vecq+\potEL\currentEL} -\sum_l\grad\brac{\frac{\elchempotl}{\temp}}\cdot\rolfluxrel - \sum_l \frac{\chempotl}{\temp}\rl~.
\end{align*}
Moreover, with \eqref{eq:NonEqThermo-mass-massBalanceTot} we have $\rol[]\derm\entropy = \dert\brac{\rol[]\entropy} + \grad\cdot\brac{\rol[]\entropy\fieldF}$ on the left-hand side, 
and on the right-hand side we take with \ref{Assump:reactionsAndCharges} the identity
\begin{align}\label{eq:NonEqThermo-entropy-reactionCriterion}
 \sum_l \frac{\chempotl}{\temp}\rl = \sum_l \frac{\elchempotl}{\temp}\rl - \frac{\potEL}{\temp}\sum_l \frac{\constCharge\zl}{\ml}\rl = \sum_l \frac{\elchempotl}{\temp}\rl
\end{align}
into account. Finally, this yields 
\begin{align}\label{eq:NonEqThermo-entropy-entropyPDE1}
  &  \dert\brac{\rol[]\entropy} + \grad\cdot\brac{ \rol[]\entropy\fieldF + \frac{1}{\temp}\vecq - \sum_l \frac{\elchempotl}{\temp}\rolfluxrel +\frac{\potEL}{\temp}\currentEL  } \nonumber\\
  &=     -\frac{1}{\temp^2}\grad\,\temp\cdot\brac{\vecq+\potEL\currentEL} + \frac{1}{\temp}\stressViscl[]:\grad\fieldF -\sum_l\grad\brac{\frac{\elchempotl}{\temp}}\cdot\rolfluxrel - \sum_l \frac{\elchempotl}{\temp}\rl~.
\end{align}
This equation is exactly the desired \enquote{electrochemical version} of the entropy evolution equation, which introduces explicit expressions for the entropy flux~$\entropyFlux$ and the entropy 
production rate~$\diss$. More precisely, comparing this equation with equation~\eqref{eq:NonEqThermo-entropy-entropyPDE} uncovers for the entropy flux the definition
\begin{subequations}
\begin{align}\label{eq:NonEqThermo-entropy-entropyFlux}
 ~\entropyFlux &:=  \frac{1}{\temp}\vecq - \sum_l \frac{\elchempotl}{\temp}\rolfluxrel +\frac{\potEL}{\temp}\currentEL~, 
\end{align}
and for the entropy production rate the definition
\begin{flalign}\label{eq:NonEqThermo-entropy-diss1}
 ~~  \underbrace{\diss}_{ \begin{subarray}{c} \text{total} \\ \text{entropy} \\ \text{production} \end{subarray} }       
 &= -\underbrace{\frac{1}{\temp^2}\grad\,\temp\cdot\brac{\vecq+\potEL\currentEL}}_{\begin{subarray}{c}  \text{electrothermal} \\ \text{part} \end{subarray}}
    +\underbrace{\frac{1}{\temp}\stressViscl[]:\grad\fieldF}_{\begin{subarray}{c}  \text{viscous} \\ \text{part}  \end{subarray}} 
    -\underbrace{\sum_l\grad\brac{\frac{\elchempotl}{\temp}}\cdot\rolfluxrel}_{\begin{subarray}{c}  \text{thermo-mixing} \\ \text{part} \end{subarray}}
    -\underbrace{\sum_l \frac{\elchempotl}{\temp}\rl~.}_{\begin{subarray}{c}  \text{electrochemical} \\ \text{part} \end{subarray}}&
\end{flalign}
\end{subequations}
This equation is of extreme importance, since knowing the precise sources of the entropy production rate reveals which processes are at the heart of irreversibility.
Moreover, according to the previous formulation of the second law of thermodynamics, we have $\diss\geq0$. Thus, when substituting constitutive laws for $\vecq$, $\stressViscl[]$, $\chempotl$, 
$\rolfluxrel$, and $\rl$ into \eqref{eq:NonEqThermo-entropy-diss1}, these constitutive laws must respect~$\diss\geq0$. Hence, $\diss\geq0$ and \eqref{eq:NonEqThermo-entropy-diss1} 
restrict the admissible choices of constitutive laws. Thus, we now have a useful criterion, which validates, whether a constitutive law respects the second law of thermodynamics. 
Next, we note that with the identity  
\begin{align*}
    \sum_l\grad\brac{\frac{\elchempotl}{\temp}}\cdot\rolfluxrel
  = -\sum_l\frac{\elchempotl}{\temp^2}\grad\,\temp\cdot\rolfluxrel + \sum_l\frac{1}{\temp}\grad\elchempotl\cdot\rolfluxrel
\end{align*}
and \eqref{eq:NonEqThermo-entropy-entropyFlux}, we rewrite the entropy production rate~$\diss$ from \eqref{eq:NonEqThermo-entropy-diss1} as
\begin{flalign}\label{eq:NonEqThermo-entropy-diss2}
 &~~~~ \underbrace{\diss}_{ \begin{subarray}{c} \text{total} \\ \text{entropy} \\ \text{production} \end{subarray} }       
 = -\underbrace{\frac{1}{\temp}\grad\,\temp\cdot\entropyFlux}_{\begin{subarray}{c}  \text{entropic flux} \\ \text{part} \end{subarray}}
   +\underbrace{\frac{1}{\temp}\stressViscl[]:\grad\fieldF}_{\begin{subarray}{c}  \text{viscous} \\ \text{part}  \end{subarray}}   
   -\underbrace{\frac{1}{\temp}\sum_l\grad\elchempotl\cdot\rolfluxrel}_{\begin{subarray}{c}  \text{mixing} \\ \text{part} \end{subarray}}
   -\underbrace{\sum_l \frac{\elchempotl}{\temp}\rl~.}_{\begin{subarray}{c}  \text{electrochemical} \\ \text{part} \end{subarray}}&    \tag{${\text{\ref{eq:NonEqThermo-entropy-diss1}}}^\prime$}
\end{flalign}
The first term of the right-hand side of this equation uncovers the remarkable fact, that flow of entropy can produce entropy. 
However, as $\grad\,\temp$ is perpendicular to isotherms, this does not occur when entropy solely flows along isotherms%
\footnote{Since isotherms are temperature contour lines, we have $\grad\,\temp\perp\entropyFlux$, which leads to $\grad\,\temp\cdot\entropyFlux=0$.}.
In this case, the temperature~$\temp$ can be considered as first integral for the entropy flow as the Lie~derivative~$\mathcal{L}_{\entropyFlux}\temp:=\grad\,\temp\cdot\entropyFlux$ vanishes. 
Thus, in particular in isothermal situations entropy flow never lead to entropy production, and generally, we deduce from the preceding equation the criterion:
\begin{align*}
  \text{entropy flow causes entropy production} \qquad\Equivalent\qquad \mathcal{L}_{\entropyFlux}\temp<0~.  
\end{align*}
\par
Next, we note that in situations without barycentric flow, without viscous effects, without reactions, and without electrics, the preceding equation for the entropy production rate reduces to 
\begin{flalign*}
 &~~ \underbrace{\diss}_{ \begin{subarray}{c} \text{total} \\ \text{entropy} \\ \text{production} \end{subarray} }       
 = -\underbrace{\frac{1}{\temp^2}\grad\,\temp\cdot\vecq~.}_{\begin{subarray}{c}  \text{heat flux} \\ \text{part} \end{subarray}}&
\end{flalign*}
Thus, we deduce with $\diss\geq0$ that the heat flux~$\vecq$ must point into the direction of the negative temperature gradient, i.e., heat must flow down the temperature gradient. 
This is exactly the mathematical formulation of the classical statement \enquote{\textit{heat must flow from hot to cold}} of the second law of thermodynamics, 
cf. \cite{DreyerMullerWeiss, oden-book, PrigogineKondepudi-book, TruesdellToupin-book}. 
\begin{remark}[Equivalent formulation of the entropy evolution, $\entropyFlux$, and $\diss$]\label{remark:NonEqThermo-entropyPDE}
 We note, that carefully reading through the above derivation of the \enquote{electrochemical} entropy evolution equation~\eqref{eq:NonEqThermo-entropy-entropyPDE1} shows, that this 
 equation is equivalent to the entropy evolution equation
 \begin{align}\label{eq:NonEqThermo-entropy-entropyPDE1a}
  &  \dert\brac{\rol[]\entropy} + \grad\cdot\brac{ \rol[]\entropy\fieldF + \frac{1}{\temp}\vecq - \sum_l \frac{\chempotl}{\temp}\rolfluxrel } \nonumber\\
  &=     -\frac{1}{\temp^2}\grad\,\temp\cdot\vecq + \frac{1}{\temp}\stressViscl[]:\grad\fieldF  
         -\sum_l\grad\brac{\frac{\chempotl}{\temp}}\cdot\rolfluxrel - \sum_l \frac{\chempotl}{\temp}\rl~.   \tag{${\text{\ref{eq:NonEqThermo-entropy-entropyPDE1}}}^\prime$}
 \end{align}
 Here, the entropy flux~$\entropyFlux$ is given by 
 \begin{subequations}
 \begin{align}
   \entropyFlux &:= \frac{1}{\temp}\vecq - \sum_l \frac{\chempotl}{\temp}\rolfluxrel~,        \label{eq:NonEqThermo-entropy-entropyFlux1a}  \tag{${\text{\ref{eq:NonEqThermo-entropy-entropyFlux}}}^\prime$}
 \end{align}
 and the entropy production rate~$\diss$ now reads as
 \begin{align}
   \diss   := -\frac{1}{\temp^2}\grad\,\temp\cdot\vecq + \frac{1}{\temp} \stressViscl[]:\grad\fieldF  
               -\sum_l\grad\brac{\frac{\chempotl}{\temp}}\cdot\rolfluxrel  - \sum_l \frac{\chempotl}{\temp}\rl~.  \label{eq:NonEqThermo-entropy-diss1a} \tag{${\text{\ref{eq:NonEqThermo-entropy-diss1}}}^{\prime\prime}$}  \\[-12.0mm] \nonumber
 \end{align}
 \end{subequations}
\hfill$\square$
\end{remark}  
\medskip
\par
So far, we considered the evolution of the total specific entropy~$\entropy$. However, in \ref{Assump:energy1a} and \ref{Assump:energy2a} we obtained a more detailed picture for the 
total specific internal energy~$\energyIntl[]$, as we assumed a decomposition into a pure substance part~$\energyIntPurel[]$ and a part due to mixing~$\energyIntMixl[]$. 
In \cref{remark:NonEqThermo-chemicalPotential}, we showed that this leads to corresponding decompositions of the chemical potentials~$\chempotl$ 
into pure substance parts~$\chempotPurel$ and parts due to mixing~$\chempotMixl$. We now assume, that a analogous decomposition holds true for the entropy~$\entropy$
\begin{enumerate}[align=left, leftmargin=*, topsep=2.0mm, label={(\labelA\arabic*)}, start=\value{countModAssump}+1]
  \item \label{Assump:entropy1}%
        \textbf{Structure of the entropy: } We suppose, the specific entropy~$\entropy~[J/(K\,kg)]$ additively decomposes into a pure substance part~$\entropyPure$ and an entropy of mixing~$\entropyMixl[]$, i.e.,
	 \begin{align}\label{eq:NonEqThermo-entropy-ansatzEntropyMixPure}
	   \entropy = \entropyPure + \entropyMixl[]~.
	 \end{align}
  \item \label{Assump:entropy1a}%
        \textbf{Structure of the specific entropy of mixing: } In particular for the specific entropy of mixing~$\entropyMixl[]$, we assume the ansatz   
	 \begin{align}\label{eq:NonEqThermo-entropy-ansatzMixDetailed}
	   \entropyMixl[]\brac{\yl[1],\ldots,\yl[L]} =\sum_l\yl\entropyMixl\brac{\yl}, \quad\text{with}~~ -\temp\entropyMixl\brac{\yl}:=\energyIntMixl\brac{\yl}~.
	 \end{align}
\setcounter{countModAssump}{\value{enumi}}%
\end{enumerate}
Regarding these assumptions, we add the following explanations and comments: 
\begin{enumerate}[align=left, label=(\roman*), leftmargin=*, topsep=2.0mm, itemsep=-1.2mm]
  \item Firstly, the ansatz $-\temp\entropyMixl=\energyIntMixl$ from \eqref{eq:NonEqThermo-entropy-ansatzMixDetailed} is well-known from mixtures of ideal gases, 
        cf. \cite{EckGarckeKnabner-book, PrigogineKondepudi-book, Atkins-book}. Thus, by adopting this relation, we assume that concerning the phenomena of mixing, 
        the considered mixtures behave as mixtures of ideal gases. 
  \item Secondly, \eqref{eq:NonEqThermo-entropy-ansatzMixDetailed} implies for the corresponding mixture variables
	\begin{align}\label{eq:NonEqThermo-entropy-mixEntropyAndEnergy}
	 -\temp\entropyMixl[] = -\temp\sum_l\yl\entropyMixl = \sum_l \yl\energyIntMixl = \energyIntMixl[]~. 
	\end{align}
  \item Thirdly, due to assumption~\ref{Assump:entropy1a}, we rigorously must distinguish between the Nemytskii~mapping~$\mathcal{N}[\entropyMixl]\tx$ and $\entropyMixl[]\brac{\yl[1],\ldots,\yl[L]}$. 
        This is analogous to \eqref{eq:NonEqThermo-entropy-nemytskii}. However, for ease of readability we henceforth omit this difference in notation. 
  \item Fourthly, from \eqref{eq:NonEqThermo-entropy-mixEntropyAndEnergy} and the definition of the chemical potentials of mixing~$\chempotMixl$ in \cref{remark:NonEqThermo-chemicalPotential}, we obtain 
	\begin{align}\label{eq:NonEqThermo-entropy-chemPotMixAndEntropy}
	\frac{\chempotMixl}{\temp}= \frac{\derd[\yl]\energyIntMixl[]}{\temp} = \frac{\derd[\yl](-\temp\entropyMixl[])}{\temp} = -\derd[\yl]\entropyMixl[]~.
	\end{align}
        From this identity, we furthermore deduce for the Nemytskii-mapping $\entropyMixl[]$
       \begin{align}\label{eq:NonEqThermo-entropy-matDerEntropyMix}
         \derm\entropyMixl[]\tx &= \derm\entropyMixl[]\brac{\yl[1]\tx,\ldots,\yl[L]\tx} \nonumber\\
                                &= \sum_l \derd[\yl]\entropyMixl\brac{\yl[1],\ldots,\yl[L]}\derm\yl\tx~.
       \end{align}
\end{enumerate}
Equipped with the above assumptions and relations, we now derive an evolution equation for the entropy of mixing~$\entropyMixl[]$. To derive this equation, we multiply 
equations~\eqref{eq:NonEqThermo-mass-massFracBalance-NonCons} by $-\temp^{-1}\chempotMixl$. This yields
\begin{align*}
   -\rol[]\,\frac{\chempotMixl}{\temp}\,\derm\yl  
 = \frac{\chempotMixl}{\temp}\grad\cdot\rolfluxrel - \frac{\chempotMixl}{\temp}\rl 
 = \grad\cdot\brac{\frac{\chempotMixl}{\temp}\rolfluxrel} - \grad\brac{\frac{\chempotMixl}{\temp}}\cdot\rolfluxrel - \frac{\chempotMixl}{\temp}\rl~.
\end{align*}
Summing over $l\in\cbrac{1,\ldots,L}$, involving \eqref{eq:NonEqThermo-entropy-chemPotMixAndEntropy}, \eqref{eq:NonEqThermo-entropy-matDerEntropyMix}, and applying on the left-hand side 
with \eqref{eq:NonEqThermo-mass-massBalanceTot} the identity $\rol[]\derm\entropyMixl[] = \dert\brac{\rol[]\entropyMixl[]} + \grad\cdot\brac{\rol[]\entropyMixl[]\fieldF}$, finally results in the desired evolution 
equation for the entropy of mixing in conservative from
\begin{flalign}\label{eq:NonEqThermo-entropy-entropyMixPDE0}
 &~~\dert\brac{\rol[]\entropyMixl[]}  + \grad\cdot\brac{\rol[]\entropyMixl[]\fieldF - \sum_l\frac{\chempotMixl}{\temp}\rolfluxrel} = -\sum_l\grad\brac{\frac{\chempotMixl}{\temp}}\cdot\rolfluxrel - \sum_l\frac{\chempotMixl}{\temp}\rl~.&
\end{flalign}
However, as we consider mixtures of charged constituents, we now replace the chemical potentials of mixing~$\chempotMixl$ by the electrochemical potentials of mixing~$\elchempotMixl$, which are 
defined in \cref{remark:NonEqThermo-electrochemicalPotential}. Thereby, we obtain for the last term on the right-hand side with \ref{Assump:reactionsAndCharges}
\begin{align*}
 \sum_l\frac{\chempotMixl}{\temp}\rl = \sum_l\frac{\elchempotMixl}{\temp}\rl~,
\end{align*}
and the parts including the drift \massfluxRol{es}~$\rolfluxrel$, we transform with \ref{Assump:electrostatics} to 
\begin{align*}
 &  \grad\cdot\brac{\sum_l\frac{\chempotMixl}{\temp}\rolfluxrel} -\sum_l\grad\brac{\frac{\chempotMixl}{\temp}}\cdot\rolfluxrel \\
 &= \grad\cdot\brac{\sum_l\frac{\elchempotMixl}{\temp}\rolfluxrel - \potEL\currentEL} - \sum_l\grad\brac{\frac{\elchempotMixl}{\temp}}\cdot\rolfluxrel - \frac{1}{\temp}\fieldEL\cdot\currentEL - \frac{\potEL}{T^2}\grad\,\temp\cdot\currentEL~.  
\end{align*}
Thus, we finally arrive at the \enquote{electrochemical} counterpart of the evolution equation~\eqref{eq:NonEqThermo-entropy-entropyMixPDE0} 
\begin{align}\label{eq:NonEqThermo-entropy-entropyMixPDE}
   &  \dert\brac{\rol[]\entropyMixl[]} + \grad\cdot\brac{\rol[]\entropyMixl[]\fieldF - \sum_l\frac{\elchempotMixl}{\temp}\rolfluxrel + \potEL\currentEL}  \nonumber \\
   &= - \frac{1}{T^2}\grad\,\temp\cdot\brac{\potEL\currentEL} - \sum_l\grad\brac{\frac{\elchempotMixl}{\temp}}\cdot\rolfluxrel - \sum_l\frac{\elchempotMixl}{\temp}\rl - \frac{1}{\temp}\fieldEL\cdot\currentEL~.
\end{align}
This equation allows us to define the entropy of mixing flux~$\entropyMixFlux$ and the entropy of mixing production rate~$\dissMix$ by
\begin{subequations}
\begin{flalign}
 ~\entropyMixFlux &:= - \sum_l \frac{\elchempotMixl}{\temp}\rolfluxrel +\frac{\potEL}{\temp}\currentEL~,& \label{eq:NonEqThermo-entropy-entropyMixFlux}\\
 ~\dissMix        &:= -\frac{1}{\temp^2}\grad\,\temp\cdot\brac{\potEL\currentEL} -\sum_l\grad\brac{\frac{\elchempotMixl}{\temp}}\cdot\rolfluxrel  
                      - \sum_l \frac{\elchempotMixl}{\temp}\rl - \frac{1}{\temp}\fieldEL\cdot\currentEL ~.                          & \label{eq:NonEqThermo-entropy-dissMix}
\end{flalign}
\end{subequations} 
\begin{remark}[Equivalent formulation of $\entropyMixFlux$ and $\dissMix$]\label{remark:NonEqThermo-dissMix}
 We note, that equation~\eqref{eq:NonEqThermo-entropy-entropyMixPDE0} shows that instead of the preceding \enquote{electrochemical} versions of the entropy flux~$\entropyMixFlux$ 
 from \eqref{eq:NonEqThermo-entropy-entropyMixFlux} and the entropy production rate~$\dissMix$ from \eqref{eq:NonEqThermo-entropy-dissMix}, 
 we can alternatively define these quantities by
 \begin{subequations}
 \begin{flalign*}
   &~\entropyMixFlux :=  - \sum_l \frac{\chempotMixl}{\temp}\rolfluxrel
     ~~~~~\text{and}~~~~~
    \dissMix := -\sum_l\grad\brac{\frac{\chempotMixl}{\temp}}\cdot\rolfluxrel - \sum_l \frac{\chempotMixl}{\temp}\rl~.   & \\[-16.mm] \nonumber 
 \end{flalign*}
 \end{subequations}
\hfill$\square$
\end{remark}
Next, we subtract equation~\eqref{eq:NonEqThermo-entropy-entropyMixPDE0} from equation~\eqref{eq:NonEqThermo-entropy-entropyPDE1a}. Together with the additive splitting $\elchempotl=\chempotMixl+\chempotPurel$ 
from \cref{remark:NonEqThermo-electrochemicalPotential} and $\entropyPure=\entropy-\entropyMixl[]$ from \ref{Assump:entropy1}, we thereby obtain the evolution equation for $\entropyPure$
\begin{flalign}\label{eq:NonEqThermo-entropy-entropyPurePDE}
  &~~\dert\brac{\rol[]\entropyPure} + \grad\cdot\brac{ \rol[]\entropyPure\fieldF + \frac{1}{\temp}\vecq - \sum_l \frac{\chempotPurel}{\temp}\rolfluxrel } \nonumber\\
  &~~=     -\frac{1}{\temp^2}\grad\,\temp\cdot\vecq  + \frac{1}{\temp}\stressViscl[]:\grad\fieldF    
         -\sum_l\grad\brac{\frac{\chempotPurel}{\temp}}\cdot\rolfluxrel - \sum_l \frac{\chempotPurel}{\temp}\rl \,. &
\end{flalign}
Hence, the pure substance entropy flux~$\entropyPureFlux$ and the pure substance entropy production rate~$\dissPure$ are given by
\begin{subequations}
\begin{flalign}
 ~\entropyPureFlux &:=  \frac{1}{\temp}\vecq - \sum_l \frac{\chempotPurel}{\temp}\rolfluxrel  ~,          & \label{eq:NonEqThermo-entropy-entropyPureFlux}\\
 ~\dissPure        &:= -\frac{1}{\temp^2}\grad\,\temp\cdot\vecq + \frac{1}{\temp}\stressViscl[]:\grad\fieldF 
                       -\sum_l\grad\brac{\frac{\chempotPurel}{\temp}}\!\cdot\rolfluxrel  - \! \sum_l \frac{\chempotPurel}{\temp}\rl \,.\hspace{-2.0mm}  & \label{eq:NonEqThermo-entropy-dissPure}
\end{flalign}
\end{subequations}
\par
Finally, the \textbf{\textit{second law of thermodynamics}} now reads with $\diss$ from \eqref{eq:NonEqThermo-entropy-diss1}, 
$\dissMix$ from \eqref{eq:NonEqThermo-entropy-dissMix}, and $\dissPure$ from \eqref{eq:NonEqThermo-entropy-dissPure} as
\begin{align}\label{eq:NonEqThermo-entropy-secondLawPureMix}
  \dissPure + \dissMix + \frac{1}{\temp}\fieldEL\cdot\currentEL = \diss \geq 0~.
\end{align}
\begin{remark}[Equivalent formulation of the second law of thermodynamics]\label{remark:NonEqThermo-secondLawMixPure}
 The preceding inequality~\eqref{eq:NonEqThermo-entropy-secondLawPureMix} is the \enquote{electrochemical} formulation of the second law of thermodynamics. 
 Returning to \cref{remark:NonEqThermo-entropyPDE}  and \cref{remark:NonEqThermo-dissMix} shows, that the above statement of the second law of thermodynamics is 
 equivalent to 
 \begin{align*}
   \dissMix + \dissPure = \diss \geq 0~.
 \end{align*}
 Here, we keep the formula~\eqref{eq:NonEqThermo-entropy-dissPure} for $\dissPure$, but we alternatively use for $\diss$ the formula from \cref{remark:NonEqThermo-entropyPDE} and for $\dissMix$ 
 the formula from \cref{remark:NonEqThermo-dissMix}. This means,  $\dissMix$ and $\diss$ are given with $\chempotl=\chempotMixl+\chempotPurel$ from \cref{remark:NonEqThermo-chemicalPotential} by 
 \begin{align*}
  &~\dissMix  = -\sum_l\grad\brac{\frac{\chempotMixl}{\temp}}\cdot\rolfluxrel - \sum_l \frac{\chempotMixl}{\temp}\rl \, ,    \\
  &~\diss     = -\frac{1}{\temp^2}\grad\,\temp\cdot\vecq + \frac{1}{\temp}\stressViscl[]:\grad\fieldF 
                 -\sum_l\grad\brac{\frac{\chempotl}{\temp}}\cdot\rolfluxrel  - \sum_l \frac{\chempotl}{\temp}\rl \, .  \\[-12.0mm] \nonumber
 \end{align*} 
\hfill$\square$
\end{remark}
\begin{remark}[Decomposition of the internal energy equation]\label{remark:NonEqThermo-splittingInternalEnergyPDE}
 Finally, we note that we can deduce an equation for the specific internal energy of mixing~$\energyIntMixl[]$ by multiplying equations~\eqref{eq:NonEqThermo-mass-massFracBalance-NonCons} 
 by the chemical potentials~$\chempotMixl$. In exactly in the same manner as we derived equation~\eqref{eq:NonEqThermo-entropy-entropyMixPDE0}, we thereby arrive at the 
 following evolution equation for~$\rol[]\energyIntMixl[]$ 
 \begin{align}\label{eq:NonEqThermo-energy-energyIntMixPDE0}
   \dert\brac{\rol[]\energyIntMixl[]}  + \grad\cdot\brac{\rol[]\energyIntMixl[]\fieldF + \sum_l\chempotMixl\rolfluxrel} = \sum_l\grad\chempotMixl\cdot\rolfluxrel + \sum_l\chempotMixl\rl~.
 \end{align}
 Moreover, by adopting the produce, which gave us equation~\eqref{eq:NonEqThermo-entropy-entropyMixPDE}, we obtain the \enquote{electrochemical version} of the evolution equation 
 \begin{align}\label{eq:NonEqThermo-energy-energyIntMixPDE}
   &  \dert\brac{\rol[]\energyIntMixl[]} + \grad\cdot\brac{\rol[]\energyIntMixl[]\fieldF + \sum_l\elchempotMixl\rolfluxrel - \potEL\currentEL}  \nonumber \\
   &= \sum_l\grad\elchempotMixl\cdot\rolfluxrel + \sum_l\elchempotMixl\rl + \fieldEL\cdot\currentEL~.
 \end{align}
 Analogous to \eqref{eq:NonEqThermo-entropy-entropyPurePDE}, we obtain the evolution equation for the specific pure substance internal energy~$\energyIntPurel[]=\energyIntl[]-\energyIntMixl[]$ 
 by subtracting equation~\eqref{eq:NonEqThermo-energy-energyIntMixPDE} from equation~\eqref{eq:NonEqThermo-energy-internalEnergyPDE}. This yields 
\begin{align}\label{eq:NonEqThermo-energy-energyIntPurePDE}
   &\dert\brac{\rol[]\energyIntPurel[]} + \grad\cdot\brac{ \rol[]\energyIntPurel[]\fieldF + \vecq - \sum_l\elchempotMixl\rolfluxrel + \potEL\currentEL } \nonumber\\
   &= \stressTotl[]:\grad\fieldF -\sum_l\grad\elchempotMixl\cdot\rolfluxrel - \sum_l\elchempotMixl\rl ~. 
\end{align}
\end{remark}
\fancyhead[L]{Conclusion}
\section{Conclusion}
In \cref{chapter:NonEqThermo} of this paper, we established the fundamental equations, that govern the evolution of mixtures of charged constituents. First of all, in \cref{sec:NonEqThermo-mass} 
we started with the equations for mass conservation. This section repeated the succinct presentation of \cite{DeGrootMazur-book} more detailed. 
Next, in \cref{sec:NonEqThermo-charge}, we derived the charge conservation equation and in \cref{sec:NonEqThermo-mom}, 
we set up the conservation equations for the barycentric momentum density following \cite{TruesdellToupin-book}. 
In \cref{sec:NonEqThermo-energy}, we adopted the presentation of \cite{DeGrootMazur-book} for the presentation of the first law of thermodynamics and for the derivation of the evolution equations 
for the electric potential energy, the kinetic energy, and the internal energy. Finally, in \cref{sec:NonEqThermo-entropy}, we again followed 
the approach of \cite{DeGrootMazur-book} to derive the evolution equation for the entropy density and to establish an explicit expression for the entropy flux 
and the entropy production rate. 
\par
The contribution of \cref{chapter:NonEqThermo} of this paper was to provide generalized nonequilibrium version of Dalton's law resp. Raoult's law for mixtures, cf. \cref{remark:NonEqThermo-dalton}. 
Moreover, we presented a more detailed picture for the internal energy and the entropy in terms of their pure substance parts and their parts due to mixing. In summary, we demonstrated that 
assumptions~\ref{Assump:entropy1} and \ref{Assump:entropy1a} not only additively decomposed the specific entropy~$\entropy$ into a pure substance part~$\entropyPure$ and a part due to mixing~$\entropyMixl[]$, 
but rather led to an additive decomposition of the evolution equation~\eqref{eq:NonEqThermo-entropy-entropyPDE1} into a pure substance part~\eqref{eq:NonEqThermo-entropy-entropyPurePDE} 
and a part due to mixing~\eqref{eq:NonEqThermo-entropy-entropyMixPDE}. The crucial point in this connection was, that we rigorously proved this decomposition of the evolution 
equation~\eqref{eq:NonEqThermo-entropy-entropyPDE1} by explicitly deriving equation~\eqref{eq:NonEqThermo-entropy-entropyMixPDE}. In \cref{remark:NonEqThermo-splittingInternalEnergyPDE}, we proved the 
same result for the specific internal energy~$\energyIntl[]$. An other essential observation concerning this decomposition of the 
entropy evolution equation was the absence of a common coupling term in equations \eqref{eq:NonEqThermo-entropy-entropyPurePDE} and \eqref{eq:NonEqThermo-entropy-entropyMixPDE}. 
This revealed, that the decomposition~$\entropy=\entropyPure+\entropyMixl[]$ from \ref{Assump:entropy1}, resulted in two \textit{decoupled} 
subprocesses, which were respectively governed by equations~\eqref{eq:NonEqThermo-entropy-entropyMixPDE} and \eqref{eq:NonEqThermo-entropy-entropyPurePDE}. Furthermore, we even obtained explicit expressions for 
the respective entropy production rates $\dissMix$ and $\dissPure$. For both parts these expressions uncovered the sources of irreversibility.
\par
Therefore, in \cref{chapter:NonEqThermo} of this paper we provided an approach, that allows to study the specific entropy of mixing~$\entropyMixl[]$ independently of the specific pure substance entropy~$\entropyPure$ 
and vice versa. This is of great value, in particular for a better understanding of the complex phenomena due to mixing in multicomponent systems.
%
%
\part{A Thermodynamical consistent Model for Electrolyte Solutions}\label{chapter:modelES}
\fancyhead[L]{Maxwell's equations}
\section{Electrostatic Limit of Maxwell's Equations}\label{sec:modelES-maxwell}
We start with the macroscopic Maxwell's equations
for linear materials.   
\begin{align*}
 & \grad\cdot \fieldMG = 0 && [\text{nonexistence of monopoles}],\\
 & \grad\cdot(\eps_r\fieldEL) = \frac{1}{\constEL}~\chargeEL && [\text{Gauss's law}], \\
 & \grad\times\fieldEL = - \dert \fieldMG && [\text{Faraday's law}],\\
 &  \grad\times\!\brac{\mu_r^{-1}\fieldMG} = \constMG~\currentEL + \frac{1}{c_0^2} ~\dert(\eps_r\fieldEL)&& [\text{Amp\`{e}re's law}].
\end{align*}
Here, $\constEL$ resp. $\constMG$ are the vacuum permittivity resp. the vacuum permeability, and $\eps_r$ resp. $\mu_r$ are the relative electric permittivity resp. the relative magnetic permeability of the medium.%
\footnote{While $\constEL$, $\constMG$ are scalar constants, $\eps_r$, $\mu_r$ may be tensors or even tensor valued functions $\eps_r\tx$, $\mu_r\tx$.} 
For a detailed derivation of these equations, we refer to \cite{LifshitzLandau-book2, LopezDavalos-book}. We note, that Maxwell's equations are relativistic equations. 
Subsequently, we derive their nonrelativistic limit. For that purpose, we introduce the nondimensional and rescaled fields~$\fieldEL^\ast$, $\fieldMG^\ast$, $\currentEL^\ast$, $\chargeEL^\ast$ by
\begin{align*}
 & E_0 \fieldEL^{\ast}(s,y):=E_0 \fieldEL^{\ast}\!\brac{ \frac{t}{\tau} , \frac{x}{l} }:=\fieldEL\tx \qquad \text{with a characteristic quantity }E_0~\sim[N/C],  \\
 & B_0 \fieldMG^{\ast}(s,y):=B_0 \fieldMG^{\ast}\!\brac{ \frac{t}{\tau} , \frac{x}{l} }:=\fieldMG\tx,\qquad \text{with a characteristic quantity } B_0~\sim[T],    \\
 & i_0 \currentEL^{\ast}(s,y):=i_0 \currentEL^{\ast}\!\brac{ \frac{t}{\tau} , \frac{x}{l} } :=\currentEL\tx  \qquad \text{with a characteristic quantity } i_0~\sim[A/m^2],   \\
 & \rho_0 \chargeEL^{\ast}(s,y):=\rho_0 \chargeEL^{\ast}\!\brac{ \frac{t}{\tau} , \frac{x}{l} } := \chargeEL\tx \qquad \text{with a characteristic quantity } \rho_0~\sim[C/m^3].
\end{align*}
We substitute these nondimensional and rescaled fields into the above Maxwell's equations for linear materials. This leads us to 
\begin{align*}
  & \grad_y\cdot \fieldMG^{\ast} = 0 && [\text{nonexistence of monopoles}],\\
  & \grad_y\cdot(\eps_r\fieldEL^{\ast}) = \frac{l \rho_0 }{\eps_0E_0} ~\chargeEL^{\ast} && [\text{Gauss's law}],\\
  & \grad_y\times\fieldEL^{\ast} = - \frac{B_0}{E_0} \frac{l}{\tau} ~\dert[s]\fieldMG^{\ast} && [\text{Faraday's law}],\\
  & \grad_y\times\!\brac{\mu_r^{-1}\fieldMG^{\ast}} = \frac{l\mu_0 i_0}{B_0}~ \currentEL^{\ast} + \frac{1}{c_0^2} \frac{E_0}{B_0}\frac{l}{\tau} ~ \dert[s](\eps_r\fieldEL^{\ast}) && [\text{Amp\`{e}re's law}].
\end{align*}
These equations show, that natural choices for $\rho_0$ and $i_0$ are
\begin{align*}
 \rho_0:= \delta_\rho \constEL\frac{E_0}{l} \qquad\text{and}\qquad i_0:=\delta_i \frac{B_0}{l\constMG} ~\text{ for some } \delta_\rho, \delta_i \in\setR_+.
\end{align*}
Furthermore, we note that~~$l/\tau\sim[m/s]$ and $E_0/B_0\sim[m/s]$ are two characteristic velocities of the considered system. More precisely, the characteristic velocity $l/\tau$ is the velocity of the 
considered system, whereas the characteristic velocity~$E_0/B_0$ is the speed of the electromagnetic fields (waves). We suppose, that these characteristic velocities are proportional to the speed of light 
in vacuum~$c_0$. This means, we have
\begin{align*}
  \frac{l}{\tau} = \delta_{\text{V}} c_0 \qquad\text{and}\qquad \frac{E_0}{B_0}  = \delta_{\text{W}} c_0 \qquad \text{ for some } \delta_{\text{V}}\geq0,\delta_{\text{W}} \geq0~.
\end{align*}
Hence, the parameter~$\delta_{\text{V}}$ describes the ratio between the speed of our system and the speed of light, and the parameter $\delta_{\text{W}}$ describes the ratio between the 
electromagnetic fields (waves) and the speed of light. By inserting the preceding relations into the above nondimensional version of Maxwell's equations, we rewrite these equations as
\begin{align*}
  & \grad_y\cdot \fieldMG^{\ast} = 0 && [\text{nonexistence of monopoles}],\\
  & \grad_y\cdot(\eps_r\fieldEL^{\ast}) = ~ \delta_\rho~ \chargeEL^{\ast} && [\text{Gauss's law}],\\
  & \grad_y\times\fieldEL^{\ast} = - \delta_{\text{V}}\delta_{\text{W}}^{-1}~ \dert[s]\fieldMG^{\ast} && [\text{Faraday's law}],\\
  & \grad_y\times\!\brac{\mu_r^{-1}\fieldMG^{\ast}} = \delta_i~ \currentEL^{\ast} + \delta_{\text{V}} \delta_{\text{W}} ~ \dert[s](\eps_r\fieldEL^{\ast}) && [\text{Amp\`{e}re's law}].
\end{align*}
We now pass to the nonrelativistic limit of Maxwell's equations: This means, we confine ourselves to systems, that move magnitudes of orders below the speed of light~$c_0$. Thus, we have a 
ratio~$\delta_{\text{V}}\ll1$ and in the nonrelativistic limit, we let $\delta_{\text{V}}\rightarrow0$. During this limit procedure, we do not touch the speed of the magnetic fields (waves),
which means, that the ratio $\delta_{\text{W}}$ remains constant. Thus, the nonrelativistic limit of Maxwell's equations is given by   
\begin{align*}
  & \grad_y\cdot \fieldMG^{\ast}= 0                                               && [\text{nonexistence of monopoles}],\\
  & \grad_y\cdot(\eps_r\fieldEL^{\ast}) = ~ \delta_\rho~ \chargeEL^{\ast}         && [\text{Gauss's law}],\\
  & \grad_y\times\fieldEL^{\ast} = 0                                              && [\text{Faraday's law}],\\
  & \grad_y\times\!\brac{\mu_r^{-1}\fieldMG^{\ast}} = \delta_i~ \currentEL^{\ast} &&  [\text{Amp\`{e}re's law}].
\end{align*}
From Faraday's law, we conclude that $\fieldEL^{\ast}$ is generated by an electrostatic potential~$\potEL^{\ast}$, i.e., we have $\fieldEL^{\ast}=-\grad_y\potEL^{\ast}$. 
Therefore, we equivalently transform Faraday's law to%
\footnote{Rigorously, we have to guarantee for the following equivalence that Poincar\'{e}'s Lemma holds true, cf. \cite{Rudin-book}} 
\begin{align*}
  & \grad_y\times\fieldEL^{\ast} = 0  \qquad\Equivalent\qquad\fieldEL^{\ast}=-\grad_y\potEL^{\ast} && [\text{Faraday's law}]~.
\end{align*}
This reveals, that Faraday's law and Gauss's law lead to Poisson's equation for the electrostatic potential~$\potEL^{\ast}$. 
In summary, after \enquote{redimensionalization} and inserting the preceding identities, the nonrelativistic limit of Maxwell's equations is given by
\begin{subequations}
\begin{align}
  & \fieldEL = -\grad\potEL                                      && [\text{Faraday's law}],              \label{eq:modelES-maxwell-limitPDE1}\\
  & \grad\cdot(\eps_r\fieldEL) = ~ \frac{1}{\constEL} \chargeEL  && [\text{Gauss's law}],                \label{eq:modelES-maxwell-limitPDE2}\\
  & \grad\cdot\fieldMG = 0                                       && [\text{nonexistence of monopoles}],  \label{eq:modelES-maxwell-limitPDE3} \\
  & \grad\times\!\brac{\mu_r^{-1}\fieldMG} = \constMG~\currentEL && [\text{Amp\`{e}re's law}].           \label{eq:modelES-maxwell-limitPDE4}
\end{align}
\end{subequations}
We are solely interested in electric effects, which are governed by equations \eqref{eq:modelES-maxwell-limitPDE1} \eqref{eq:modelES-maxwell-limitPDE2}.
As these equations are decoupled from the magnetic equations~\eqref{eq:modelES-maxwell-limitPDE3}, \eqref{eq:modelES-maxwell-limitPDE4}, we henceforth omit 
the equations for $\fieldMG$. Thereby, we obtain the electrostatic limit of Maxwell's equations. Here, we refer to the electrostatic limit of Maxwell's equations, 
as the nonrelativistic limit, and additionally neglecting magnetic effects.
\begin{subequations}
\begin{align}
  & \fieldEL = -\grad\potEL                                      && [\text{Faraday's law}], \label{eq:modelES-maxwell-eslimitPDE1}\\
  & \grad\cdot(\eps_r\fieldEL) = ~ \frac{1}{\constEL} \chargeEL  && [\text{Gauss's law}].   \label{eq:modelES-maxwell-eslimitPDE2}
\end{align}
\end{subequations}
Combining these equations leads equivalently to  
\begin{subequations}
\begin{align} \label{eq:modelES-maxwell-eslimitPDE1-mod}
  & -\grad\cdot(\eps_r\grad\potEL) = \frac{1}{\constEL}~ \chargeEL && [\text{Poisson's equation} = \text{Gauss's law}+\text{Faraday's law}]. \tag{\ref{eq:modelES-maxwell-eslimitPDE1}${}^\prime$}
\end{align}
\end{subequations}
This proves, that the macroscopic Maxwell's equations for linear media reduce in the electrostatic limit to Poisson's equation for the electrostatic potential~$\potEL$. 
\medskip
\par
\textit{Henceforth, we assume that the electric phenomena inside the considered electrolyte solutions are sufficiently captured by equation~\eqref{eq:modelES-maxwell-eslimitPDE1-mod}.} 
\medskip
\par
Thus, to account for the electric phenomena inside the considered electrolyte solutions, we solely combine Poisson's equation~\eqref{eq:modelES-maxwell-eslimitPDE1-mod} with the 
remaining conservation laws from \cref{chapter:NonEqThermo}.
\medskip
\par
\begin{remark}[Energy minimization]
 We note, that Poisson's equation~\eqref{eq:modelES-maxwell-eslimitPDE1-mod} is the  Euler-Lagrange equation corresponding to the electrostatic energy functional
 \begin{align*}
  \mathcal{F}\brac{\potEL}:=\Intdx{ \frac{\eps_r}{2}\abs{\grad\potEL}^2 - \frac{1}{\constEL}\chargeEL\potEL ~}~.
 \end{align*}
 Thus, the electrostatic limit is governed by energy minimizing principles, cf. \cite{Dacorogna-book, Guisti-book, LopezDavalos-book}. 
\hfill$\square$
\end{remark}
\begin{remark}[Instantaneous equilibrium assumption]\label{remark:modelES-limitPDEs}
 The nonrelativistic Maxwell's equations~\eqref{eq:modelES-maxwell-limitPDE1}--\eqref{eq:modelES-maxwell-limitPDE4} coincide with the well-known equations of electrostatics and magnetostatics, 
 cf. \cite{EckGarckeKnabner-book, LopezDavalos-book, LifshitzLandau-book2}. However, electrostatics and magnetostatics investigate equilibrium states. Thus, the electrostatic fields and the magnetostatic fields 
 are temporal constant. In contrast to this, the preceding nonrelativistic Maxwell's equations~\eqref{eq:modelES-maxwell-limitPDE1}--\eqref{eq:modelES-maxwell-limitPDE4} 
 are formulated for temporal variable fields. Hence, there are ongoing temporal dynamics. Nevertheless, the structure of the nonrelativistic Maxwell's equations coincides with the equations of 
 electrostatics and magnetostatics. Thus, the temporal dynamics are assumed to take place such that the electromagnetic fields instantaneously switch from one equilibrium state to another one. 
 This assumption is appropriate for time scales, which are orders of magnitudes above the relaxation time for electromagnetic phenomena. In terms of the above parameters $\delta_{\text{V}}$ 
 and $\delta_{\text{W}}$, this applies for $\delta_{\text{V}}\ll\delta_{\text{W}}$. This means, that, e.g., the function $t\mapsto\potEL(t,\cdot)$ is a one-parameter collection of equilibrium 
 potentials~$\potEL(t,\cdot)$. Hence, the dynamics in $t$ do not resolve temporal dynamics in between two equilibrium states. In this connection, we recall that we assumed~$\dert\potEL=0$ in 
 equation~\eqref{eq:NonEqThermo-energy-electricEnergyPDE}, whereas now, we generally have $\dert\potEL\neq0$. However, equation~\eqref{eq:NonEqThermo-energy-electricEnergyPDE} exactly resolves 
 nonequilibrium dynamics in between two equilibrium states. As $\dert\potEL$ does not resolve these dynamics, we continue to neglect $\dert\potEL$ in \eqref{eq:NonEqThermo-energy-electricEnergyPDE}.   
\hfill$\square$
\end{remark}
\begin{remark}[Nonrelativistic limit equations]\label{remark:modelES-maxwellScaling}
 The nonrelativistic Maxwell's equations~\eqref{eq:modelES-maxwell-limitPDE1}--\eqref{eq:modelES-maxwell-limitPDE4} may change, if we use a different scaling for $\delta_{W}$. 
 More precisely, choosing the ansatz~$\delta_{W}:=\delta_{V}^\alpha$, we previously set~$\alpha=0$, and even for $\alpha\in(0,1)$, we come to the same limit equations. However, for $\alpha=1$, 
 we obtain a different limit of Faraday's law, which reads as $\grad\times\fieldEL = - \dert\fieldMG$. Hence, in this case, the magnetic effects do not decouple from the electric effects in the 
 nonrelativistic limit of Maxwell's equations.
\hfill$\square$
\end{remark}
%
\begin{remark}[Electromagnetic potentials]
  Combining Helmholtz's decomposition, cf. \cite{Monk-book}, and the nonexistence of monopoles shows that 
  \begin{align*}
    \grad\cdot \fieldMG= 0 \qquad\Equivalent\qquad \fieldMG=\grad\times\vecA && [\text{nonexistence of monopoles}].
  \end{align*}
  Hence, we can express the magnetic field in terms of a vector potential~$\vecA$. Commonly, $(\vecA,\potEL)$ are known as the electromagnetic potentials, cf. \cite{LopezDavalos-book}. 
  Furthermore, we can transform Maxwell's equations such that the resulting \enquote{potential equations} are solely solved by $(\vecA,\potEL)$.
  To this end, we combine equations~\eqref{eq:modelES-maxwell-limitPDE1} and \eqref{eq:modelES-maxwell-limitPDE2} to obtain Poisson's equation for $\potEL$. To compute~$\vecA$, it suffices to solve 
  equation~\eqref{eq:modelES-maxwell-limitPDE4}. Thus, the nonrelativistic limit of Maxwell's equations transform to the \enquote{potential equations} 
  \begin{align*}
    & -\grad\cdot(\eps_r\grad\potEL) = \frac{1}{\eps_0}~ \chargeEL && [\text{Poisson's equation}],\\
    & -\Delta\vecA = \mu_r\mu_0~ \currentEL && [\text{Amp\`{e}re's law}]. 
  \end{align*}
  Here, we assumed $\mu_r$ to constant and we involved the identity~$\grad\times\grad\times\vecA=\grad(\grad\cdot\vecA)-\Delta\vecA$ together with Coulomb's gauge~$\grad\cdot\vecA=0$, 
  cf. \cite{LopezDavalos-book, LifshitzLandau-book2}. Thus, for constant $\mu_r$, Maxwell's equations reduce in the nonrelativistic limit to two decoupled elliptic equations for the electromagnetic potentials. 
  Whereas, in the relativistic case, Maxwell's equations transform with the Lorentz gauge, cf. \cite{LopezDavalos-book, LifshitzLandau-book2}, to two coupled hyperbolic wave equations for the 
  electromagnetic potentials~$(\vecA,\potEL)$. This reveals, that in the nonrelativistic limit Maxwell's equations switch from hyperbolic to elliptic. 
\hfill$\square$
\end{remark}
\fancyhead[L]{Governing equations}
\section{The Governing Equations}\label{sec:modelES-pdes}
First of all, we note that subsequently assumptions~\ref{Assump:Domain}--\ref{Assump:entropy1a} from \cref{chapter:NonEqThermo} continue to hold true. 
Thus, in particular, we henceforth suppose that the considered electrolyte solutions are multicomponent mixtures of $L$ different charged constituents, 
which are indexed such that the $L$th chemical species is the solvent. 
\medskip
\par
For the sake of completeness and to henceforth avoid permanent cross-referencing to \cref{chapter:NonEqThermo}, we now briefly list the general equations 
from \cref{chapter:NonEqThermo} and \cref{sec:modelES-maxwell}, which govern the dynamics of electrolyte solutions.
\\[3.0mm]
\begin{subequations}
\textbf{1. Electric potential equation: } According to \ref{Assump:electrostatics}, we have $\fieldEL=-\grad\potEL$ for the electric field~$\fieldEL$, and due to \eqref{eq:NonEqThermo-charge-DefFreeCharge} 
and \eqref{eq:modelES-maxwell-eslimitPDE1-mod}, the electric potential~$\potEL$ solves 
\begin{align}\label{eq:modelES-pdes-poisson}
  -\grad\cdot(\eps_r\grad\potEL) = \frac{1}{\constEL}~ \chargeEL \qquad\text{with }\quad \chargeEL=\sum_l\frac{\constCharge\zl}{\ml}\rol. 
\end{align}
\textbf{2. Mass conservation equations: } For $l\in\cbrac{1,\ldots,L-1}$, we have
\begin{align}
  & \dert\rol + \grad\cdot\brac{ \rol\fieldF + \rolfluxrel} ~=~\rl~, \label{eq:modelES-pdes-massBalance}    \\[2.0mm]
  & \dert\rol[] + \grad\cdot\brac{ \rol[] \fieldF } = 0~.            \label{eq:modelES-pdes-massBalanceTot}
\end{align}
Furthermore, the \concRol~$\rol[L]$ of the solvent and the \massfluxRol~$\rolfluxrel[L]$ of the solvent are given according to \cref{remark:NonEqThermo-independentDriftFluxes} 
and \cref{remark:NonEqThermo-independentVariables} by 
\begin{align}\label{eq:modelES-pdes-sumCondition}
 \rol[L] = \rol[] -\sum_{l=1}^{L-1} \rol \qquad\text{ and }\qquad \rolfluxrel[L] = - \sum_{l=1}^{L-1} \rolfluxrel~. 
\end{align}
\textbf{3. Momentum conservation equations: } For the barycentric momentum density holds according to \eqref{eq:NonEqThermo-mom-momBalanceTot} 
\begin{align}\label{eq:modelES-pdes-momBalanceTot}
  \dert(\rol[]\fieldF) + \grad\cdot\brac{\rol[]\fieldF\otimes\fieldF} = \grad\cdot\stressTotl[] + \chargeEL\fieldEL ~. 
\end{align}
Moreover, according to~\ref{Assump:mixtureStress2}, the mixture stress tensor~$\stressTotl[]$ is given by
\begin{align}\label{eq:modelES-pdes-mixtureStress}
 \stressTotl[] = -\pressHydrl[]\mathds{1} + \stressViscl[]~. 
\end{align}
This equation defines the mixture pressure~$\pressHydrl[]$ and the viscous mixture stress tensor~$\stressViscl[]$. Furthermore, in \eqref{eq:NonEqThermo-mom-extendedDaltonsLaw} and 
\eqref{eq:NonEqThermo-mom-viscStressMix} we obtained more detailed expressions for these quantities. However, in \cref{chapter:NonEqThermo}, we distinguished between the mixture pressure~$\pressHydrl[]$
and the total mixture pressure~$\pressTotl[]$, which was defined in \ref{Assump:thermoPressure} by $\pressTotl[]:=-\frac{1}{n}\trace{\stressTotl[]}$. 
According to \eqref{eq:NonEqThermo-mom-pressTotpressHydr}, these pressures are related by 
\begin{align}\label{eq:modelES-pdes-pressures}
 \pressHydrl[] = \pressTotl[] + \frac{1}{n}\trace{\stressViscl[]}~.  
\end{align}
\textbf{4. Energy conservation equation: } In \cref{sec:NonEqThermo-energy}, we proposed the following ansatz for the total energy density~$\rol[]\energyTotl[]$ in \eqref{eq:NonEqThermo-energy-ansatzTot-a}:  
\begin{align*}
  \rol[]\energyTotl[] = \rol[]\energyIntPurel[] + \rol[]\energyIntMixl[] + \chargeEL\potEL + \frac{1}{2}\rol[]\abs{\fieldF}^2\,. 
\end{align*}
For the total energy density, we formulated the \textit{first law of thermodynamics} in \eqref{eq:NonEqThermo-energy-firstLaw} as
\begin{align*}
  \dert\brac{\rol[]\energyTotl[]} + \grad\cdot\!\brac{ \rol[]\energyTotl[]\fieldF + \vecq + \potEL\currentEL - \stressTotl[]\fieldF } = 0 ~.
\end{align*}
In particular, the decomposition $\energyIntl[]=\energyIntMixl[]+\energyIntPurel[]$ of the internal energy~$\energyIntl[]$ into a pure substance part~$\energyIntPurel[]$ and a 
part due to mixing~$\energyIntMixl[]$ implied in \cref{remark:NonEqThermo-chemicalPotential} resp. \cref{remark:NonEqThermo-electrochemicalPotential} 
the splittings $\chempotl = \chempotMixl + \chempotPurel$ resp.~$\elchempotl = \elchempotMixl + \chempotPurel$ of the chemical potentials~$\chempotl$ resp.~the electrochemical potentials~$\elchempotl$ 
into their respective parts due mixing~$\chempotMixl, \elchempotMixl$ and their pure substance parts~$\chempotPurel$.
\par
Furthermore, each part of the total energy density~$\rol[]\energyTotl[]$ is subject to an evolution equation. More precisely, $\rol[]\energyIntMixl[]$ solves \eqref{eq:NonEqThermo-energy-energyIntMixPDE},  
$\chargeEL\potEL$ solves \eqref{eq:NonEqThermo-energy-electricEnergyPDE}, and $\frac{1}{2}\rol[]\abs{\fieldF}^2$ solves \eqref{eq:NonEqThermo-energy-kineticEnergyPDE}. 
However, we derived these equations by suitable manipulations of \eqref{eq:modelES-pdes-massBalance} and~\eqref{eq:modelES-pdes-momBalanceTot}. Thus, their information content 
is essentially contained in \eqref{eq:modelES-pdes-massBalance} and~\eqref{eq:modelES-pdes-momBalanceTot}. On the other hand, the evolution equation~\eqref{eq:NonEqThermo-energy-energyIntPurePDE} 
for $\rol[]\energyIntPurel[]$ is independent of the other governing equations. For this reason, we add   
\begin{align}\label{eq:modelES-pdes-energyIntPurePDE}
  &~~\dert\brac{\rol[]\energyIntPurel[]} + \grad\cdot\brac{ \rol[]\energyIntPurel[]\fieldF + \vecq - \sum_l\elchempotMixl\rolfluxrel + \potEL\currentEL } \nonumber\\
  &= \stressTotl[]:\grad\fieldF -\sum_l\grad\elchempotMixl\cdot\rolfluxrel - \sum_l\elchempotMixl\rl ~ 
\end{align}
to the set of governing equations.\\[3.0mm]
\textbf{5. Entropy evolution equation: } For the specific entropy~$\entropy$, we have $\entropy = \entropyPure + \entropyMixl[]$ due to~\ref{Assump:entropy1}. 
Here, $\entropyPure$ is the specific pure substance part and $\entropyMixl[]$ the specific entropy of mixing. In \cref{sec:NonEqThermo-entropy}, we proved that the 
evolution equation~\eqref{eq:NonEqThermo-entropy-entropyPDE1} for the entropy density~$\rol[]\entropy$ decomposes into the evolution equation~\eqref{eq:NonEqThermo-entropy-entropyMixPDE} 
for $\rol[]\entropyMixl[]$ and the evolution equation~\eqref{eq:NonEqThermo-entropy-entropyPurePDE} for $\rol[]\entropyPure$. However, we derived 
these equations based on Gibbs relation~\eqref{eq:NonEqThermo-entropy-gibbsNonEq} and \eqref{eq:modelES-pdes-massBalance}. Hence, the information content of these equations is essentially contained 
in \eqref{eq:modelES-pdes-massBalance}, \eqref{eq:modelES-pdes-massBalanceTot}, and \eqref{eq:modelES-pdes-energyIntPurePDE}. For this reason, we exclude the equations for $\rol[]\entropy$, 
$\rol[]\entropyMixl[]$, and $\rol[]\entropyPure$ from the set of governing equations. 
\par
Nevertheless, concerning the \textit{second law of thermodynamics}~\eqref{eq:NonEqThermo-entropy-secondLaw}, we established due to \eqref{eq:NonEqThermo-entropy-reactionCriterion}, \eqref{eq:NonEqThermo-entropy-entropyFlux},
and \eqref{eq:NonEqThermo-entropy-diss2} the formulation  
\begin{flalign}\label{eq:modelES-pdes-diss}
  ~~~0 \leq& -\frac{1}{\temp^2}\grad\,\temp\cdot\brac{\vecq+\potEL\currentEL-\sum_l\elchempotl\rolfluxrel} + \frac{1}{\temp}\stressViscl[]:\grad\fieldF 
             -\frac{1}{\temp}\sum_l\grad\elchempotl\cdot\rolfluxrel  - \sum_l \frac{\chempotl}{\temp}\rl ~.& 
\end{flalign}
Subsequently, we use exactly this inequality to validate the constitutive ansatzes for the drift \massfluxRol{es}~$\rolfluxrel$, the mass production rates~$\rl$, 
the viscous stress tensor~$\stressViscl[]$, the heat flux~$\vecq$, and for the internal energies~$\energyIntl$, which determine the electrochemical potentials~$\elchempotl$.  
\end{subequations}
\fancyhead[L]{Internal energy}
\section{Constitutive Ansatz for the Internal Energy}\label{sec:modelES-energy}
In \ref{Assump:energy1a} we assumed $\rol\energyIntl = \rol\energyIntPurel + \rol\energyIntMixl$ for the total internal energy densities.
Thereby, we arrived in \eqref{eq:NonEqThermo-energy-energyIntMix} and \eqref{eq:NonEqThermo-energy-energyIntPure} for the pure substance internal energy density~$\rol[]\energyIntPurel[]$ 
and the internal energy of mixing density~$\rol\energyIntMixl[]$ at
\begin{align}\label{eq:modelES-energy-repeatDefIntEnergies}
  \rol[]\energyIntMixl[] = \rol[]\sum_l \yl \energyIntMixl
  \qquad\text{resp.}\qquad 
  \rol[]\energyIntPurel[] = \rol[]\sum_l \yl \energyIntPurel + \frac{\rol[]}{2}\sum_l\yl\abs{\fieldF_l-\fieldF}^2~.   
\end{align}
In continuation of the previous assumptions, we now introduce the crucial ansatzes for the specific internal energies~$\energyIntMixl$, $\energyIntPurel$, and $\energyIntPurel[]$.
\begin{enumerate}[align=left, leftmargin=*, topsep=1.0mm, itemsep=-2.0mm, label={(\labelA\arabic*)}, start=\value{countModAssump}+1]
  \item \label{Assump:ansatzEnergyMix}%
        \textbf{Internal energy of mixing: } For the specific internal energies of mixing~$\energyIntMixl\sim[J/kg]$ from \ref{Assump:energy1a}, we supposed $\energyIntMixl=\energyIntMixl(\yl)$ 
        in \ref{Assump:energy2a}. In accordance with this functional dependency, we now assume similar to \cite{KrautleExPaper, Krautle-habil}, with a given real number~$\beta_l$, the ansatz
	 \begin{align*}
	   \yl\energyIntMixl(\yl) := \frac{\constBoltz\temp}{\ml}\, \yl\brac{ \beta_l -1+\ln(\yl)} + \frac{\constBoltz\temp}{\ml} \exp(-\beta_l) \qquad\sim\sqbrac{J \, kg^{-1}}.
	 \end{align*}
  \item \label{Assump:ansatzEnergyPure}%
        \textbf{Pure substance internal energy: } For the specific pure substance internal energies~$\energyIntPurel$ from \ref{Assump:energy1a}, 
        we supposed $\energyIntMixl=\energyIntMixl(\entropy,\volumeSpecific)$ in \ref{Assump:energy2a}. In accordance with this functional dependency, we now assume, with 
        a general specific energy function $\energyIntPureAnsatz(\entropyPure,\volumeSpecific)\sim[J/kg]$, the ansatz   
	 \begin{align*}
	   \yl\energyIntPurel(\entropy,\volumeSpecific) := \yl\energyIntPureAnsatz(\entropyPure,\volumeSpecific) \qquad\sim\sqbrac{J \, kg^{-1}}.
	 \end{align*}
  \item \label{Assump:ansatzEnergyPureTot}%
        \textbf{Total pure substance internal energy: } The velocities~$\fieldF_1,\ldots,\fieldF_L,\fieldF$ we treated according to \eqref{eq:NonEqThermo-entropy-energyIntPure} as parameters 
        for the total specific pure substance internal energy~$\energyIntPurel[]\sim[J/kg]$ from \eqref{eq:modelES-energy-repeatDefIntEnergies}. Henceforth, we assume, that we can neglect 
        this parameter dependency, i.e., instead of \eqref{eq:modelES-energy-repeatDefIntEnergies}, we suppose together with \ref{Assump:ansatzEnergyPure} the ansatz 
        \begin{align*}
         \rol[]\energyIntPurel[] = \rol[]\energyIntPurel[]\brac{\entropy, \volumeSpecific, \yl[1],\ldots, \yl[L]} 
                           = \rol[]\sum_l\yl\energyIntPurel\brac{\entropy, \volumeSpecific} 
                           = \rol[]\sum_l\yl\energyIntPureAnsatz\brac{\entropyPure, \volumeSpecific}
                           = \rol[]\energyIntPureAnsatz\brac{\entropy, \volumeSpecific}~.
        \end{align*} 
\setcounter{countModAssump}{\value{enumi}}%
\end{enumerate}
In \cref{remark:modelES-modelComputation}, we present a possible choice of~$\energyIntPureAnsatz$. Moreover, it is important that analogously to \cref{remark:NonEqThermo-chemicalPotential}, 
we obtain for the chemical potentials from~\ref{Assump:ansatzEnergyMix}--\ref{Assump:ansatzEnergyPureTot}  the crucial ansatzes
\begin{align}\label{eq:modelES-energy-chemPot}
   \chempotl = \chempotPurel + \chempotMixl 
   \quad\text{with}\quad \chempotPurel = \energyIntPureAnsatz\brac{\entropyPure, \volumeSpecific}  
   \quad\text{and}\quad     \chempotMixl = \frac{\constBoltz\temp}{\ml} \, \brac{\beta_l + \ln(\yl)} ~. 
\end{align}
\begin{remark}[Limitation of the ansatzes]
 Assumption~\ref{Assump:ansatzEnergyPure} for $\energyIntPurel$ is motivated by the fact, that we have a uniquely defined specific entropy~$\entropy$ and a uniquely defined 
 specific volume~$\volumeSpecific$ inside the mixture. Thus, the energetic contribution caused by these variables should be the same for all constituents.   
 Furthermore, replacing \eqref{eq:modelES-energy-repeatDefIntEnergies} by \ref{Assump:ansatzEnergyPureTot} is admissible as long as the kinetic contributions 
 due to the drift velocities~$\fieldF_l-\fieldF$ are small compared to the entropic contribution plus the volumetric contribution. 
 \hfill$\square$ 
\end{remark}
\fancyhead[L]{Reaction rates}
\section{Constitutive Ansatz for the Reaction Rates}\label{sec:modelES-reactions}
In this section, we briefly repeat the basic chemical definitions. For a detailed introduction to chemical reactions, 
we refer to \cite{FeinbergNetwork1, FeinbergNetwork2, RubinReactions, SmithMissen-book, PrigogineKondepudi-book, Upadhay-book}.
Henceforth, we consider general chemical reactions, which transform some constituents of the mixture into other ones. 
These chemical reactions can be described by stoichiometric equations. Provided we denote the involved constituents of the mixture by $C_l$, e.g., 
the stoichiometric equation for the $j$th chemical reaction may look like
\begin{align}\label{eq:modelES-reactions-reacEQ}
 \tilde{s}_{1j}\,C_1 + \tilde{s}_{2j}\,C_2 + \tilde{s}_{4j}\,C_4 ~\rightleftharpoons~ \tilde{s}_{6j}\,C_6~, \qquad\text{with}~~\slj\in\setN. 
\end{align}
We formally rearrange this stoichiometric equation to 
\begin{align*}
 0 ~\rightleftharpoons~ -\tilde{s}_{1j}\,C_1 - \tilde{s}_{2j}\,C_2 - \tilde{s}_{4j}\,C_4  + \tilde{s}_{6j}\,C_6 ~, \qquad\text{with}~~\slj\in\setN.
\end{align*}
This equation shows, that the constituents $C_1$, $C_2$, $C_4$, $C_6$ participate in the $j$th chemical reaction.
More precisely, $C_1$, $C_2$, $C_4$ are the so-called reactants and $C_6$ is the so-called product of the $j$th chemical reaction. Furthermore,
\begin{align*}
 \slj[1]:=-\tilde{s}_{1j},~~\slj[2]:=-\tilde{s}_{2j},~~\slj[4]:=-\tilde{s}_{4j},~~\text{and}~~\slj[6]:=\tilde{s}_{6j}
\end{align*}
are the dimensionless stoichiometric coefficients of the involved constituents. Additionally, we define for the remaining constituents that are not affected by the $j$th chemical reaction, 
the stoichiometric coefficients by $\slj=0$. Thus, the $j$th chemical reaction is described by the reaction vector
\begin{align*}
\vecs_j:= \brac{\slj[1], \slj[2], 0, \slj[4], 0, \slj[6], 0, \ldots, 0}^\top  \in\setZ^L~.
\end{align*}
In case of $J\in\setN$ reactions, we define the so-called stoichiometric matrix~$S\in\setZ^{L\times J}$ by
\begin{align*}
  S = \brac{\vecs_1,\ldots,\vecs_J} 
    =\left(\begin{array}{c c c}
         s_{11} & \ldots & s_{1J}  \\
         \vdots & \ddots & \vdots  \\
         \vdots & \ddots & \vdots  \\
         s_{L1} & \ldots & s_{LJ} 
      \end{array}\right)
 \in \setZ^{L\times J}~.
\end{align*}
Here, e.g., the $j$th column is given by the above reaction vector~$\vecs_j$. Thus, each column of~$S$ describes a chemical reaction. Henceforth, we assume that we have for the stoichiometric matrix~$S$
\begin{align}
  \text{rank}\brac{S} = J < L~.
\end{align}
This assumption implies linear independency of the chemical reactions, i.e., none of the chemical reactions can be reproduced by arbitrary combinations of the remaining ones. 
Since the maximal number of linear independent chemical reactions is bounded by the number of available constituents, we furthermore restrict~$J<L$.
\medskip
\par
Following \cite{Atkins-book, EckGarckeKnabner-book, SmithMissen-book, PrigogineKondepudi-book, Upadhay-book}, we define for the exemplary chemical reaction~\eqref{eq:modelES-reactions-reacEQ}, 
with the mass fractions~$\yl$, the corresponding mathematical reaction rate~$\Rj$ by
\begin{align*}
 \Rj = \Rj^f - \Rj^b := \kfj{\yl[1]}^{-\slj[1]}{\yl[2]}^{-\slj[1]}{\yl[4]}^{-\slj[4]} - \kbj{\yl[6]}^{\slj[6]} \qquad\sim\sqbrac{m^{-3}s^{-1}}.
\end{align*}
Here, $\Rj^f\sim[1/(m^3s)]$ is the so-called forward reaction rate, which models the \enquote{$\rightharpoonup$}-reaction in ~\eqref{eq:modelES-reactions-reacEQ}, and $\Rj^b\sim[1/(m^3s)]$ 
is the so-called backward reaction rate, which describes in~\eqref{eq:modelES-reactions-reacEQ} the \enquote{$\leftharpoondown$}-reaction. Furthermore, $\kfj\sim[1/(m^3s)]$ is the so-called 
forward rate constant and $\kbj\sim[1/(m^3s)]$ the so-called backward rate constant of the $j$th chemical reaction. 
Generally, we suppose that for each of the $J$ chemical reactions, the corresponding mathematical reaction rate is given by
\begin{align}
 \Rj = \Rj^f - \Rj^b := \Rfjmal - \Rbjmal \qquad\sim\sqbrac{m^{-3}s^{-1}},
\end{align}
where again $\kfj\sim[1/(m^3s)]$ denotes the forward rate constant and $\kbj\sim[1/(m^3s)]$ the backward rate constant. 
Next, we define for the $j$th chemical reaction rate~$\Rj$ the so-called equilibrium constant $\Kj$ by
\begin{align*}
  \Kj= \frac{\kfj}{\kbj} \Hence \Kj = \prod_{\slj\neq 0} \brac{\yl[i]}^{\slj} \text{in case of }~ \Rj=0~.
\end{align*} 
This shows that in chemical equilibrium, i.e., $\Rj=0$, the product of the right-hand side is constant with constant value~$\Kj$. 
Furthermore, the equation $\Kj = \prod_{\slj\neq 0} \brac{\yl}^{\slj}$ is exactly the equilibrium mass action law, cf. \cite{EckGarckeKnabner-book, PrigogineKondepudi-book, Upadhay-book}. 
This is the reason, why we refer to the reaction rates~$\Rj$ as reaction rates according to mass action law. Next, we obtain the total reaction rate~$\Rl^{tot}$ for the $l$th constituent by multiplying 
the elementary reaction rates~$\Rj[1],\ldots,\Rj[L]$ by the stoichiometric coefficient~$\slj$ summing over $j$. Here, $\slj$ is the stoichiometric coefficient of the $l$th constituent in $j$th reaction. 
Thus, the total reaction rate~$\Rl^{tot}$ of the $l$th constituent is given by 
\begin{align}\label{eq:modelES-reactions-ansatzReactionTot}
 \Rl^{tot}:=\sum_j\slj\Rj = \sum_j\slj\sqbrac{\Rfjmal - \Rbjmal} \qquad\sim\sqbrac{m^{-3}s^{-1}}.
\end{align}
We now state the fundamental relation between the reaction rates~$\Rj\sim[1/(m^3s)]$ and the mass production rates~$\rl\sim[kg/(m^3s)]$. 
More precisely, we suppose for the mass production rates~$\rl$ the \textit{\textbf{constitutive ansatz}} 
\begin{align}\label{eq:modelES-reactions-ansatz}
 \rl:=\ml\Rl^{tot} = \ml\sum_j\slj\sqbrac{\Rfjmal - \Rbjmal} \quad\sim\sqbrac{kg\,m^{-3}s^{-1}}.
\end{align}
\medskip
\par
Next, we demonstrate that the mass production rates~$\rl$ from \eqref{eq:modelES-reactions-ansatz} are subject to the mass conservation property~\ref{Assump:reactions}. 
More precisely, the mass conservation property~\ref{Assump:reactions} applies due to
\begin{align}\label{eq:modelES-reactions-massCrit}
  \sum_l \ml\slj =0~ \qquad \text{for}\quad j\in\cbrac{1,\ldots,L}.
\end{align}
As the general structure of mass production rates~$\rl$ from \eqref{eq:modelES-reactions-ansatz} is contained in the exemplary mass production rate corresponding to~\eqref{eq:modelES-reactions-reacEQ}, 
it suffices to concentrate on this example. To this end, we multiply the components of the reaction vector~$\vecs_j$ by the respective 
molecular masses~$\ml[1]$, $\ml[2]$, $\ml[4]$, and $\ml[6]$. Thereby, we obtain the mass transfer vector
\begin{align*}
 \brac{\ml[1]\slj[1], \ml[2]\slj[2], 0, 0, \ml[4]\slj[4], 0, \ml[6]\slj[6], 0, \ldots, 0}^\top~\in\setZ^L~.
\end{align*}
Note that due to the stoichiometry~\eqref{eq:modelES-reactions-reacEQ}, \eqref{eq:modelES-reactions-ansatz}, $\slj[6]$ molecules of the product~$C_6$ possess the molecular weight   
\begin{align*}
 \slj[6]\ml[6]= \abs{\slj[1]}\ml[1]+\abs{\slj[2]}\ml[2]+\abs{\slj[4]}\ml[4]~.
\end{align*}
Thus, summing over the components of the mass transfer vector, leads together with the definition of the stoichiometric coefficients~$\slj$ and $\ml[6]$ to
\begin{align*}
 \sum_l \ml\slj  = -\ml[1]\abs{\slj[1]} - \ml[2]\abs{\slj[2]} -\ml[4]\abs{\slj[4]} + \ml[6]\slj[6] = 0~.
\end{align*}
\medskip
\par
Additionally, the mass production rates~$\rl$ from \eqref{eq:modelES-reactions-ansatz} are subject to the charge conservation property~\ref{Assump:reactionsAndCharges} due to  
\begin{align}\label{eq:modelES-reactions-chargeCrit}
  \sum_l \zl\slj =0~.
\end{align}
Again, it suffices to verify this criterion for the exemplary mass production rate corresponding to~\eqref{eq:modelES-reactions-reacEQ}. For that purpose, we assume for a moment, 
that in the exemplary chemical reaction~\eqref{eq:modelES-reactions-reacEQ}, the constituents~$C_1$ and $C_2$ are electrically charged chemical species, whereas $C_4$ is a electrically neutral. 
Thus, we have the valency $\zl[4]=0$. Multiplying the reaction vector~$\vecs_j$ by the respective valencies~$\zl[1]$, $\zl[2]$, $\zl[4]$, and $\zl[6]$, we obtain the charge transfer vector
\begin{align*}
 \brac{\zl[1]\slj[1], \zl[2]\slj[2], 0, 0, 0, 0, \zl[6]\slj[6], 0, \ldots, 0}^\top~\in\setZ^L~.
\end{align*}
According to the stoichiometry~\eqref{eq:modelES-reactions-ansatz}, the valency of $C_6$ is given by 
\begin{align*}
 \zl[6]\slj[6] := \zl[1]\slj[1]+ \zl[2]\slj[2]~.
\end{align*}
This is owing to the fact that chemical reactions solely transfer electric charges, and not create charges. Hence, summing over the components of the charge transfer vector results in
\begin{align*}
 \sum_l \zl\slj  = -\zl[1]\abs{\slj[1]} - \zl[2]\abs{\slj[2]} + \zl[6]\slj[6] = 0~.
\end{align*}
\medskip
\par
It now remains to show that the constitutive ansatz~\eqref{eq:modelES-reactions-ansatz} is in accordance with the second law of thermodynamics. From \eqref{eq:modelES-pdes-diss}, 
we know that this is the case, if the sufficient condition
\begin{align*}
 - \sum_l \frac{\chempotl}{\temp}\rl  \geq 0 
\end{align*}
holds true. Together with $\chempotl=\chempotMixl +\chempotPurel$ from \cref{remark:NonEqThermo-electrochemicalPotential}, we furthermore strengthen this criterion to  
\begin{align}\label{eq:modelES-reactions-thermoCrit}
  - \sum_l \frac{\chempotPurel}{\temp}\rl  \geq 0~ \qquad\text{ and }\qquad  -\sum_l \frac{\chempotMixl}{\temp}\rl    \geq 0~.
\end{align}
Note, that the chemical potentials are given in \eqref{eq:modelES-energy-chemPot} by  
\begin{align*}
 \chempotPurel = \energyIntPureAnsatz 
 \qquad\text{~~and~~}\qquad 
 \chempotMixl = \frac{\constBoltz\temp}{\ml} \, \brac{\beta_l + \ln(\yl)}~. 
\end{align*}
In particular, the definition of the pure substance chemical potentials~$\chempotPurel$ immediately results with the mass conservation property~\ref{Assump:reactions} in
\begin{align*}
 - \sum_l \frac{\chempotPurel}{\temp}\rl   =  - \frac{\energyIntPureAnsatz}{\temp}\sum_l \rl = 0~.
\end{align*}
This proves the first inequality in \eqref{eq:modelES-reactions-thermoCrit}. Furthermore, the preceding equation reveals that these ansatzes for the pure substance chemical potentials~$\chempotPurel$
never lead to production of specific pure substance entropy~$\entropyPure$. 
\par
As to the second inequality in \eqref{eq:modelES-reactions-thermoCrit}, we follow the ideas of \cite{KrautleExPaper}. For that purpose, we firstly define with the equilibrium constants~$\Kj$ 
the vector~$\vecK\in\setR^J$ by
\begin{align*}
 \vecK:= \brac{-\ln\Kj[1],\ldots,-\ln\Kj[J]} \in\setR^J~,
\end{align*}
and we collect the constants~$\beta_l$ from \ref{Assump:ansatzEnergyMix} in a vector $\vecbeta\in\setR^L$. We fix these constants by choosing them such that~$\vecbeta$ solves the linear equation system 
\begin{align*}
    S^\top \vecbeta = \vecK~.
\end{align*}
Due to $\text{rank}\brac{S} = J < L$, this linear equation system has at least one solution $\vecbeta$, which can be chosen, e.g., as 
$\vecbeta := \min\cbrac{~|\vecv|:~~S^\top \vecv = \vecK }$. Equipped with these definitions, we deduce the fundamental equivalence:
\begin{align*}
\Rj \gtreqqless 0   &  \Equivalent  \Rj^f \gtreqqless \Rj^b  
                      ~\Equivalent~ \ln\Rj^f \gtreqqless \ln\Rj^b
		      ~\Equivalent~  0 \gtreqqless -\ln\Kj + \sum_{l=1}^L \slj\ln\yl                     \nonumber \\
		    &  \Equivalent   0 \gtreqqless \sum_{l=1}^L \slj\brac{ \beta_l + \ln\yl} \qquad\text{for } j \in\cbrac{1,\ldots,J}.
\end{align*}
Hence, we obtain the estimates
\begin{align*}
\Rj \sum_{l=1}^L \slj\brac{ \beta_l + \ln\yl} ~\leq~0 \qquad \text{for } j \in\cbrac{1,\ldots,J}~,
\end{align*}
which finally lead us with $\chempotMixl$ from \eqref{eq:modelES-energy-chemPot} and $\rl$ from \eqref{eq:modelES-reactions-ansatz} to
\begin{align*}
   -\sum_l \frac{\chempotMixl}{\temp} ~\rl 
  = -\sum_l \frac{\constBoltz\temp}{\ml\temp} \brac{\beta_l + \ln(\yl)} ~\ml\sum_l\slj\Rj 
  = -\constBoltz \sum_j \Rj \sum_l \slj \brac{\beta_l + \ln(\yl)} \geq 0~.
\end{align*} 
\fancyhead[L]{Diffusion fluxes}
\section{Constitutive Ansatzes for the Diffusion Fluxes}\label{sec:modelES-diffusion}
As to the drift \massfluxRol{es}~$\rolfluxrel$, we firstly recall the sum condition~\eqref{eq:modelES-pdes-sumCondition}
\begin{align*}
 \rolfluxrel[L] = - \sum_{l=1}^{L-1} \rolfluxrel 
 \qquad\Equivalent\qquad
 \sum_l \rolfluxrel = 0~.
\end{align*}
Moreover, the constitutive ansatzes for~$\rolfluxrel$ are in accordance with the second law of thermodynamics~\eqref{eq:modelES-pdes-diss}, 
if the sufficient condition 
\begin{align*}
  - \frac{1}{\temp}\sum_l \grad\elchempotl\cdot\rolfluxrel ~\geq~ 0
  \qquad\Equivalent\qquad
 -\sum_l \grad\elchempotl\cdot\rolfluxrel ~\geq~ 0
\end{align*}
holds true. Recalling the splitting $\elchempotl = \elchempotMixl + \chempotPurel$ from \cref{remark:NonEqThermo-electrochemicalPotential}, 
we strengthen this condition to
\begin{align}\label{eq:modelES-diffusion-criterion}
  -\sum_l \grad\chempotPurel\cdot\rolfluxrel ~\geq~ 0~
  \qquad\text{and}\qquad
  -\sum_l \grad\elchempotMixl\cdot\rolfluxrel ~\geq~ 0~.
\end{align}
\par
The first inequality of \eqref{eq:modelES-diffusion-criterion} follows immediately by inserting the ansatz~\eqref{eq:modelES-energy-chemPot} 
for the pure substance chemical potentials~$\chempotPurel$. More precisely, we obtain for the drift \massfluxRol{es} with the above sum condition 
\begin{align*}
   -\sum_l \grad\chempotPurel\cdot\rolfluxrel = -\grad\energyIntPureAnsatz\cdot\brac{\sum_l \rolfluxrel} = 0~. 
\end{align*}
Hence, the chosen ansatzes~\eqref{eq:modelES-energy-chemPot} for the pure substance chemical potentials~$\chempotPurel$ never cause production 
of specific pure substance entropy~$\entropyPure$. 
Concerning the second inequality of \eqref{eq:modelES-diffusion-criterion}, we substitute the above sum condition. Thereby, we transform the left-hand side to  
\begin{align*}
 & -\sum_l \grad\elchempotMixl\cdot\rolfluxrel = -\sum_{l=1}^{L-1} \grad\brac{\elchempotMixl-\elchempotMixl[L]}\cdot\rolfluxrel~. 
\end{align*}
Following \cite{DeGrootMazur-book, Lavenda-book, PrigogineKondepudi-book},  we choose for the drift \massfluxRol{es}~$\rolfluxrel$ of the solutes, i.e., for $l\in\{1,\ldots,L-1\}$, 
the \textbf{\textit{constitutive ansatzes}} 
\begin{align}\label{eq:modelES-diffusion-ansatz}
 \rolfluxrel :=  -\ml\rol\mobl\grad\brac{\elchempotl-\elchempotl[L]} = -\ml\rol\mobl\grad\brac{\elchempotMixl-\elchempotMixl[L]} \qquad\sim\sqbrac{kgm^{-2}s^{-1}}~.
\end{align}
Here, $0\leq\mobl\sim[m/(Ns)]$ are the so-called mobilities, which are connected to the so-called diffusion coefficients~$0\leq\Dl\sim[m^2/s]$ 
in the \textit{Einstein--Smoluchowski relation} by
\begin{align}\label{eq:modelES-diffusion-smoluchowskiEinstein}
 \mobl(\temp) ~=~ \frac{\Dl}{\constBoltz\temp} \qquad\sim\sqbrac{mN^{-1}s^{-1}}~.  
\end{align}
Generally, the mobilities~$\mobl$ describe the capability of the $l$th chemical species to react to a driving force density. 
More precisely, in the above ansatz the induced drift \massfluxRol~$\rolfluxrel$ and its generating body force density~$\ml\rol\grad\brac{\elchempotMixl-\elchempotMixl[L]}\sim[N/m^3]$ are proportional 
to each other, where the constant of proportionality is given by the mobility~$\mobl$, cf. \cite{GreenInRamos-book, Lyklema-book2, Masliyah-book}.
Hence, the mobilities reflect the magnitude of the induced particle movement generated by a driving force, cf. \cite{Masliyah-book, Probstein-book, Russel-book}. 
Note, that in particular  \eqref{eq:modelES-diffusion-ansatz} shows that the drift \massfluxRol{es}~$\rolfluxrel$ of the solutes are generated by 
their electrochemical potentials of mixing~$\elchempotMixl$ and the electrochemical potential of mixing~$\elchempotMixl[L]$ of the 
solvent. Thus, the ansatzes~\eqref{eq:modelES-diffusion-ansatz} account for solute-solvent interactions. Moreover, the constitutive ansatz for the solvent drift \massfluxRol~$\rolfluxrel[L]$ is 
determined by \eqref{eq:modelES-diffusion-ansatz} and the above sum condition, cf. \cref{remark:NonEqThermo-independentDriftFluxes}. 
Next, we insert the definitions $\elchempotMixl = \chempotMixl + \frac{\constCharge\zl}{\ml}\potEL $ of the electrochemical potentials of mixing from \cref{remark:NonEqThermo-electrochemicalPotential} 
into the ansatzes~\eqref{eq:modelES-diffusion-ansatz}. This leads with \eqref{eq:modelES-energy-chemPot} and \ref{Assump:electrostatics} to
\begin{align*}
  \rolfluxrel &= -\ml\rol\mobl \grad\brac{\elchempotMixl-\elchempotMixl[L]} \\[2.0mm]
              &= -\ml\rol[]\yl\mobl \grad\brac{ \frac{\constBoltz\temp}{\ml}\sqbrac{\beta_l +\ln(\yl) } - \frac{\constBoltz\temp}{\ml[L]}\sqbrac{\beta_L + \ln(\yl[L])}  + \frac{\constCharge\zl}{\ml}\potEL - \frac{\constCharge\zl[L]}{\ml}\potEL } \\[2.0mm]
              &= -\ml\rol[]\yl\mobl \grad\brac{ \frac{\constBoltz\temp}{\ml}\sqbrac{\beta_l +\ln(\yl) } - \frac{\constBoltz\temp}{\ml[L]}\sqbrac{\beta_L + \ln(\yl[L])}} + \constCharge\rol\mobl\sqbrac{\zl - \frac{\ml\zl[L]}{\ml[L]}}\fieldEL .
\end{align*}
Furthermore, calculating the remaining derivatives, reveals with \eqref{eq:modelES-diffusion-smoluchowskiEinstein} 
\begin{align*}
  \rolfluxrel &=   -\overbrace{\rol[]\Dl\grad\yl+\frac{\ml\rol\Dl}{\ml[L]\yl[L]}\grad\yl[L]}^{\begin{subarray}{c} \text{mixing} \\ \text{induced diffusion} \end{subarray}} 
                   +\overbrace{\frac{\constCharge\rol\Dl}{\constBoltz\temp} \sqbrac{\zl - \frac{\ml\zl[L]}{\ml[L]}}\fieldEL}^{\begin{subarray}{c} \text{electric} \\ \text{induced diffusion} \end{subarray}}  \\[4.0mm]
              &~~~~-\underbrace{\rol\Dl\sqbrac{\beta_l+\ln(\yl)}\grad\ln(\temp) + \frac{\ml\rol\Dl}{\ml[L]}\sqbrac{\beta_L+\ln(\yl[L])}\grad\ln(\temp)~. }_{\begin{subarray}{c} \text{thermal} \\ \text{induced diffusion} \end{subarray}}
\end{align*}
Finally, the ansatzes~\eqref{eq:modelES-diffusion-ansatz} lead immediately to 
\begin{align*}
 &  -\sum_{l=1}^{L-1} \grad\brac{\elchempotMixl-\elchempotMixl[L]}\cdot\rolfluxrel  \\
 &= \sum_{l=1}^{L-1} \ml\rol\mobl \grad\brac{\elchempotMixl-\elchempotMixl[L]}\cdot\grad\brac{\elchempotMixl-\elchempotMixl[L]} \\
 &= \sum_{l=1}^{L-1} \ml\rol\mobl \abs{\grad\brac{\elchempotMixl-\elchempotMixl[L]}}^2 \geq 0~.
\end{align*}
This proves \eqref{eq:modelES-diffusion-criterion}. Hence, ansatzes~\eqref{eq:modelES-diffusion-ansatz} are thermodynamical consistent.
\begin{remark}[Solute-Solute interactions]
 Instead of \eqref{eq:modelES-diffusion-ansatz}, we can choose the ansatzes
 \begin{align*}
  \rolfluxrel := -\sum_{k=1}^{L-1} \rol\mobl[lk] \grad\brac{\elchempotMixl[k]-\elchempotMixl[L]} \qquad\text{for}~~ l\in\{1,\ldots,L-1\}~.
 \end{align*}
 Here, $\mobl[lk]$ is the mobility of the $l$th chemical species with respect to the forces coming from 
 the $k$th electrochemical potential of mixing~$\elchempotMixl[k]$. Hence, in addition to solute-solvent interactions, 
 these ansatzes account for cross effects between the solutes of the mixture. In particular, these ansatzes are the natural choices 
 for modeling cross diffusion.  
\hfill$\square$
\end{remark}
\begin{remark}[Changing the model]
 The above ansatzes~\eqref{eq:modelES-diffusion-ansatz} reveal, that the constitutive ansatzes for the drift \massfluxRol{es} are 
 determined by the constitutive ansatzes for the chemical potentials. Thus, a crucial starting point for generalizations of the model is 
 to find admissible generalizations for the chemical potentials, cf. \cite{Burger12}.
\hfill$\square$
\end{remark}
\fancyhead[L]{Viscous Stress Tensor}
\section{Constitutive Ansatz for the Viscous Stress Tensor}\label{sec:modelES-viscous}
First of all, we recall that according to \ref{Assump:mixtureStress2}, for the viscous stress tensor must hold
\begin{align*}
 \stressViscl[] = \stressViscl[]^\top~.
\end{align*}
Next, we henceforth assume that the rheology of the mixture is sufficiently well described by considering the mixture as newtonian fluid. Thus, following, e.g., 
\cite{EckGarckeKnabner-book, DeGrootMazur-book, Reiner1945, Rivlin1, Rivlin2} we suppose for the viscous stress tensor~$\stressViscl[]$ the \textbf{\textit{newtonian constitutive ansatz}}
\begin{align}\label{eq:modelES-viscous-ansatz}
\stressViscl[] := \eta \sqbrac{\grad\fieldF + (\grad\fieldF)^\top } + \eta_v(\grad\cdot\fieldF)\mathds{1} \qquad\sim\sqbrac{Jm^{-3}}.
\end{align}
Here, $\eta=\eta(\rol[],\temp)\sim[Ns/m^2]$ is the so-called shear viscosity, as the first term models shear effects. Whereas, $\eta_v=\eta_v(\rol[],\temp)\sim[Ns/m^2]$ is the so-called bulk viscosity, 
since the second term describes volume effects, cf. \eqref{eq:NonEqThermo-mass-specificVolPDE}. For the mixture stress tensor~$\stressTotl[]$, this newtonian ansatz results in
\begin{align*}
\stressTotl[] = -\pressHydrl[]\mathds{1}_n + \eta\sqbrac{\grad\fieldF + (\grad\fieldF)^\top } + \eta_v (\grad\cdot\fieldF)\mathds{1}~.
\end{align*}
Obviously, the newtonian ansatz ensures the symmetry of~$\stressViscl[]$ and $\stressTotl[]$. Moreover, the trace of $\stressViscl[]$ is given by 
 \begin{align*}
 \trace{\stressViscl[]} &= \trace{\eta \sqbrac{\grad\fieldF + (\grad\fieldF)^\top } +\eta_v\,(\grad\cdot\fieldF)\mathds{1}} \\
                        &= \eta\,\grad\cdot\fieldF + \eta\,\grad\cdot\fieldF + n\eta_v\,\grad\cdot\fieldF  =  (2\eta+ n\eta_v)\,\grad\cdot\fieldF.
\end{align*} 
\begin{remark}[Traceless newtonian stress tensor]\label{remark:modelES-tracelessNewtonianStresstensor}
 The preceding equation reveals, that we have to enforce    
 \begin{align*}
  2\eta+n\eta_v=0 \qquad\Equivalent\qquad \eta_v=-\frac{2\eta}{n}~,
 \end{align*}
 to obtain a traceless tensor~$\stressViscl[]$. Alternatively, the viscous mixture stress tensor~$\stressViscl[]$ is traceless in 
 incompressible situations, which are characterized by $\grad\cdot\fieldF\equiv 0$. In both cases, the total mixture pressure~$\pressTotl[]$ 
 coincides with the mixture pressure~$\pressHydrl[]$, cf. \eqref{eq:modelES-pdes-pressures}.
\hfill$\square$
\end{remark}
\begin{remark}[Validity of the newtonian ansatz]
 It is important to note, that the newtonian ansatz for $\stressViscl[]$ remains valid, if the barycentric flow on the considered spatial scales is not affected by the size of the constituents and 
 their molecular interactions. Consequently, the newtonian ansatz for $\stressViscl[]$ restricts both, the size of the constituents, and their molecular interactions. In case these assumptions are violated, 
 the microscopic structure of the mixture influences the barycentric flow. This leads to a viscoelastic rheology of the mixture. In these situations, we have to choose among the various constitutive laws for 
 viscoelastic materials instead, cf. \cite{Larson-Rheology-book, Russel-book, Tadmor-book, Tadros-book}. 
\hfill$\square$
\end{remark}
\par
Concerning the thermodynamic consistency, we obtain from \eqref{eq:modelES-pdes-diss} the sufficient condition that the newtonian ansatz~\eqref{eq:modelES-viscous-ansatz} is in 
accordance with the second law of thermodynamics, if this ansatz leads to 
\begin{align}\label{eq:modelES-viscous-entropyCrit1}
 \frac{1}{\temp}\stressViscl[]:\grad\fieldF \geq 0
 \quad\Equivalent\quad 
  \stressViscl[]:\grad\fieldF \geq 0~.
\end{align}
We now recall some facts from linear algebra, cf. \cite{Knabner-LA-book}. Firstly, arbitrary matrices $A,\setR^{n\times n}$ can be decomposed 
into a symmetric part~$A^s$ and a skew symmetric part~$A^a$, i.e.,
\begin{align*}
 A= A^s + A^a, \qquad\text{with }  A^s=\frac{1}{2}\brac{A+A^\top} ~~\text{and}~~ A^a = \frac{1}{2}\brac{A-A^\top}~.
\end{align*}
Secondly, for arbitrary matrices $A,B\in\setR^{n\times n}$ holds 
\begin{align*}
 A:B:=\trace{A^\top B}, \qquad\text{and } \trace{AB}=0 \text{ for } A \text{ symmetric }, B \text{ skew symmetric}.
\end{align*}
By inserting the newtonian ansatz~\eqref{eq:modelES-viscous-ansatz} into \eqref{eq:modelES-viscous-entropyCrit1}, we transform the left-hand side of this inequality 
with the preceding linear algebra facts to
\begin{align}\label{eq:modelES-viscous-entropyCrit2}
 &   \stressViscl[]:\grad\fieldF 
  ~=~\stressViscl[]:\brac{\grad\fieldF}^s  \nonumber\\
 &=  \frac{\eta}{2}\, \sqbrac{\grad\fieldF + (\grad\fieldF)^\top }:\sqbrac{\grad\fieldF + (\grad\fieldF)^\top } +\frac{\eta_v}{2}\,(\grad\cdot\fieldF)\,\trace{\grad\fieldF+(\grad\fieldF)^\top}\nonumber\\
 &=  \frac{\eta}{2}\, \abs{\grad\fieldF + (\grad\fieldF)^\top }^2 + \eta_v\,\abs{\grad\cdot\fieldF}^2. 
\end{align}
Next, recall the elementary inequalities 
\begin{align*}
 (a+b)^2 \leq 2(a^2+ b^2) 
 \qquad\text{and}\qquad
 (a+b+c)^2 \leq 3(a^2+ b^2 + c^2)
 \quad\text{for } a,b,c\geq0~. 
\end{align*}
With these elementary inequalities, we obtain for arbitrary matrices $A,\setR^{n\times n}$, $n\in\{2,3\}$
\begin{align*}
 \abs{A+A^\top}^2 &= \sum_{ij=1}^n \abs{A_{i,j}+A_{ji}}^2 ~\geq~ \sum_{i=1}^n \abs{2\,A_{ii}}^2 ~=~ 4 \sum_{i=1}^n \abs{A_{ii}}^2 \\
                  &\geq \frac{4}{n} \brac{\sum_{i=1}^n A_{ii}}^2 ~=~ \frac{4}{n} \abs{\trace{A}}^2. 
\end{align*}
Substituting this into \eqref{eq:modelES-viscous-entropyCrit2}, we arrive with $\trace{\grad\fieldF}=\grad\cdot\fieldF$ at the inequality
\begin{align*}
 \stressViscl[]:\grad\fieldF &=     \frac{\eta}{2}\, \abs{\grad\fieldF + (\grad\fieldF)^\top }^2 + \eta_v\,\abs{\grad\cdot\fieldF}^2 \\
                             &\geq  \frac{2\eta}{n}\, \abs{\grad\cdot\fieldF}^2 + \eta_v\,\abs{\grad\cdot\fieldF}^2  = \sqbrac{\frac{2\eta}{n} +  \eta_v}  \abs{\grad\cdot\fieldF}^2. 
\end{align*}
Altogether, we have shown for the newtonian ansatz~\eqref{eq:modelES-viscous-ansatz} the following criterion for thermodynamic consistency:
\begin{align*}
 \stressViscl[]:\grad\fieldF \geq 0~, \qquad\text{ if }\quad \frac{2\eta}{n} +  \eta_v \geq 0~.
\end{align*}
Hence, in particular $\eta_v=-\frac{2}{n}\eta$ from \cref{remark:modelES-tracelessNewtonianStresstensor} leads to a thermodynamic consistent ansatz. 
More precisely, this choice of the bulk viscosity~$\eta_v$ is exactly the borderline case of thermodynamic admissible choices of~$\eta_v$.
\fancyhead[L]{Heat flux}
\section{Constitutive Ansatz for the Heat Flux}\label{sec:modelES-heat}
Regarding the heat flux~$\vecq$, we obtain from the second law of thermodynamics~\eqref{eq:modelES-pdes-diss} that a sufficient condition for thermodynamically admissible constitutive ansatzes is
\begin{align*}
 -\frac{1}{\temp^2}\grad\,\temp\cdot\brac{\vecq+\potEL\currentEL-\sum_l\elchempotl\rolfluxrel} \geq 0 
 \quad\Equivalent\quad 
 -\grad\,\temp\cdot\brac{\vecq+\potEL\currentEL-\sum_l\elchempotl\rolfluxrel} \geq 0~.
\end{align*}
Note, that this condition is due to $\elchempotl=\elchempotMixl+\chempotPurel$, and the observation
\begin{align*}
 \sum_l\frac{\chempotPurel}{\temp}\rolfluxrel = \frac{\energyIntPureAnsatz}{\temp} \sum_l \rolfluxrel = 0
\end{align*}
from \cref{sec:modelES-diffusion} equivalent to
\begin{align*}
  -\grad\,\temp\cdot\brac{\vecq+\potEL\currentEL-\sum_l\elchempotMixl\rolfluxrel} \geq 0~.
\end{align*}
Obviously, this criterion is guaranteed, if we suppose the following \textbf{\textit{extended version of Fourier's law}}
\begin{align}\label{eq:modelES-heat-ansatz}
 \vecq := -\kappa \grad\,\temp - \potEL\currentEL +\sum_l\elchempotMixl\rolfluxrel ~.
\end{align}
Here, $0\leq\kappa\sim[J/K]$ is the heat capacity of the mixture, which might be a tensor valued function $\kappa\tx\in\setR^{n\times n}$. In the tensor valued case, $0\leq\kappa$ is to be understood 
in the sense that $\kappa\tx$ are positive definite tensors.  Moreover, $\kappa$ is assumed to be symmetric tensor valued function due to Onsager reciprocal relations, cf. \cite{DeGrootMazur-book, PrigogineKondepudi-book, Lavenda-book}. 
Finally, the preceding ansatz reveals the cross effects
\begin{align*}
    \underbrace{\vecq}_{\begin{subarray}{c} \text{total} \\ \text{heat flow} \end{subarray}} 
  = -\underbrace{\kappa \grad\,\temp}_{\begin{subarray}{c} \text{temperature} \\ \text{driven} \\ \text{heat flow} \end{subarray}} 
    -\underbrace{\potEL\currentEL}_{\begin{subarray}{c} \text{electric} \\ \text{induced} \\ \text{heat flow} \end{subarray}} 
    +\underbrace{\sum_l\elchempotMixl\rolfluxrel}_{\begin{subarray}{c} \text{mixing} \\ \text{induced} \\ \text{total heat flow} \end{subarray}} ~.
\end{align*}
\fancyhead[L]{Mathematical model}
\section{Mathematical Model for Electrolyte Solutions}\label{sec:modelES-model}
In this section, we present a mathematical model for electrolyte solutions. This mathematical model is the condensed output of \cref{chapter:NonEqThermo}, \cref{sec:modelES-maxwell},
and \cref{sec:modelES-energy} -- \cref{sec:modelES-heat}. More precisely, in \cref{chapter:NonEqThermo} and \cref{sec:modelES-maxwell}, we obtained the general governing equations, and in 
\cref{sec:modelES-energy} -- \cref{sec:modelES-heat}, we presented the involved constitutive ansatzes. As we showed the thermodynamical consistency of these constitutive ansatzes, 
the following mathematical model is a thermodynamically consistent model.
\medskip
\par
More precisely, we now repeat the governing equations from \cref{sec:modelES-pdes}, and we combine these equations with the constitutive ansatzes from \cref{sec:modelES-energy} -- \cref{sec:modelES-heat}.\\[5.0mm]
\begin{subequations}
\textbf{1. Electric potential equation: } We have $\fieldEL=-\grad\potEL$ for the electric field~$\fieldEL$, and the electric potential~$\potEL$ solves   
\begin{align}
  -\grad\cdot(\eps_r\grad\potEL) = \frac{1}{\constEL}~ \chargeEL \qquad\text{with }\quad\chargeEL=\sum_l\frac{\constCharge\zl}{\ml}\rol[]\yl~. \label{eq:modelES-model-poisson}
\end{align}
\textbf{2. Mass conservation equations: } For $l\in\cbrac{1,\ldots,L-1}$, we have with $\rol=\rol[]\yl$
\begin{align}
  & \dert(\rol[]\yl) + \grad\cdot\brac{ \rol[]\yl\fieldF + \rolfluxrel} ~=~\rl~, \label{eq:modelES-model-massBalance}    \\[2.0mm]
  & \dert\rol[] + \grad\cdot\brac{ \rol[] \fieldF } = 0~.                        \label{eq:modelES-model-massBalanceTot}
\end{align}
The solvent concentration~$\rol[L]$ is obtained from these equations by 
\begin{align}\label{eq:modelES-model-solventConc}
 \rol[L] = \rol[] - \sum_{l=1}^{L-1} \rol \quad\Equivalent\quad \yl[L] = 1 - \sum_{l=1}^{L-1} \yl~, 
\end{align}
and the mass production rates~$\rl$ are given according to \eqref{eq:modelES-reactions-ansatz} by
\begin{align}\label{eq:modelES-model-ansatzReactions}
 \rl = \sum_j \ml\slj \sqbrac{\Rfjmal - \Rbjmal}. 
\end{align}
Furthermore, for the drift \massfluxRol{es}~$\rolfluxrel$ the ansatzes~\eqref{eq:modelES-diffusion-ansatz} read as
\begin{flalign}\label{eq:modelES-model-ansatzMassFlux}
 \rolfluxrel = \!-\ml\rol\mobl \grad\brac{\elchempotMixl\!-\elchempotMixl[L]}.  
\end{flalign}
Here, we have for the electrochemical potentials of mixing due to \cref{remark:NonEqThermo-electrochemicalPotential} and \eqref{eq:modelES-energy-chemPot}
\begin{align}\label{eq:modelES-model-ansatzChemPot}
  \elchempotMixl = \frac{\constBoltz\temp}{\ml} \, (\beta_l + \ln(\yl) ) + \frac{\constCharge\zl}{\ml} \potEL, 
  \quad\text{with constants } \beta_l \text{ from \cref{sec:modelES-reactions}}.   
\end{align}
\textbf{3. Momentum conservation equations: } For the barycentric momentum density holds
\begin{align}\label{eq:modelES-model-momBalanceTot}
  \dert(\rol[]\fieldF) + \grad\cdot\brac{\rol[]\fieldF\otimes\fieldF} = -\grad\pressHydrl[] + \grad\cdot\stressViscl[] + \chargeEL\fieldEL ~. 
\end{align}
Moreover, the newtonian ansatz~\eqref{eq:modelES-viscous-ansatz} for the viscous mixture stress tensor is given by
\begin{align}\label{eq:modelES-model-mixtureStress1}
 \stressViscl[] = \eta\sqbrac{\grad\fieldF + (\grad\fieldF)^\top} + \eta_v (\grad\cdot\fieldF)\mathds{1}~. 
\end{align}
\textbf{4. Internal energy evolution: } We have $\energyIntl[]=\energyIntMixl[]+\energyIntPurel[]$ for the internal energy~$\energyIntl[]$, where the internal energy of mixing is 
given due to \ref{Assump:ansatzEnergyMix} by
\begin{align*}
 \energyIntMixl[] = \sum_l \sqbrac{\frac{\constBoltz\temp}{\ml}\, \yl\brac{ \beta_l -1+\ln(\yl)} + \frac{\constBoltz\temp}{\ml} \exp(-\beta_l) }.
\end{align*}
The evolution of the pure substance internal energy density~$\rol[]\energyIntPurel[]$ is subject to 
\begin{flalign}\label{eq:modelES-model-energyIntPurePDE}
  &~~\dert\brac{\rol[]\energyIntPurel[]} + \grad\cdot\brac{ \rol[]\energyIntPurel[]\fieldF + \vecq - \sum_l\elchempotMixl\rolfluxrel + \potEL\currentEL } \nonumber\\
  &~~=-\pressHydrl[]\grad\cdot\fieldF + \stressViscl[]:\grad\fieldF -\sum_l\grad\elchempotMixl\cdot\rolfluxrel - \sum_l\elchempotMixl\rl ~. &
\end{flalign}
Here, the ansatz~\eqref{eq:modelES-heat-ansatz} for the heat flux~$\vecq$ reads as
\begin{align}\label{eq:modelES-model-heat}
  \vecq = -\kappa \grad\,\temp - \potEL\currentEL +\sum_l\elchempotMixl\rolfluxrel~.
\end{align}
\end{subequations}
\par
The preceding set of equations~\eqref{eq:modelES-model-poisson}, \eqref{eq:modelES-model-massBalance}, \eqref{eq:modelES-model-massBalanceTot}, \eqref{eq:modelES-model-momBalanceTot}, 
and \eqref{eq:modelES-model-energyIntPurePDE} is exactly the mathematical model, which we propose for electrolyte solutions. This model is thermodynamically consistent, as we proved in 
\cref{sec:modelES-reactions} -- \cref{sec:modelES-heat} for the involved constitutive laws~\eqref{eq:modelES-model-ansatzReactions}, \eqref{eq:modelES-model-ansatzMassFlux}, 
\eqref{eq:modelES-model-ansatzChemPot}, \eqref{eq:modelES-model-mixtureStress1},  and \eqref{eq:modelES-model-heat}, that they are subject to the second law of thermodynamics~\eqref{eq:modelES-pdes-diss}. 
\begin{remark}[Computation of the model]\label{remark:modelES-modelComputation}
  The preceding model contains the unknowns 
  \begin{align*}
  \brac{\potEL,\yl[1],\ldots,\yl[L-1],\rol[],\fieldF,\pressHydrl[],\temp} \in \setR^{L+3+n}. 
  \end{align*}
  To compute these $L+3+n$~unknowns, we solve the $L+n+2$ equations~\eqref{eq:modelES-model-poisson}, \eqref{eq:modelES-model-massBalance}, \eqref{eq:modelES-model-massBalanceTot}, 
  \eqref{eq:modelES-model-momBalanceTot}, and \eqref{eq:modelES-model-energyIntPurePDE}. Thus, to close the model, we have to apply an additional constitutive law for the mixture pressure~$\pressHydrl[]$.
  For that purpose, e.g., the ideal gas law or \vdW\ equation of state can be used. 
  Moreover, we supposed $\energyIntPurel[]=\energyIntPureAnsatz(\entropyPure,\volumeSpecific)$ in \ref{Assump:ansatzEnergyPureTot}. Thus, to close equation~\eqref{eq:modelES-model-energyIntPurePDE},
  we have to specify the function~$\energyIntPureAnsatz$, e.g., by the simple ansatz  
  \begin{flalign*}
    ~~~\energyIntPureAnsatz(\entropyPure,\volumeSpecific) &:= \energyIntPureAnsatz_1(\entropyPure) + \pressHydrl[]\volumeSpecific  \quad\sim[J\,kg^{-1}], 
     \qquad\text{with}\quad \energyIntPureAnsatz_1(\entropyPure) := \frac{\constBoltz\temp_r}{\ml[a]} \, \exp\brac{\frac{\ml[a]\entropyPure}{\constBoltz}} \quad\sim[J\,kg^{-1}]. &
  \end{flalign*}
  Here, $\temp_r$ is a given reference temperature, $\ml[a]:=\frac{1}{L}\sum_l\ml$ is the average molecular mass, and the specific pure substance entropy~$\entropyPure$ is defined by 
  \begin{align*}
   \entropyPure := \frac{\constBoltz}{\ml[a]} \ln\!\brac{\frac{\temp}{\temp_r}} \quad\sim[J K^{-1} \,kg^{-1}].
  \end{align*}
  These ansatzes lead directly to 
  \begin{align*}
      \energyIntPureAnsatz_1 
    = \frac{\constBoltz\temp_r}{\ml[a]} \, \exp\brac{\frac{\ml[a]\entropyPure}{\constBoltz}} 
    = \frac{\constBoltz\temp_r}{\ml[a]} \, \exp\brac{ \ln\!\brac{\frac{\temp}{\temp_r}} } 
    = \frac{\constBoltz\temp}{\ml[a]} \qquad\sim[J\,kg^{-1}].
  \end{align*}
  Furthermore, from these ansatzes we rediscover the thermodynamic definitions~\eqref{eq:NonEqThermo-entropy-temp}, \eqref{eq:NonEqThermo-entropy-press} of the temperature~$\temp$ and 
  the pressure~$\pressHydrl[]$ via
  \begin{align*}
       \derd[\entropy] \energyIntPureAnsatz_1
    &= \derd[\entropy] \sqbrac{ \frac{\constBoltz\temp_r}{\ml[a]} \, \exp\brac{\frac{\ml[a]\entropyPure}{\constBoltz}} } 
     = \temp_r \, \exp\brac{\frac{\ml[a]\entropyPure}{\constBoltz}} 
     = \temp_r \, \exp\brac{ \ln\!\brac{\frac{\temp}{\temp_r}} } 
     = \temp~,  \\
       \derd[\volumeSpecific] (\pressHydrl[]\volumeSpecific)
    &= \pressHydrl[]~.                                          
  \end{align*}
  Substituting these ansatzes into \eqref{eq:modelES-pdes-energyIntPurePDE}, involving the definition~$\volumeSpecific=\rol[]^{-1}$ from \cref{sec:NonEqThermo-mass}, 
  and replacing the heat flux~$\vecq$ by \eqref{eq:modelES-model-heat}, yields the temperature equation
  \begin{align*}
    & \dert\brac{\rol[]\frac{\constBoltz\temp}{\ml[a]}+\pressHydrl[]} + \grad\cdot\brac{ \rol[]\frac{\constBoltz\temp}{\ml[a]}\fieldF+\pressHydrl[]\fieldF - \kappa\grad\,\temp } \\[2.0mm]
    &=\stressTotl[]:\!\grad\fieldF -\sum_l\grad\elchempotMixl\!\cdot\rolfluxrel - \!\sum_l\elchempotMixl\rl \,. 
  \end{align*}
  However, in the preceding sections, we kept the ansatz for $\energyIntPurel[]$ on the abstract level $\energyIntPurel[]=\energyIntPureAnsatz$, 
  as \cref{sec:modelES-reactions} and \cref{sec:modelES-diffusion} revealed that due to \ref{Assump:ansatzEnergyPure}, the ansatzes for $\energyIntPureAnsatz$ lead to vanishing contributions in the second law of thermodynamics~\eqref{eq:modelES-pdes-diss}.  
\hfill$\square$
\end{remark} 
\fancyhead[L]{\pnp}
\section{The \pnp\ with Convection}\label{sec:modelES-pnp}
Subsequently, we show that the model from \cref{sec:modelES-model} contains the well-known family of \pnp{s}, 
cf. \cite{Allaire10, Burger12, DreyerGuhlkeMuller, HyongEtAll, RayMunteanKnabner, Samohyl, Masliyah-book, Roubicek2005-1}. 
We start this task by imposing the following additional assumptions.
\begin{enumerate}[align=left, leftmargin=*, topsep=2.0mm, itemsep=-1.0mm, label={(PNP\arabic*)}, start=1]
 \item We confine ourselves to isothermal situation, i.e., $\temp\equiv const$.%
       \label{PNP:assump1}
 \item We restrict ourselves to incompressible electrolyte solutions, i.e., the mixture density~$\rol[]$ does not change with varying pressure~$\pressHydrl[]$. 
       This is commonly modeled in terms of $\rol[]\equiv const$, which transforms equation~\eqref{eq:modelES-model-massBalanceTot} to the well-known 
       incompressibility constraint~$\grad\cdot\fieldF=0$, cf. \cite{Madja-book, Temam-book}.%
       \label{PNP:assump2}
 \item We suppose an electrically neutral solvent, i.e., $\zl[L]=0$ and thus $\elchempotMixl[L]=\chempotMixl[L]$.%
       \label{PNP:assump3}
 \item We assume $\rl[L]=0$, which means that the solvent is nonreactive.%
       \label{PNP:assump4}
 \item We limit ourselves to dilute electrolyte solutions. Here, the mass fraction of the solvent~$\yl[L]$ is order of magnitudes above the sum of the solute mass fractions~$\yl$, i.e.,
       $\yl[L] \ll \sum_{l=1}^{L-1}\yl$. Hence, we have for the solvent mass fraction~$\yl[L]$ the expression $\yl[L]\tx=\yl[L]^\ast+\delta\yl[L]\tx$, where $\yl[L]^\ast$ is a given constant value, 
       and $\delta\yl[L]\tx$ captures the small variations. Therefore, together with~\ref{PNP:assump4}, we have $\yl[L]\approx const$.%
       \label{PNP:assump5}
\end{enumerate}
Assumptions~\ref{PNP:assump2} and \ref{PNP:assump5} result together with \eqref{eq:modelES-model-ansatzChemPot} for the electrochemical potentials of mixing in 
\begin{align*}
 \elchempotMixl[L]=\frac{\constBoltz\temp}{\ml} \,(\beta_l + \ln(\yl[L]^\ast+\delta\yl[L]))\approx const 
  \quad\Hence\quad \grad\elchempotMixl[L]\approx 0~.
\end{align*}
Thus, we henceforth neglect the contributions from the electrochemical potentials of mixing~$\elchempotMixl[L]$ of the solvent. This leads for the drift \massfluxRol{es}~$\rolfluxrel$ from \eqref{eq:modelES-model-ansatzMassFlux} 
with \ref{PNP:assump1}, \ref{PNP:assump2} and \eqref{eq:modelES-diffusion-smoluchowskiEinstein} to 
\begin{align*}
 \rolfluxrel &\approx -\ml\rol\mobl \grad\elchempotMixl 
              =       -\ml\rol[]\yl\mobl \grad\brac{\frac{\constBoltz\temp}{\ml} \, (\beta_l + \ln(\yl) ) + \frac{\constCharge\zl}{\ml}\potEL} \\
             &=       -\rol[]\mobl\constBoltz\temp\grad\yl -  \rol[]\yl\mobl\constCharge\zl \,\grad\potEL 
              =       -\Dl\grad\rol + \frac{\constCharge\zl\Dl}{\constBoltz\temp} \,\rol \fieldEL~.
\end{align*}
We note, that mass fluxes~$\rolfluxrel$, which are solely generated by the electrochemical potential of mixing~$\elchempotMixl$ are commonly called \textit{Nernst--Planck fluxes}, 
cf. \cite{Hunter-book, Masliyah-book, Probstein-book, Russel-book}. For this reason, we subsequently refer to the mass conservation equations as the \textit{Nernst--Planck equations}. 
Next, we deduce from \cref{sec:modelES-viscous}, that assumption~\ref{PNP:assump2} leads to
\begin{align*}
 \stressViscl[]=\eta\sqbrac{\grad\fieldF+(\grad\fieldF)^\top} \qquad\text{and}\qquad \trace{\stressViscl[]}=0~.
\end{align*}
Hence, substituting the ansatz for the $\stressViscl[]$ into the barycentric momentum balance equations~\eqref{eq:modelES-model-massBalanceTot}, reduces these equations 
together with \eqref{eq:modelES-model-massBalanceTot} and \eqref{PNP:assump2} to the famous Navier--Stokes equations, cf. \cite{EckGarckeKnabner-book, DeGrootMazur-book, Temam-book, Wilmanski-book}. 
Finally, as we confine ourselves to isothermal situations, we omit the temperature equation~\eqref{eq:modelES-model-energyIntPurePDE}. 
Altogether, the simplified mathematical model is given by the 
following set of equations:
\\[4.0mm]
\begin{subequations}
\textbf{1. Poisson's equation: } We have $\fieldEL=-\grad\potEL$ for the electric field~$\fieldEL$, and the electric potential~$\potEL$ solves   
\begin{align}\label{eq:modelES-pnp-poisson}
  -\grad\cdot(\eps_r\grad\potEL) = \frac{1}{\constEL}~ \chargeEL \qquad\text{with }\quad \chargeEL=\sum_l\frac{\constCharge\zl}{\ml}\rol~. 
\end{align}
\textbf{2. Nernst--Planck equations: } For $l\in\cbrac{1,\ldots,L-1}$, we have 
\begin{align}\label{eq:modelES-pnp-massBalance}
   \dert\rol + \grad\cdot\brac{ \rol\fieldF -\Dl\grad\rol - \frac{\constCharge\Dl\zl}{\constBoltz\temp}\,\rol\fieldEL } ~=~\rl~,     
\end{align}
Here, the ansatzes for the mass production rates~$\rl$ are given by \eqref{eq:modelES-model-ansatzReactions}, and the solvent concentration~$\rol[L]$
is obtained with~\eqref{eq:modelES-model-solventConc}, \ref{PNP:assump2}. This \enquote{postprocessing} calculation of $\rol[L]$ is a good verification of the crucial assumption~\ref{PNP:assump5}.
\\[3.0mm]
\textbf{3. Navier--Stokes equations: } For the barycentric momentum density holds
\begin{align}
  & \grad\cdot\fieldF = 0~,                                                                                                                 \label{eq:modelES-pnp-massBalanceTot}\\
  &\rol[] \dert\fieldF + \rol[]\grad\cdot\brac{\fieldF\otimes\fieldF} = -\grad\pressHydrl[] + 2\eta\Delta\fieldF + \chargeEL\fieldEL ~.   \label{eq:modelES-pnp-momBalanceTot}
\end{align}
\end{subequations}
This system of equations is the so-called \textbf{\textit{\nspnp}},
and in particular for electrolyte solutions at rest, i.e., $\fieldF\equiv 0$, this system is known as the so-called \textbf{\textit{\pnp}}. 
Thus, the \pnp\ captures dilute, incompressible, isothermal, and newtonian electrolyte solutions at rest. Furthermore, the  \pnp\ is the 
common standard model for the investigating the interplay between diffusion processes and electrostatic effects, cf. \cite{Elimelech-book, Juengel-book, Lyklema-book2, Masliyah-book, Probstein-book, Russel-book}.
Note, that the \pnp\ is also known as the drift-diffusion equations. Moreover, in case of $L=3$ and $\zl[1]=1=-\zl[2]$, the \pnp\ reduces to the van-Rosenbrock equations 
resp.~the semiconductor device equations. In particular the semiconductor device equations have been intensively analytically studied and great parts of the analytical theory for the \pnp\ have 
been developed in the context of semiconductor device equations. 
\begin{remark}[Ansatz for Computation]
 The preceding model contains the unknowns
 \begin{align*}
  \brac{\potEL,\rol[1],\ldots,\rol[L-1],\pressHydrl[],\fieldF} \in \setR^{L+1+n}.
 \end{align*}
 To compute these unknowns, we solve the $L+1+n$ equations~\eqref{eq:modelES-pnp-poisson}, \eqref{eq:modelES-pnp-massBalance}, \eqref{eq:modelES-pnp-massBalanceTot}, and \eqref{eq:modelES-pnp-momBalanceTot}.
 Hence, the \nspnp\ is a closed system.
\hfill$\square$ 
\end{remark}
\begin{remark}[\spnp]
In case of a fully developed laminar flow, the rheology is sufficiently well described by the stationary Stokes equations. Hence, it is admissible to replace the Navier--Stokes equations 
by the stationary Stokes equations. This leads us to the so-called \textbf{\textit{\spnp}}, which describes dilute, isothermal, newtonian, and incompressible electrolyte solutions, 
which are restricted to fully developed laminar barycentric flow. These systems are used as so-called pore-scale models%
\footnote{For the notion of pore-scales resp. field-scales, we refer to \cite{bear-book} \label{footnote:fieldscale} } 
for electrolyte solutions in porous media, cf. \cite{RayMunteanKnabner, RayFrankNoorden_pnp, Allaire10}.
\hfill$\square$ 
\end{remark}
\medskip
\par
A field-scale%
\footnote{See footnote~\ref{footnote:fieldscale} }  
model for electrolyte solutions in porous media is the so-called \textbf{\textit{\dpnp}}. This model captures dilute, isothermal, newtonian, and incompressible electrolyte solutions 
in porous media on field scales. This system can be obtained as homogenization limit from pore-scale \spnp{s}, cf. \cite{RayMunteanKnabner, RayFrankNoorden_pnp, Allaire10}.
More precisely, the Darcy-Poison-Nernst--Planck system is given by the following set of equations. 
\\[4.0mm]
\begin{subequations}
\textbf{1. Poisson's equation: } We have $\fieldEL=-\grad\potEL$ for the electric field~$\fieldEL$, and the electric potential~$\potEL$ solves   
\begin{align}\label{eq:modelES-pnp-dpnpPoisson}
  -\grad\cdot(\eps_r\grad\potEL) = \frac{\theta}{\constEL}~ \chargeEL \qquad\text{with }\quad \chargeEL=\sum_l\frac{\constCharge\zl}{\ml}\rol~. 
\end{align}
Here, the porosity~$\theta$ occurs during the homogenization procedure as an additional parameter.\\[3.0mm]
\textbf{2. Nernst--Planck equations: } For $l\in\cbrac{1,\ldots,L-1}$, we have with $\rol=\rol[]\yl$
\begin{align}\label{eq:modelES-pnp-dpnpMassBalance}
   \theta\dert\rol + \grad\cdot\brac{ \rol\fieldF -\Dl\grad\rol - \frac{\constCharge\Dl\zl}{\constBoltz\temp}\,\rol\fieldEL } ~=~\theta\rl~,     
\end{align}
Here, the ansatzes for the mass production rates~$\rl$ are given by \eqref{eq:modelES-model-ansatzReactions}.
\\[3.0mm]
\textbf{3. Extended Darcy's Law: } For the barycentric momentum density holds
\begin{align}
  & \grad\cdot\fieldF = 0~,                                                        \label{eq:modelES-pnp-dpnpMassBalanceTot}\\
  & \fieldF = \permeabH\mu^{-1} \brac{-\grad\pressHydrl[] + \chargeEL\fieldEL}~. \label{eq:modelES-pnp-dpnpMomBalanceTot}
\end{align}
Here, $\permeabH\tx\in\setR^{n\times n}$ is the tensor valued permeability function  and $\mu$ the dynamic viscosity.
\end{subequations}
\medskip
\par
Although the \dpnp\ is a field scale model, it contains an electroosmotic force term in extended Darcy's law~\eqref{eq:modelES-pnp-dpnpMomBalanceTot}. 
This is remarkable, as electroosmotic flows are generated only in very small electric double layers around the solid matrix, cf. \cite{Hunter-book, Masliyah-book, ZholkovskijEtAll}. 
The physical reason why electroosmotic flows nevertheless become visible even on field scales are by far dominating surface effects in porous media. 
\begin{remark}[Reformulation of the \dpnp]\label{remark:modelES-dpnp}
  Note, that by dividing Nernst--Planck equations by $\ml$, we obtain with $\rol=\ml\nl$ from \cref{sec:NonEqThermo-mass}
  \begin{align*}
    \theta\dert\nl + \grad\cdot\brac{ \nl\fieldF -\Dl\grad\nl - \frac{\constCharge\Dl\zl}{\constBoltz\temp}\,\nl\fieldEL } ~=~\frac{\theta}{\ml}\rl~.    
  \end{align*}
  Furthermore, the definition of the mass production rates~$\rl$ and the $l$th total reaction rate~$\Rl^{tot}$ in \cref{sec:modelES-reactions} reveal 
  \begin{align*}
    \frac{\theta}{\ml}\rl = \frac{\theta}{\ml} \ml \sum_j\slj\Rj = \theta\sum_j\slj\Rj = \theta\Rl^{tot}~.    
  \end{align*}
  Thus, we can rewrite the Nernst--Planck equations as 
  \begin{align*}
    \theta\dert\nl + \grad\cdot\brac{ \nl\fieldF -\Dl\grad\nl - \frac{\constCharge\Dl\zl}{\constBoltz\temp}\,\nl\fieldEL } ~=~\theta\Rl^{tot}~.    
  \end{align*}
  Furthermore, from \cref{sec:NonEqThermo-charge}, we recall
  \begin{align*}
  \chargeEL = \sum_l \frac{\constCharge\zl\rol}{\ml} = \sum_l \constCharge\zl\nl~.
  \end{align*}
  Thus, we can equivalently reformulate the \dpnp\ with the number concentrations~$\nl$. 
\hfill$\square$ 
\end{remark}
%
%
%
%
\fancyhead[L]{Conclusion}
\section{Conclusion}
In \cref{chapter:modelES} of this paper, we presented a thermodynamically consistent mathematical model for electrolyte solutions. This model is based on the general governing equations for 
mixtures of charged constituents, which we derived by means of nonequilibrium thermodynamics in \cref{chapter:NonEqThermo}. These equations were shortly summarized in \cref{sec:modelES-pdes}. 
Furthermore, we combined these nonrelativistic equations with the electrostatic limit of Maxwell's equations from \cref{sec:modelES-maxwell}, and we applied several constitutive ansatzes in 
\cref{sec:modelES-energy} -- \cref{sec:modelES-heat}. Thereby, we transformed the general governing equations into a specific physical model for electrolyte solutions. Most importantly, 
we proved for all constitutive laws the thermodynamical consistency, i.e., all constitutive laws respect the second law of thermodynamics~\eqref{eq:modelES-pdes-diss}.
Next, in \cref{sec:modelES-model}, we summarized the resulting mathematical model. Furthermore, by applying suitable simplifying assumptions, we showed in \cref{sec:modelES-pnp}, 
that the well-known and widely used family of \pnp{s} is contained in the model from \cref{sec:modelES-model}. More precisely, the choices of the constitutive ansatzes 
in \cref{sec:modelES-energy} -- \cref{sec:modelES-heat}  were exactly motivated by the goal, to obtain a model, that contains the family of \pnp{s}.
\medskip
\par
In summary, the first contribution of \cref{chapter:modelES} of this paper was to identify in \cref{sec:modelES-maxwell}, in which situations the electric phenomena are sufficiently captured by Poisson's equation. 
Secondly, the main contribution of \cref{chapter:modelES} of this paper was to embed the family of \pnp{s} in the general framework of nonequilibrium thermodynamics. Thereby, we provided a 
thermodynamical verification and we clearly revealed the assumptions and restrictions, which are implicitly contained in \pnp{s}. 
Therefore, we uncovered the limitations of the classical \pnp{s}, and by means of the model from \cref{sec:modelES-model}, 
we additionally presented a possible thermodynamically consistent extension of \pnp{s} to more general situations. 
\medskip
\par
Finally, we note that the presented mathematical model from \cref{sec:modelES-model} is subject to an \enquote{arrow of time}. In \cref{chapter:NonEqThermo}, we already mentioned 
that the second law of thermodynamics is commonly considered to restrict admissible direction of physical processes. More precisely, as the second law of thermodynamic 
states that entropy only can be produced~($\diss\geq0$), we know that irreversible processes~($\diss>0$) never return to their initial states. Illustrative speaking, 
this introduces an \enquote{arrow of time}. To rigorously show this, we restrict ourselves to nonreactive electrolyte solutions and we recall the abstract mass balance 
equations~\eqref{eq:modelES-model-massBalance} 
\begin{align}\label{eq:modelES-conclusion-eq1}
 \dert\rol +\grad\cdot\brac{\rol\fieldF+\rolfluxrel} =0~. \tag{$\ast^1$}
\end{align}
Next, we consider the rescaled functions
\begin{align*}
 \rol\!\brac{\frac{t}{\tau},\frac{x}{l}} \qquad\text{for } (\tau,l)\in\{(1,1), (-1-1)\}~.
\end{align*}
For $(\tau,l)=(1,1)$, we obtain~$\rol^f:=\rol\tx$, which describes the forward-in-time processes, whereas for $(\tau,l)=(-1,-1)$, we obtain~$\rol^b:=\rol(-t,-x)$, which describes the 
backward-in-time processes. Inserting these functions into \eqref{eq:modelES-conclusion-eq1} shows that both, $\rol^f$, and $\rol^b$ solve \eqref{eq:modelES-conclusion-eq1}. Hence, these equations 
are symmetric in time. However, by means of the constitutive ansatzes for $\rolfluxrel$, we transformed above mass balance equations transform to~\eqref{eq:modelES-pnp-dpnpMassBalance}. 
\begin{align}\label{eq:modelES-conclusion-eq2}
 \dert\rol +\grad\cdot\brac{\rol\fieldF-\Dl\grad\rol + \frac{\constCharge\zl\Dl}{\constBoltz\temp}\,\rol\fieldEL} =0~. \tag{$\ast^2$}
\end{align}
Note, that the thermodynamic verification of the constitutive ansatzes for $\rolfluxrel$ in \cref{sec:modelES-diffusion} revealed that these ansatzes lead to production of entropy of mixing. 
Hence, this constitutive ansatz is one source of irreversibility. Moreover, substituting $\rol^f$ and $\rol^b$ into equations~\eqref{eq:modelES-conclusion-eq2}, shows that $\rol^b$ is not 
a solution. Hence, equations~\eqref{eq:modelES-conclusion-eq2} are asymmetric in time. This proves that irreversibility breaks the time symmetry of equations~\eqref{eq:modelES-conclusion-eq1} and introduces 
an \enquote{arrow of time} in equations~\eqref{eq:modelES-conclusion-eq2}. The same analysis holds true for the constitutive laws from \cref{sec:modelES-reactions}, \cref{sec:modelES-viscous}, 
and \cref{sec:modelES-heat}. This reveals the physical meaning of the well-known scaling properties of hyperbolic equations of type~\eqref{eq:modelES-conclusion-eq1} and parabolic equations of 
type~\eqref{eq:modelES-conclusion-eq1}.
%
\section*{Acknowledgements}
M. Herz was supported by the Elite Network of Bavaria. Furthermore, we would like to thank Wolfgang Dreyer for reading an early version of this preprint and for his constructive criticism.
%
\addcontentsline{toc}{part}{References}
\fancyhead[L]{References}
\fancyhead[R]{}
\printbibliography   
\end{document}